\documentclass[12pt]{amsart}
\usepackage{amssymb}

\newcommand{\x}{\times}
\newcommand{\<}{\langle}
\renewcommand{\>}{\rangle}

\renewcommand{\a}{\alpha}

\renewcommand{\d}{\delta}
\newcommand{\D}{\Delta}
\newcommand{\e}{\varepsilon}
\newcommand{\g}{\gamma}

\renewcommand{\l}{\lambda}

\newcommand{\var}{\varphi}
\newcommand{\s}{\sigma}

\renewcommand{\th}{\theta}

\newcommand{\dist}{\mathop{\mathrm{dist}}\nolimits}
\newcommand{\diam}{\mathop{\mathrm{diam}}\nolimits}
\newcommand{\dotcup}{\mathop{{\cup\hspace{-.2cm}\cdot\ }}\nolimits}
\newcommand{\ball}{\mathop{\mathrm{ball}}\nolimits}
\newcommand{\Iso}{\mathop{\mathrm{Iso}}\nolimits}
\newcommand{\Cov}{\mathop{\mathrm{Cov}}\nolimits}
\newcommand{\Fin}{\mathop{\mathrm{Fin}}\nolimits}
\newcommand{\Scv}{\mathop{\mathrm{Scv}}\nolimits}

\begin{document}

\title{Gromov--Hausdorff distance for \\ quantum metric spaces}
\author{Marc A. Rieffel}
\address{Department of Mathematics \\
University of California \\ Berkeley, CA 94720-3840}
\email{rieffel\@math.berkeley.edu}
\date{January 9, 2003}
\thanks{The research reported here was
supported in part by National Science Foundation grant DMS99-70509.}

\begin{abstract}
By a quantum metric space we mean a $C^*$-algebra (or more generally an
order-unit space) equipped with a generalization of the usual 
Lipschitz seminorm on
functions which one
associates to an ordinary metric.  We develop for compact
quantum metric spaces a version of Gromov--Hausdorff distance.  We show that
the basic theorems of the classical theory have natural quantum analogues.
Our main example involves the quantum tori, $A_{\th}$.  We show, for
consistently defined ``metrics'', that if a sequence $\{\th_n\}$ of
parameters converges to a parameter $\th$, then the sequence $\{A_{\th_n}\}$
of quantum tori converges in quantum Gromov--Hausdorff distance to $A_{\th}$.
\end{abstract}

\maketitle

\section{Introduction}

When one looks at the theoretical physics literature which deals with string
theory and related topics, one finds various statements to the effect that
some sequence of operator algebras converges to another operator algebra.  A
mathematician specializing in operator algebras will immediately suspect
that one is dealing here with continuous fields of operator algebras
\cite{D}.  But closer inspection shows that in many situations the
framework of the physicists involves various lengths (see 
references below), and
that the physicists are quite careful about the bookkeeping for these
lengths as they discuss the convergence of their sequences of algebras.
(This is hardly surprising, since the physicists want action functionals,
such as Yang--Mills functionals.)  All of this suggests that there are
metric considerations involved in their convergence of algebras, and that one
is perhaps dealing with some kind of convergence for corresponding
``quantum'' metric spaces.

Within the mathematical literature, the only widely used notion of
convergence of ordinary metric spaces of which I am aware is that
given by the Gromov--Hausdorff distance between metric spaces \cite{G1},
\cite{G2}.  The aim of the present article is to introduce a corresponding
``quantum Gromov--Hausdorff distance'' for ``quantum metric spaces'', and
to develop its basic properties.  As our main example we will consider the
quantum tori \cite{R1} \cite{R3}.  We will see that, 
for a consistent choice of ``metrics'', if a
sequence of parameters $\{\th_n\}$ converges to a parameter $\th$, then the
corresponding sequence of quantum tori, $\{A_{\th_n}\}$, converges in
quantum Gromov--Hausdorff distance to $A_{\th}$.  We remark that recently
quantum tori have found considerable employment in string theory
(\cite{CDS}, \cite{Sw1}, \cite{Sw2}, \cite{K}, \cite{SW} and references
therein).

In this article we will deal only with compact quantum metric spaces.  (Most
of the examples in the string-theory literature are compact.)  I have
already introduced the notion of a compact quantum metric space in
\cite{R4}, \cite{R5} (but without using that terminology), following up on
glimpses of such a notion given by Connes \cite{C1}, \cite{C2} in connection
with his theory of quantum Riemannian geometry defined by Dirac
operators.  The basic
definitions and facts will be reviewed in Section~2.  But, very briefly, in
the quantum case the role of the metric on an ordinary metric space is
played by a generalization of the usual Lipschitz seminorm on functions
which is defined by an ordinary metric.  Our quantum spaces are unital
$C^*$-algebras, or, more generally, order-unit spaces.  It will be crucial
for us that such a ``Lipschitz seminorm'' defines an ordinary metric on the
state-space of our quantum space, in generalization of the Kantorovich
metric \cite{Kn} \cite{KR} 
on the probability measures on an ordinary metric space \cite{R4},
\cite{R5}.

I plan to discuss elsewhere further examples of quantum Gromov--Hausdorff
convergence for situations pertinent to the quantum physics literature.  One
class of examples \cite{R7} 
involves sequences of matrix algebras of increasing
dimension, equipped with consistent ``metrics'', which converge to ordinary
compact metric spaces.  A number of examples in the quantum physics
literature are of this type.
The case in which the ordinary compact space is the $2$-sphere appears in a
number of places, and within this context the matrix algebras are often
referred to as ``fuzzy spheres''.  See, for example \cite{GK1}, \cite{BV},
sections~$2.2$ of \cite{Ta}, and references therein.   
Examples of bookkeeping with lengths is
found, for example, in sections~$7.2$--$3$ of \cite{M}.  (A nice
exposition of some of the relations between string theory and
non-commutative geometry, with much bookkeeping of lengths, can be found in
\cite{K}, but it contains little discussion of convergence of algebras.)
Approximating the sphere by matrix algebras is popular because the symmetry 
group $SU(2)$ acts on the matrix algebras as well as on the sphere, 
whereas traditional ``lattice'' approximations coming from choosing a finite 
number of points on the sphere break that symmetry. In \cite{R7} we use 
the material of the present paper together with ideas from Berezin 
quantization to show how matrix algebras do indeed converge to the sphere 
for quantum Gromov-Hausdorff distance, with analogous results for any
integral coadjoint orbit of any compact Lie group.

The quantum fuzzy
sphere is discussed in \cite{GMS}.
The fuzzy $4$-sphere also makes an occasional appearance --- see \cite{HL}
and references therein.  The case in which the space is the $2$-torus has
also received considerable attention.  See, for example, \cite{Dj},
\cite{AM}, and references therein.  My incomplete calculations make me 
optimistic that in a natural way one can show that matrix algebras converge
to tori (of any  dimension). What will distinguish these matrix algebras
from those converging to the sphere is the ``metric" structure which
is placed on them. This is probably related to the ideas of ``change
of topology'' which one finds in the string-theory literature \cite {MS} 
\cite{Blc} \cite{BBC} \cite{Hyk}.

The case of higher-genus surfaces is
mentioned fleetingly in the seminal paper \cite{CDS}, as well as in
\cite{BM}.  Within the mathematical literature, Berezin--Toeplitz
quantization for compact K\"ahler manifolds has been extensively explored.
(See \cite{Sch} and references therein.)  It leads to sequences of
matrix algebras.  It will be an interesting challenge to see how generally
Berezin--Toeplitz quantization might mesh with our 
quantum metric-space theory.

Very recently, building on the present paper, Hanfeng Li has shown
\cite{Li} for the Connes-Landi-DuboisViolette spheres  $\{S_\th\}$ (and
related quantum manifolds) \cite{CoL} \cite{CoD} with
their Dirac operators, that
they form compact quantum metric spaces, and that if a sequence  $\th_n$ of
parameters converges to a parameter $\th$, then the sequence $S_{\th_n}$
converges to $S_\th$ for quantum Gromov-Hausdorff distance. Also, very
recently David Kerr has developed \cite{Ker} a matricial version
of quantum Gromov-Hausdorff distance.

At an extremely speculative level, there is the popular BFSS conjecture
\cite{BF} in string theory, which conjectures that the putative
``$M$-theory'' which is supposed to unify the various versions of string
theory is a ``suitable'' limit of theories on matrix algebras.  See
\cite{SW} for relations with non-commutative geometry, including bookkeeping
with lengths, especially in Section~7.  One can wonder whether quantum
Gromov--Hausdorff distance might have a bit to say in clarifying
``suitable''.

When one goes back to the string-theory literature to see whether
our notion of quantum Gromov-Hausdorff distance is of use, it is
clear that, while it may clarify matters a bit, it is nevertheless
quite inadequate. The reason is that the string-theorists need
the whole apparatus of bundles, connections, action functionals, etc.
(including in the non-commutative setting). Thus what seems to
be needed is a definition of when (quantum) spaces together with all of
their apparatus are close together. No such definition seems to
have been given so far, even for ordinary spaces, probably
because it is not clear what apparatus to include once one
goes beyond manifolds. (But see \cite{Hw1}, \cite{Hw2}, \cite{GRS}.)
However, hints of what apparatus to include can
be found in the literature on the collapsing of Riemannian manifolds.
This involves bundles, and forms, over the (possibly singular) limit
space. (See \cite{Fu} \cite{G2} \cite{Sa} \cite{Lt1} 
\cite{Lt2} \cite{Sbl} \cite{KSo}
and the references therein, especially to papers of Cheeger, Fukaya
and Gromov.) However, no characterization has been given for the 
limit structures which occur in this 
Riemannian setting, so they do not have an autonomous existence
yet. Thus much less is this so for quantum spaces. It is natural to
speculate that Connes' axioms for non-commutative Riemannian
manifolds \cite{Cn3} \cite{GVF} will play an important role in elucidating
this matter.

The rest of the contents of this paper are as follows.
In Section~3 we recall the definition of ordinary Gromov--Hausdorff
distance, and, motivated by this definition, we develop the material
concerning quotients (= ``subsets'') which we will need in the quantum
setting.  Then in Section~4 we give our definition of quantum
Gromov--Hausdorff distance, and prove that it satisfies the triangle
inequality.  We also discuss there the fact that 
when our definition is applied to ordinary metric spaces,
it does not in general give the ordinary Gromov--Hausdorff distance.  
(But we show that the
quantum distance is never greater than the ordinary distance.) A 
specific example of this failure, found by Hanfeng Li, is presented
in Appendix 1.

In Section~5 we develop a useful technique for estimating quantum
Gromov--Hausdorff distances, involving ``bridges'', and we use it to present
some simple examples.  Bridges are used extensively in later sections.

In the situation of ordinary compact metric spaces, if the classical
Gromov--Hausdorff
distance between two metric spaces is $0$, they need not be the same
set-theoretically, but they will be isometric to each other.  In Section~7
we show that the analogous fact holds for quantum Gromov--Hausdorff
distance.  In preparation for this we discuss isometries between compact
quantum metric spaces in Section~6.

In \cite{R4} several ways of constructing ``metrics'' on $C^*$-algebras were
described when one has an action of a compact group on a $C^*$-algebra
together with a length function on the group.  In Section~8, in preparation
for our discussion of quantum tori, we give an important way of using this
group-action construction to approximate the $C^*$-algebra with metric, 
for Gromov--Hausdorff distance,
by natural finite dimensional compact quantum metric
spaces.  Sections~9--11 then carry out the discussion of quantum tori.  This
involves, in particular, a discussion of continuous fields of quantum metric
spaces.

A basic fact in the classical theory is that the metric space of
isometry-classes of compact metric spaces, equipped with the
Gromov--Hausdorff distance, is complete.  In Section~12 we show that the
same is true in the quantum case.  In the classical case Gromov then gives
an important necessary and sufficient criterion for when a subset of this
complete space is compact.  This criterion has had many applications in
Riemannian geometry \cite{G2}, \cite{Be}, \cite{Sa}, \cite{Sh}.  Section~13
is devoted to formulating and demonstrating the corresponding quantum
criterion. 

Finally, in the brief second appendix we answer 
question $11.1$ of \cite{R5} by
observing that every ``metric'' for a compact quantum metric space can be
obtained by the ``Dirac operator'' construction.

A substantial part of the research reported here was carried out while I
visited the Institut de Math\'ematique de Luminy, Marseille, for three
months.  I would like to thank Gennady Kasparov, Etienne Blanchard, Antony
Wassermann, and Patrick Delorme very much for their warm hospitality and
their mathematical stimulation during my very enjoyable visit.

\section{Compact quantum metric spaces}

As in \cite{R5} we will work with order-unit spaces.  Typically they will
arise as real linear subspaces of the vector space of self-adjoint operators
on a Hilbert space, which contain the identity operator (the order unit).
In fact, any order-unit space can be realized in this way.  We will be most
interested in the order-unit spaces which arise as the space of all
self-adjoint elements of a unital $C^*$-algebra (an algebra of operators on
a Hilbert space which is closed under taking adjoints of operators, and
closed for the operator norm).  But we will see in Sections~8--11 that it is
technically very useful to include order-unit spaces which do not arise in
this way from $C^*$-algebras.

There is an attractive abstract characterization of order-unit spaces due to
Kadison \cite{Kd}, \cite{A}.  An order-unit space is a real partially
ordered vector space, $A$, with a distinguished element $e$ (the order unit)
which satisfies:

\begin{itemize}
\item[1)] (Order unit property)  For each $a \in A$ there is an $r \in
{\mathbb R}$ such that $a \le re$.
\item[2)] (Archimedean property) If $a \in A$ and if $a \le re$ for all $r
\in {\mathbb R}$ with $r > 0$, then $a \le 0$.
\end{itemize}
The norm on an order-unit space is given by
\[
\|a\| = \inf\{r \in {\mathbb R}: -re \le a \le re\}.
\]
Thus $A$ is a normed vector space, and we can consider its dual, $A'$,
consisting of the bounded linear functionals, equipped with the dual norm,
$\|\cdot\|'$.

By a {\em state} of an order-unit space $(A,e)$ we mean a $\mu \in A'$
such that $\mu(e) = 1 = \|\mu\|'$.  States are automatically positive.  We
denote the collection of all the states of $A$, i.e.\ the state-space of
$A$, by $S(A)$.  It is a bounded closed convex subset of $A'$, and so is
compact for the $w^*$-topology on $A'$.  Unless the contrary is specified,
it will always be the $w^*$-topology which we use on $S(A)$.  Each $a \in A$
defines a continuous affine function on $S(A)$ by $a(\mu) = \mu(a)$. It
is a basic theorem of Kadison \cite{Kd} \cite{A} that this representation
of $A$ as affine functions is isometric for the supremum norm
on affine functions.

In \cite{R4}, \cite{R5} we recalled that the metric on an ordinary compact
metric space is determined by the Lipschitz seminorm it defines on
functions.  This suggested that for ``non-commutative spaces'', that is
$C^*$-algebras, the way to specify a ``metric'' is by means of a seminorm
playing the role of a Lipschitz seminorm.  For such a seminorm, $L$, we
defined an ordinary metric, $\rho_L$, on $S(A)$ by
\[
\rho_L(\mu,\nu) = \sup\{|\mu(a) -\nu(a)|: L(a) \le 1\}.
\]
(This may take value $+\infty$ in the absence of further hypotheses.)  This
is a generalization of the Kantorovich metric on the probability measures on
an ordinary compact metric space \cite{Kn} \cite{KR}, 
and within the context of Dirac
operators it was introduced into non-commutative geometry by Connes in
\cite{C1}, \cite{C2}.  In \cite{R5} we extended these ideas to consider
``Lipschitz seminorms'' on order-unit spaces.  Since $C^*$-algebras are over
the complex numbers, while order-unit spaces are over the real numbers, we
should point out that for $C^*$-algebras we require a Lipschitz seminorm,
$L$, to satisfy $L(a^*) = L(a)$ for $a \in A$.  Under this condition it
suffices to take the above supremum just over self-adjoint elements of $A$
when defining $\rho_L$.  One sees this as follows.  (See also lemma~1 of
\cite{IKM}.)  Let $\mu,\nu \in S(A)$ and $\d > 0$ be given.  Then there is
an $a \in A$ such that, after multiplying it by a complex number of
modulus~$1$, we have $\mu(a) - \nu(a) \ge \rho_L(\mu,\nu) - \d$ and $L(a)
\le 1$.  Let $b = (a+a^*)/2$.  Then $L(b) \le 1$ since $L(a^*) = L(a)$,
while still $\mu(b) - \nu(b) \ge \rho_L(\mu,\nu) - \d$.

As in \cite{R5}, we must require that $L$ be such that $\rho_L$ is nice.
This is formulated by:

\bigskip
\noindent
{\bf 2.1 DEFINITION.}  Let $(A,e)$ be an order-unit space.  By a {\em
Lip-norm} on $A$ we mean a seminorm, $L$, on $A$ with the following
properties:

\begin{itemize}
\item[1)] For $a \in A$ we have $L(a) = 0$ iff $a \in {\mathbb R}e$.
\item[2)] The topology on $S(A)$ from the metric $\rho_L$ is the
$w^*$-topology.
\end{itemize}

\bigskip
We remark that condition 2 implies that $\rho_L$ takes only finite
values on $S(A)$. For if $\rho(\mu_0, \nu_0) = +\infty$ for 
some $\mu_0, \nu_0 
\in S(A)$, then $\{\mu: \rho(\mu, \nu_0) < +\infty\}$ is a proper 
subset of $S(A)$ which is both open and closed, which is not possible
since $S(A)$, being convex, is always connected.
In definition $5.1$ of \cite{R5}, where we first gave the definition of a
Lip-norm, we also required that $L$ be lower semi-continuous for the norm on
$A$.  (Even for ordinary metric spaces $L$ is usually not continuous.  That
is why we do not require that $A$ be complete for its norm.)  In theorem
$4.2$ of \cite{R5} we saw that if $L$ is not lower semi-continuous, it can
always be replaced by the largest lower semi-continuous seminorm smaller than
$L$, as this will give the same metric on $S(A)$.  We omit here the
requirement of lower semi-continuity in Definition $2.1$ only for the
convenience that at many intermediate stages of argument we will then not
need to verify, or adjust to obtain, this lower semi-continuity.

In definition $5.1$ of \cite{R5}, condition $2$ of Definition $2.1$ was
stated in a form which appears quite different, but which is shown in
theorem $1.9$ of \cite{R4} to be equivalent, and which is more useful for
verifying examples.  We will recall this other form below in Theorem $4.5$.

We are now prepared to make:

\bigskip
\noindent
{\bf 2.2 DEFINITION.} By a {\em compact quantum metric space} we mean a pair
$(A,L)$ consisting of an order-unit space $A$ with a Lip-norm $L$ defined on
it.

\bigskip
We remark that the use of the word ``quantum'' here may seem a bit of a
stretch.  Given the $C^*$-algebraic origins, at first thought the terms
``non-commutative'' might seem more appropriate.  But since we are using
general order-unit spaces, which have no algebra structure in general, there
is nothing present which can be ``non-commutative''.  On the other hand, the
state space $S(A)$ (with its metric $\rho_L$) will play the central role in
our story, and states play a central role in quantum physics, so it seems to
me that use of the term ``quantum'' is not unreasonable here.

Let $A'$ denote again the Banach-space dual of $A$.  Let $A'{}^{\circ}$
denote the subspace of elements $\l \in A'$ for which $\l(e) = 0$.  As seen
in lemma $2.1$ of \cite{R5}, the ball of radius $2$ about $0$ in
$A'{}^{\circ}$ coincides with $\{\mu-\nu: \mu,\nu \in S(A)\}$.  Define $L'$
on $A'{}^{\circ}$ by the usual formula for a dual norm, namely $L'(\l) =
\sup\{|\l(a)|: L(a) \le 1\}$.  A simple computation, given in lemma $4.3$ of
\cite{R5}, shows that $\rho_L(\mu,\nu) = L'(\mu-\nu)$.  Since $L$ is a
Lip-norm so that $\rho_L$ gives $S(A)$ the $w^*$-topology, it follows
that $\rho_L$ is bounded on $S(A)$, so that $L'$ is bounded on the ball of
radius $2$.  In other words, $L'$ is bounded with respect 
to $\|\cdot\|'$, and
so there is a smallest constant, $r_L$, 
such that $L' \le r_L\|\cdot\|'$.  This
constant is called the {\em radius} of $(A,L)$, because it is exactly half the
diameter of the metric space $(S(A),\rho_L)$, as shown in proposition $2.2$
of \cite{R5}.  There we used the equivalent relation that $\|{\tilde
a}\|^{\sim} \le r_LL^{\sim}({\tilde a})$, where ${\tilde a}$ denotes the
image of $a$ in ${\tilde A} = A/{\mathbb R}e$, 
while $\|\cdot\|^{\sim}$ and $L^{\sim}$
denote the corresponding quotient norms.  We often prefer to use the radius
rather than the diameter so as to avoid factors of $2$ in the above
formulas.  (The above observations
work for any seminorm $L$ for which $L(e) = 0$,
though then we must allow $r_A = +\infty$.)  At times we will denote the
diameter of $(A,L)$ by $\diam(A,L)$.

\section{Quotients (= ``subsets'')}

We now recall the definition of classical Gromov--Hausdorff distance, since
we will model our quantum version on it.  Let $(Z,\rho)$ be an ordinary
compact metric space, and let $X$ be a closed subset of $Z$.  For any $r \in
{\mathbb R}$, $r > 0$, we define the $r$-neighborhood, ${\mathcal
N}_r^{\rho}(X)$, of $X$ for $\rho$ by
\[
{\mathcal N}_r^{\rho}(X) = \{z \in Z:
\mbox{ there is } x \in X \mbox{ with } \rho(z,x) < r\}.
\]
For closed
subsets $X$ and $Y$ of $Z$, the Hausdorff distance, $\dist_H^{\rho}(X,Y)$,
between them for $\rho$, is defined to be:
\[
\dist_H^{\rho}(X,Y) = \inf\{r: X \subseteq {\mathcal N}_r^{\rho}(Y)
\mbox{ and } Y \subseteq {\mathcal N}_r^{\rho}(X)\}.
\]

Suppose now that $(X,\rho_X)$ and $(Y,\rho_Y)$ are independent compact
metric spaces, not viewed as subsets of some larger metric space.  There are
several equivalent definitions of the Gromov--Hausdorff distance between
them \cite{G1}, \cite{G2}, \cite{Be}, \cite{Sa}.  The following is most
convenient for our purposes.  Let $X \dotcup Y$ denote the disjoint union of
$X$ and $Y$, with corresponding compact topology for which $X$ and $Y$ are
closed subspaces.  Let ${\mathcal M}(\rho_X,\rho_Y)$ denote the set of all
metrics on $X \dotcup Y$ giving the topology of $X \dotcup Y$ and whose
restrictions to $X$ and $Y$ are $\rho_X$ and $\rho_Y$, respectively.  Then
the Gromov--Hausdorff distance between $(X,\rho_X)$ and $(Y,\rho_Y)$ is
defined by
\[
\dist_{GH}(X,Y) = \inf\{\dist_H^{\rho}(X,Y): \rho \in {\mathcal
M}(\rho_X,\rho_Y)\}.
\]
Note that for simplicity of notation we are not explicitly indicating the
metrics $\rho_X$ and $\rho_Y$ on the left-hand side.

To formulate a version of the above definition for compact quantum metric
spaces, we need first to see how to translate to the quantum situation the
statement above that the restriction of $\rho$ to $X$ is $\rho_X$.  Thus, let
$(Z,\rho)$ be an arbitrary compact metric space, let $X$ be a closed subset
of $Z$, and let $\rho_X$ be the restriction of $\rho$ to $X$.  Let $C(Z)$
denote the algebra of {\em real}-valued continuous functions on $Z$, and
similarly for $C(X)$.  Let $L_X$ be the Lipschitz ``seminorm'' for $\rho_X$
on $C(X)$, defined for $g \in C(X)$ by
\[
L_X(g) = \sup\{|g(x) - g(y)|/\rho_X(x,y): x \ne y\}.
\]
This supremum can have value $+\infty$.  We can work with this, or with the
dense subalgebra where $L_X$ is finite (the Lipschitz functions), as
convenient.  We will usually do the latter.  Let $L$ be the corresponding
Lipschitz seminorm for $\rho$ on $Z$.  Let $\pi: C(Z) \rightarrow C(X)$
denote the process of restricting functions on $Z$ to $X$.  It is well-known
and easily checked that
\[
L_X(\pi(f)) \le L(f)
\]
for $f \in C(Z)$.  It is also well-known \cite{W2} that for any $g \in C(X)$
there is an $f \in C(Z)$ such that $\pi(f) = g$ and $L(f) = L_X(g)$.  (This
latter {\em fails} for complex-valued functions, see example $1.5.7$ of
\cite{W2}, thus providing one important reason for our emphasis on real
vector spaces.)  All of this says that $L_X$ is exactly the quotient
seminorm from $L$ for $\pi$.

For our purposes, the appropriate morphisms between order-unit spaces are
linear positive maps which preserve the order-units, and consequently are
of norm $1$.  (See proposition II$.1.3$ of \cite{A}, and surrounding
text, for some other possibilities.)  Let $A$ and $B$ be order-unit spaces,
and let $\var: A \rightarrow B$ be a morphism.  We then have the dual
mapping, $\var': B' \rightarrow A'$, which is of norm $1$ and
carries $B'{}^{\circ}$ into
$A'{}^{\circ}$.  For any $\nu \in S(B)$ we have $\var'(\nu) \in S(A)$, and
we obtain in this way a continuous (for the $w^*$-topologies) affine
mapping, $S(\var)$, from $S(B)$ into $S(A)$.  In particular, $S(\var)(S(B))$
will be a closed convex subset of $S(A)$.

For later purposes (e.g., Proposition $12.7$), it is useful to treat
seminorms $L$ on $A$ more general than Lip-norms; we will only require that
$L(a) = 0$ exactly if $a \in {\mathbb R}e$.  As in \cite{R5}, we call these
``Lipschitz seminorms''.  They define a metric, $\rho_L$, on $S(A)$ by the
same formula as before, except that $\rho_L$ can take the value $+\infty$.
In particular, $\rho_L$ need not give the $w^*$-topology, though as seen in
proposition $1.4$ of \cite{R4} the $\rho_L$-topology is always finer than
the $w^*$-topology.

\bigskip
\noindent
{\bf 3.1 PROPOSITION.} {\em Let $A$ and $B$ be order-unit spaces, and let
$\pi: A \rightarrow B$ be a morphism which is surjective, so that $S(\pi)$
is an injection of $S(B)$ into $S(A)$.  Let $L$ be a Lipschitz seminorm on
$A$, and let $L_B$ be the corresponding quotient seminorm on $B$, defined by}
\[
L_B(b) = \inf\{L(a): \pi(a) = b\}.
\]
{\em Then $\pi'$ is an isometry for the norms $L_B'$ and $L'$ 
on $B'{}^{\circ}$ and
$A'{}^{\circ}$, and
$S(\pi)$ is an isometry for the corresponding metrics $\rho_{L_B}$ and
$\rho_L$.  If $L$ is actually a Lip-norm, then so is $L_B$.}

\bigskip
It is important to note that in the setting of Proposition $3.1$ the norm on
$B$ (from the order-unit) will not, in general, be the quotient of that on
$A$, though it is always less than or equal to the quotient norm.  (See
proposition II$.1.6$ of \cite{A} and the example following it.)  Put another
way, in the classical setting described earlier, for any $g \in C(X)$ the
extension theorem says that actually there is an $f \in C(Z)$ restricting to
$g$ such that both $L(f) = L_X(g)$ and $\|f\|_{\infty} = \|g\|_{\infty}$.
But we do not have this stronger statement in our present setting.  However,
this does not seem to cause us difficulties.

We also remark that in the setting of Proposition $3.1$ it is not true in
general that if $L$ is lower semi-continuous then $L_B$ is also, even if $L$
is a Lip-norm:

\bigskip
\noindent
{\bf 3.2 EXAMPLE.} Let $L$ be the Lip-norm for the usual metric on the
interval $Z = [0,3]$.  Let $h$ be the function on $Z$ which interpolates
linearly between the points $(0,0)$, $(1,0)$, $(2,1)$, $(3,-1)$.  Notice
that $L(h) = 2$.  Let $A$ consist of the functions on $Z$ which are the sum
of a polynomial and a scalar multiple of $h$.  We view $L$ as defined just
on $A$.  It is easily seen that the restriction of any Lip-norm to a
subspace containing the order-unit is again a Lip-norm (e.g., use Theorem
$4.5$), which is lower semi-continuous if the original Lip-norm is.  Thus, $L$
on
$A$ is a lower semi-continuous Lip-norm.

Let $B$ be the order-unit space of functions on $[0,2]$ obtained by
restricting the functions in $A$ to $[0,2]$, and let $\pi$ be the
corresponding restriction map.  Let $L_B$ be the quotient of $L$ on $B$.
Then $L_B$ is not lower semi-continuous.  To see this, let $g = \pi(h)$.
Note that $h$ is the only preimage of $g$ under $\pi$ in $A$.  Thus, $L_B(g)
= 2$.  Somewhat as in example $3.5$ of \cite{R5}, for each $n \ge 1$ let
$f_n$ be the continuous function on $Z$ with value $1$ on $[1,3]$, value $0$
on $[0,1-1/n]$, and linear in between.  Let $q_n$ be a polynomial such that
$\|f_n-q_n\|_{\infty} < 1/2n$ on $Z$, and let $p_n$ be the anti-derivative
of
$q_n$ for which $p_n(0) = 0$.  Then
\[
L(p_n) = \|p'_n\|_{\infty} \le 1 + 1/2n.
\]
Furthermore, for $t \in [0,1]$ we have
\[
|p_n(t) - h(t)| = \left| \int_0^t q_n(s)ds\right| \le 1/2n,
\]
while for $t \in [1,2]$ we have
\[
|p_n(t) - h(t)| \le (1/2n) + \left| \int_1^t (q_n(s) - f_n(s))ds\right| \le
1/n.
\]
This says that $\{\pi(p_n)\}$ converges uniformly to $g$.  But
$L_B(\pi(p_n)) \le L(p_n) \le 1 + 1/2n$ while $L_B(g) = 2$.  Thus, $L_B$ is
not lower semi-continuous.

\bigskip
We can avoid the above pathology if we work with closed Lip-norms.  (We
recall from definition $4.5$ of \cite{R5} that a Lip-norm $L$ on $A$ is {\em
closed} if its ``unit ball'' $\{a: L(a) \le 1\}$ is closed in the
completion, ${\bar A}$, of $A$. Closed Lip-norms are automatically
lower semi-continuous.)  But before showing this, we give:

\bigskip
\noindent
{\it\bfseries Proof of Proposition 3.1.} The first part is basically 
just the familiar fact that
when forming quotients of normed spaces, at the level of dual spaces one
obtains isometries.  Since here we work with seminorms, we recall the usual
argument for this, using our present notation.

Let $\l \in B'{}^{\circ}$. 
For any $a \in A$ we clearly have $L_B(\pi(a)) \le
L(a)$, and so if $L(a) \le 1$ we have
\[
|\pi'(\l)(a)| = |\l(\pi(a))| \le L_B'(\l).
\]
Consequently, $L'(\pi'(\l)) \le L_B'(\l)$.  But let $\d > 0$ 
be given, and let $b \in B$ with
$L_B(b) \le 1$.  Then there is an $a \in A$ such that $\pi(a) = b$ and $L(a)
\le \a L(b)$, where $\a = 1 + \d$.  Thus, $L(a/\a) \le L_B(b) \le 1$.
Consequently,
\[
L'(\pi'(\l)) \ge |\pi'(\l)(a/\a)| = |\l(b)|/\a.
\]
Taking the supremum over $b \in B$ with $L_B(b) \le 1$, we see that
\[
L'(\pi'(\l) \ge L_B'(\l)/(1+\d).
\]
Since $\d$ is arbitrary, we see that $\pi'$ is indeed an isometry. But
then $S(\pi)$ is also an isometry, since $\rho_L(\mu, \nu) = L'(\mu - \nu)$
by lemma 4.3 of \cite{R5}, and similarly for $\rho_{L_B}$.

Suppose now that $L$ is a Lip-norm.  Since $\pi'$ is $w^*$-continuous and
injective, and $S(B)$ is compact, $\pi'$ is a homeomorphism of $S(B)$ onto
$\pi'(S(B))$ in $S(A)$.  Because $L$ is a Lip-norm, $\rho_L$ gives the
$w^*$-topology on $S(A)$, and so the restriction of $\rho_L$ to $\pi'(S(B))$
gives the relative topology of $\pi'(S(B))$.  Since $S(\pi)$ is an isometry,
it follows that $\rho_{L_B}$ gives the $w^*$-topology on $S(B)$.  Thus
$L_B$ is a Lip-norm.\qed

\bigskip
For an example of the use of a quotient Lip-norm in a
non-commutative context in the spirit of what we do here,
see equation 3.19 of \cite{ZS}.

Because of Proposition $3.1$ we will find it convenient in a number of
places to identify $S(B)$ with its image in $S(A)$, and view the metric on
$S(B)$ as just the restriction of $\rho_L$ to $S(B)$.  Under the conditions
of Proposition $3.1$ we will say that $L$ {\em
induces} $L_B$.

As promised earlier, we have:

\bigskip
\noindent
{\bf 3.3 PROPOSITION.} {\em Let $\pi: A \rightarrow B$ be a morphism of
order-unit spaces which is surjective.  Let $L$ be a Lip-norm on $A$ which
is closed (and so lower semi-continuous).  Then the quotient Lip-norm,
$L_B$, is closed (so lower semi-continuous).}

\bigskip
\noindent
{\it\bfseries Proof.} Let $(B^c,L_B^c)$ denote the closure of $(B,L_B)$, as
in definition $4.5$ of \cite{R5} (so that the ``unit ball'' for $L_B^c$ is
the closure in ${\bar B}$ of that for $L_B$).  Let $d \in B^c$ with $\|d\|_B
\le 1$.  Then there is a sequence, $\{b_n\}$, in $B$ which converges in
${\bar B}$ to $d$ and for which $\{L_B(b_n)\}$ is bounded, say by $k$.  We
can assume that $\|b_n\| \le \|b\| + k$ for each $n$.  Pick some $\d > 0$.
Then we can find a sequence, $\{a_n\}$, in $A$ such that $\pi(a_n) = b_n$
and $L(a_n) \le L_B(b_n) + \d \le k + \d$ for each $n$.  The following lemma
will be useful later, as well as here.

\bigskip
\noindent
{\bf 3.4 LEMMA.} {\em Let $(A,L)$ be a compact quantum metric space of
radius $r$, and let $\pi: A \rightarrow B$ be a morphism of order-unit
spaces.  For any $a \in A$ we have}
\[
\|a\| \le \|\pi(a)\| + 2rL(a).
\]

\bigskip
\noindent
{\it\bfseries Proof.} Fix some $\nu \in S(B)$.  Then for any $\mu \in S(A)$
we have
\begin{eqnarray*}
|\mu(a)| &\le &|\mu(a) - \pi'(\nu)(a)| + |\pi'(\nu)(a)| \\
&\le &\rho_L(\mu,\pi'(\nu))L(a) + |\nu(\pi(a))| \le 2rL(a) + \|\pi(a)\|.
\end{eqnarray*}
\qed

\bigskip
\noindent
We return to the proof of Proposition $3.3$.  It follows from Lemma $3.4$
that $\{a_n\}$ is a bounded sequence.  Because $L$ is a Lip-norm, any
sequence in $A$ which is bounded for both $\|\cdot\|$ and $L$ is totally
bounded.  (See Theorem $4.5$.)  Thus, $\{a_n\}$ has a subsequence, which we
still denote by $\{a_n\}$, which converges in ${\bar A}$, say to $a_*$.
Since $(A,L)$ is closed, $a_* \in A$.  But $\pi$ is continuous, 
and so $\pi(a_*)
= d$.  Hence $d \in B$.  Thus $B^c = B$, so that $(B,L_B)$ is closed.\qed

\bigskip
But again, we do not insist that our Lip-norms be closed, since that is not
how examples are usually presented to us.

We have seen, in essence, that a surjection $\pi: A \rightarrow B$ of
order-unit spaces gives a continuous affine injection, $\pi'$, of $S(B)$
into $S(A)$, so that $\pi(S(B))$ is a closed convex subset of $S(A)$.  We
now show the converse.  (For the moment, Lip-norms are not involved.)  Let
$A$ be an order-unit space, and let $K$ be a closed convex subset of
$S(A)$.  We denote the space of affine continuous functions on $K$ by
$Af(K)$.  View the elements of $A$ as (affine) functions on $S(A)$, and let
$B$ consist of their restrictions to $K$, with $\pi$ the restriction map of
$A$ into $Af(K)$, and so onto $B$.  As a subspace of $Af(K)$ containing the
order-unit, $B$ has a natural structure as an order-unit space.

\bigskip
\noindent
{\bf 3.5 PROPOSITION.} {\em Let $A$ be an order-unit space, and let $K$ be a
closed convex subset of $S(A)$.  Let $\pi$ be the ``restriction'' map from
$A$ into $Af(K)$, and let $B = \pi(A)$, with its natural order-unit space
structure.  Then $\pi'(S(B)) = K$.}

\bigskip
\noindent
{\it\bfseries Proof.} Since $\pi$ is surjective, $\pi'$ is injective.  By
the definition of $B$, every $\mu \in K$ defines a state, $\mu_B$, of $B$.
For $a \in A$ we then have $\pi'(\mu_B)(a) = \mu_B(\pi(a)) = \mu(a)$.  Thus,
$K \subseteq \pi'(S(B))$.

Suppose now that $\mu \in S(A)$ and $\mu \notin K$.  Because $K$ is convex
and closed, by the Hahn--Banach theorem there is an $a \in A$ and a
$t \in {\mathbb R}$ such that $\mu(a) < t \le \mu_1(a)$ for all $\mu_1 \in
K$.  Upon replacing $a$ by
$a - te$, we obtain $a \in A$ such that $\mu(a) < 0 \le \mu_1(a)$ for all
$\mu_1 \in K$.  Thus, $\pi(a) \ge 0$ in $B$.  If we had $\mu = \pi'(\nu)$
for some $\nu \in S(B)$, we would then have $\mu(a) = \nu(\pi(a)) \ge 0$.
Thus, $\mu \notin \pi'(S(B))$, and so $\pi'(S(B)) = K$ as desired.\qed

\bigskip
We remark that if $K_1$ and $K_2$ are two closed convex subsets of $S(A)$
with $K_1 \nsupseteq K_2$ then the above Hahn--Banach argument shows that
there will be an $a \in A$ which is non-negative on $K_1$ but takes a
strictly negative value on at least one point of $K_2$.  In particular,
suppose that $K$ is a closed convex subset of $S(A)$ which spans all of
$A'$.  Then $\pi: A \rightarrow B$ is bijective, but it will not be an
isomorphism of order-unit spaces if $K \ne S(A)$, because $\pi^{-1}$ will
not be positive.

Given an order-unit space $A$, there is an evident notion of isomorphism
between pairs $(\pi,B)$ such that $B$ is an order-unit space and $\pi$ is a
surjective morphism from $A$ onto $B$.  We will refer to such pairs as
``quotients'' of $A$.  The above comments can be combined to obtain:

\bigskip
\noindent
{\bf 3.6 PROPOSITION.} {\em Let $A$ be an order-unit space.  There is a
natural bijection between isomorphism classes of quotients of $A$ and closed
convex subsets of $S(A)$.}

\bigskip
All of this makes it natural to think of the closed convex subsets of $S(A)$
as corresponding to the ``quantum closed subsets'' of $A$.  If $L$ is a
Lip-norm on
$A$, then $\dist_H^{\rho_L}$ gives a metric on these quantum subsets.  It
is easily seen (and well-known) that the Hausdorff limit of convex subsets is
again convex.  Thus, the set of quantum closed subsets of $A$ is closed in
the space of all closed subsets of $S(A)$ for the Hausdorff metric.  Since
the latter space is compact (proposition $6.1$ of the appendices of
\cite{Sa}) we see that the space of quantum closed subsets of $A$ is a
compact metric space for the Hausdorff metric.

At various points later we will have use for the following transitivity
property, which is really just a (probably known) property of quotients of
general seminorms.

\bigskip
\noindent
{\bf 3.7 PROPOSITION.} {\em Let $A$, $B$ and $C$ be order-unit spaces, let
$\pi_1$ be a morphism of $A$ onto $B$, and let $\pi_2$ be a morphism of $B$
onto $C$.  Set $\pi = \pi_2 \circ \pi_1$, so that $\pi$ is a morphism of $A$
onto $C$.  Let $L$ be a Lipschitz seminorm on $A$, let $L_B$ be its quotient
for $\pi_1$, and let $L_C$ be its quotient for $\pi$.  Let $L_C^B$ be the
quotient of $L_B$ for $\pi_2$.  Then $L_C^B = L_C$.}

\bigskip
\noindent
{\it\bfseries Proof.} Let $c \in C$.  If $a \in A$ and $\pi(a) = c$, then
$\pi_2(\pi_1(a)) = c$, so that $L_C^B(c) \le L_B(\pi_1(a)) \le L(a)$.  It
follows that $L_C^B \le L_C$.  But let $\e > 0$ be given.  Then there exists
$b \in B$ with $\pi_2(b) = c$ and $L_B(b) \le L_C^B(c) + \e$.  But then
there exists $a \in A$ with $\pi_1(a) = b$ and $L(a) \le L_B(b) + \e$.  Then
$\pi(a) = c$ and $L(a) \le L_C^B(c) + 2\e$.  Thus $L_C(c) \le L_C^B(c) +
2\e$.  Since $\e$ is arbitrary, we obtain the desired conclusion.\qed

\bigskip
We will need later a (partial) converse to Proposition $3.1$.  We will
consider two compact quantum metric spaces, $(A,L_A)$ and $(B,L_B)$, and a
surjection $\pi: A \rightarrow B$, and we will consider, in terms of the
dual norms $L'_A$ and $L'_B$, how $L_B$ relates to the quotient of $L_A$.
For this purpose we let $(L_A)_B$ denote the quotient of $L_A$.  Where
convenient we will identify $B'{}^{\circ}$ with its image in $A'{}^{\circ}$
via $\pi'$.  We approach our goal in steps.

\bigskip
\noindent
{\bf 3.8 PROPOSITION.} {\em With notation as above, suppose that $L'_A(\l)
\le L'_B(\l)$ for all $\l \in B'{}^{\circ}$.  If $L_B$ is lower
semi-continuous, then $(L_A)_B \ge L_B$.}

\bigskip
\noindent
{\it\bfseries Proof.} For any $a \in A$ and $\l \in B'{}^{\circ}$ we have
\[
|\<\pi(a),\l\>| = |\<a,\pi'(\l)\>| \le L_A(a)L'_A(\pi'(\l)) \le
L_A(a)L'_B(\l).
\]
Thus, if $L_A(a) \le 1$, then $|\<\pi(a),\l\>| \le L'_B(\l)$, so that
$\pi(a) \in {}^O(\ball_{L'_B})$, where $\ball_{L'_B} = \{\l:
L'_B(\l) \le 1\}$, and where ${}^O$ denotes ``prepolar'' \cite{Cw}.  By the
bipolar theorem \cite{Cw}, ${}^O(\ball_{L'_B}) = (\ball_{L_B})^-$.  Since
$L_B$ is lower semi-continuous, $\ball_{L_B}$ is closed in $B$, and so
$L_B(\pi(a)) \le 1$.  Thus, $L_B(\pi(a)) \le L_A(a)$.  It follows that $L_B
\le (L_A)_B$.\qed

\bigskip
\noindent
{\bf 3.9 PROPOSITION.} {\em With notation as above, suppose that $L'_A(\l)
\ge L'_B(\l)$ for all $\l \in B'{}^{\circ}$.  If $L_A$ is closed, then
$(L_A)_B \le L_B$.}

\bigskip
\noindent
{\it\bfseries Proof.} The hypothesis says that
\[
\ball_{L'_A} \cap B'{}^{\circ} \subseteq \ball_{L'_B}.
\]
Let $\l \in B'{}^{\circ}$.  Now $\l \in \ball_{L'_A} \cap B'{}^{\circ}$ iff
for all
$a
\in A$ with
$L_A(a) \le 1$ we have
\[
|\<\pi(a),\l\>| = |\<a,\pi'(\l)\>| \le 1,
\]
iff $\l \in (\pi(\ball_{L_A}))^O$, the polar \cite{Cw}.  That is,
\[
\ball_{L'_A} \cap B'{}^{\circ} = (\pi(\ball_{L_A}))^O.
\]
Consequently,
\[
{}^O((\pi(\ball_{L_A}))^O) \supseteq {}^O(\ball_{L'_B}),
\]
and so by the bipolar theorem \cite{Cw},
\[
(\pi(\ball_{L_A}))^- \supseteq (\ball_{L_B})^- \supseteq \ball_{L_B}.
\]
Let $b \in \ball_{L_B}$.  It follows that there is a sequence, $\{a_n\}$, in
$\ball_{L_A}$ such that $\|b - \pi(a_n)\| \le 1/n$ for each $n$.  Note that
then $\|\pi(a_n)\| \le \|b\| + 1$ for each $n$.  By Lemma $3.4$ the sequence
$\{a_n\}$ is bounded.  Because $L_A$ is a closed Lip-norm, there is a
subsequence which converges to some $a
\in A$ such that $L_A(a) \le 1$.  (See Theorem $4.5$.)  But $\pi$ is
continuous, so $\pi(a) = b$.  Thus, $(L_A)_B(b) \le 1$.  This shows that
$(L_A)_B \le L_B$.\qed

\bigskip
Upon combining these two propositions, we obtain:

\bigskip
\noindent
{\bf 3.10 COROLLARY.} {\em With notation as above, suppose that $L'_A(\l) =
L'_B(\l)$ for all $\l \in B'{}^{\circ}$.  If $L_A$ is closed and $L_B$ is
lower semi-continuous, then $L_B = (L_A)_B$.}

\bigskip
We remark that it then follows from Proposition $3.3$ that $L_B$ is closed.

\section{Quantum Gromov--Hausdorff distance}

  From our descriptions of the classical Gromov--Hausdorff distance given in
the previous section, it is easy to guess how to proceed in the quantum
case.  Let $(A,L_A)$ and $(B,L_B)$ be compact quantum metric spaces.  The
natural generalization of forming the disjoint union $X \dotcup Y$ is to
form the direct sum, $A \oplus B$, of vector spaces, with $(e_A,e_B)$ as
order-unit, and with evident order structure, to obtain an order unit
space.  We have the evident projections from $A \oplus B$ onto $A$ and $B$,
which are order-unit space morphisms.  The natural generalization of metrics
on $X \dotcup Y$ which restrict to $\rho_X$ and $\rho_Y$ consists of
Lip-norms on $A \oplus B$ which induce $L_A$ and $L_B$.

\bigskip
\noindent
{\bf 4.1 NOTATION.} We will denote by ${\mathcal M}(L_A,L_B)$ the set of
Lip-norms on $A \oplus B$ which induce $L_A$ and $L_B$.

\bigskip
For any given $L \in {\mathcal M}(L_A,L_B)$ we have its metric, $\rho_L$, on
$S(A \oplus B)$.  As discussed somewhat before Proposition $3.1$, we view
$S(A)$ and $S(B)$ as (closed, convex) subsets of $S(A \oplus B)$.  Thus we
can consider the Hausdorff distance between them, that is,
$\dist_H^{\rho_L}(S(A),S(B))$.  We are now ready for the main definition of
this paper.

\bigskip
\noindent
{\bf 4.2 DEFINITION.} Let $(A,L_A)$ and $(B,L_B)$ be compact quantum metric
spaces.  We define the {\em quantum Gromov--Hausdorff distance} between
them, denoted $\dist_q(A,B)$, by
\[
\dist_q(A,B) = \inf\{\dist_H^{\rho_L}(S(A),S(B)): L \in {\mathcal
M}(L_A,L_B)\}.
\]

\bigskip
Let us see in what way $\dist_q$ has the properties of a distance.  It is
clearly symmetric in $A$ and $B$.  We are about to show that it satisfies
the triangle inequality.  In Section~7 we will show that if $\dist_q(A,B) =
0$, then $A$ and $B$ are isometrically isomorphic, just as in the classical
case.  Thus, $\dist_q$ is really a metric on the isometric isomorphism
classes.  We will even show, in Section~12, that the set of isometric
isomorphism classes is complete for this metric.

\bigskip
\noindent
{\bf 4.3 THEOREM (THE TRIANGLE INEQUALITY).} {\em Let $(A,L_A)$, $(B,L_B)$,
and $(C,L_C)$ be compact quantum metric spaces.  Then}
\[
\dist_q(A,C) \le \dist_q(A,B) + \dist_q(B,C).
\]

\bigskip
\noindent
{\it\bfseries Proof.} Let $\e > 0$ be given.  Then we can find an
$L_{AB} \in {\mathcal M}(L_A,L_B)$ such that
\[
\dist_H^{\rho_{L_{AB}}}(S(A),S(B)) \le \dist_q(A,B) + \e.
\]
Similarly we can find a corresponding $L_{BC} \in {\mathcal M}(L_B, L_C)$.  
For the
next lemma, and in many places later, the symbol $\vee$ will mean
``maximum''.

\bigskip
\noindent
{\bf 4.4  LEMMA.} {\em Define $L$ on $A \oplus B \oplus C$ by}
\[
L(a,b,c) = L_{AB}(a,b) \vee L_{BC}(b,c).
\]
{\em Then $L$ is a Lip-norm, and it induces $L_{AB}$, $L_{BC}$, $L_A$,
$L_B$, and $L_C$ for the evident quotient maps.}

\bigskip
\noindent
{\it\bfseries Proof.} We first verify the inducing statements.  Let $\d > 0$
be given.  For any given $(a,b) \in A \oplus B$ we can find $c \in C$ such
that $L_{BC}(b,c) \le L_B(b) + \d$.  But $L_B(b) \le L_{AB}(a,b)$, and so
\[
L(a,b,c) \le L(a,b) + \d.
\]
Thus $L$ induces $L_{AB}$.  In the same way one sees that $L$ induces
$L_{BC}$.  The fact that $L$ then induces $L_A$, $L_B$, and $L_C$ follows
from the transitivity established in Proposition $3.7$.

We now show that $L$ is a Lip-norm.  Condition~1 of Definition $2.1$ is
easily checked.  We must verify Condition~2, which asserts that $\rho_L$
gives the $w^*$-topology on $S(A \oplus B \oplus C)$.  For this purpose
we use the criterion given in theorem $1.9$ of \cite{R4}.  Since we will
also use this criterion later, we recall it here for the reader's
convenience.  For this we use the definition of radius recalled at the end
of Section~2.  As in theorem $1.9$ of \cite{R4}, or section~5 of \cite{R5},
let
\[
{\mathcal B}_1 = \{a: L(a) \le 1 \mbox{ and } \|a\| \le 1\}.
\]
Then we have:

\bigskip
\noindent
{\bf 4.5 THEOREM} (essentially theorem $1.9$ of \cite{R4}).  {\em Let $L$ be
a seminorm on the order-unit space $A$ such that $L(a) = 0$ iff $a \in
{\mathbb R}e$.  Then $\rho_L$ gives $S(A)$ the $w^*$-topology exactly if}

\begin{itemize}
\item[i)] {\em $(A,L)$ has finite radius, and}
\item[ii)] {\em ${\mathcal B}_1$ is totally bounded in $A$ for
$\|\cdot\|_A$.}
\end{itemize}

\bigskip
We now apply this criterion to complete the proof of Lemma $4.4$.  Let
$d_{AB}$ and $d_{BC}$ denote the diameters of $S(A \oplus B)$ and $S(B
\oplus C)$ for $L_{AB}$ and $L_{BC}$, respectively.  Pick any $\nu \in S(B)$.
Then for any $\mu \in S(A)$ and $\zeta \in S(C)$ we have
\[
\rho_L(\mu,\zeta) \le \rho_L(\mu,\nu) + \rho_L(\nu,\zeta) =
\rho_{L_{AB}}(\mu,\nu) + \rho_{L_{BC}}(\nu,\zeta) \le d_{AB} + d_{BC},
\]
since $\rho_L$ restricts to $\rho_{L_{AB}}$ and $\rho_{L_{BC}}$ by
Proposition $3.1$.  Thus, the $\rho_L$-diameter of $S(A) \cup S(B) \cup
S(C)$ is no bigger than $d_{AB} + d_{BC}$.  But $S(A \oplus B \oplus C)$ is
the closed convex hull of $S(A) \cup S(B) \cup S(C)$, and $\rho_L$ is convex
(definition $9.1$ of \cite{R5}) since it just comes from the dual norm $L'$
of $L$ by lemma $4.3$ of \cite{R5}.  It follows easily that the diameter of
$S(A \oplus B \oplus C)$ for $\rho_L$ is no bigger than $d_{AB} + d_{BC}$.
Thus $L$ has finite radius.

To apply Theorem $4.5$ we must show next that
\[
{\mathcal B}_1 = \{(a,b,c): L(a,b,c) \le 1 \mbox{ and } \|(a,b,c)\| \le 1\}
\]
is totally bounded in $A \oplus B \oplus C$.  But if $(a,b,c) \in {\mathcal
B}_1$ then $L_{AB}(a,b) \le 1$ and $\|(a,b)\| \le 1$, while $L_C(c) \le
L_{BC}(c) \le 1$ and $\|c\| \le 1$.  That is,
\[
{\mathcal B}_1 \subseteq {\mathcal B}_1^{AB} \x {\mathcal B}_1^C
\]
for the evident notation.  But ${\mathcal B}_1^{AB}$ and ${\mathcal B}_1^C$
are totally bounded, from which it is easy to see that ${\mathcal B}_1$ is
also.\qed

\bigskip
  From Propositions $3.1$ and $3.7$ we then immediately obtain:

\bigskip
\noindent
{\bf 4.6 LEMMA.} {\em Let $L_{AC}$ be the quotient of the above $L$ for the
evident quotient map from $A \oplus B \oplus C$ onto $A \oplus C$.  Then
$L_{AC}$ is a Lip-norm which induces $L_A$ and $L_C$, that is, $L_{AC} \in
{\mathcal M}(L_A,L_C)$.}

\bigskip
We complete the proof of Theorem $4.3$ by showing that $L_{AC}$ gives the
desired estimate for $\dist_q(A,C)$.  Now $\rho_L$ restricted to $S(A \oplus
B)$ is $\rho_{L_{AB}}$, and so
\[
\dist_H^{\rho_L}(S(A),S(B)) \le \dist_q(A,B) + \e.
\]
In the same way
\[
\dist_H^{\rho_L}(S(B),S(C)) \le \dist_q(B,C) + \e,
\]
and so
\[
\dist_H^{\rho_L}(S(A),S(C)) \le \dist_q(A,B) + \dist_q(B,C) + 2\e.
\]
But by Proposition $3.1$ the restriction of $\rho_L$ to $S(A \oplus C)$ is
$\rho_{L_{AC}}$, and so
\[
\dist_H^{\rho_{L_{AC}}}(S(A),S(C)) \le \dist_q(A,B) + \dist_q(B,C) + 2\e.
\]
Since $L_{AC} \in {\mathcal M}(L_A,L_C)$, it follows that
\[
\dist_q(A,C) \le \dist_q(A,B) + \dist_q(B,C) + 2\e.
\]
Since $\e$ is arbitrary, we obtain the triangle inequality for $\dist_q$.\qed

\bigskip
It is important to emphasize that when our definition of
quantum Gromov-Hausdorff distance is applied 
to ordinary compact metric spaces it does not in general agree
with ordinary Gromov-Hausdorff distance. By this we mean the
following. Let $(X,\rho)$ be an ordinary 
compact metric space, let $A$ denote its
subalgebra of real-valued 
Lipschitz functions, viewed as an order-unit space, and
let $L$ denote the Lipschitz norm on $A$.  Thus $(A,L)$ can be 
considered to be a
compact quantum metric space.  Now suppose that we have two ordinary compact
metric spaces, $(X_j,\rho_j)$ for $j = 1,2$, with associated
$(A_j,L_j)$'s.  Let us compare $\dist_{GH}(X_1,X_2)$ with
$\dist_q(A_1,A_2)$.  Suppose that 
$\rho$ is a metric on $X_1 \dotcup X_2$ whose
restrictions to $X_j$ for
$j = 1,2$ coincide with $\rho_j$.  Let $L$ be the Lip-norm on
$A_1 \oplus A_2$ corresponding to $\rho$, 
with its associated metric $\rho_L$ on $S(A_1 \oplus
A_2)$.  Now $X_j$ is 
naturally identified with the set of extreme points of $S(X_j)$.
Let
$d_{\rho} = \dist_H^{\rho}(X_1,X_2)$.  Given $x \in X_1$ there is a $y \in
X_2$ with $\rho_L(x,y) = \rho(x,y) \le d_{\rho}$, and conversely.  Because
$\rho_L$ is a convex metric (definition $9.1$ of \cite{R5}) it follows that
$\dist_H^{\rho_L}(S(A_1),S(A_2)) \le d_{\rho}$.  Consequently,
$\dist_q(A_1,A_2) \le d_{\rho}$.  Since $\dist_{GH}(X_1,X_2) =
\inf_{\rho}\{d_{\rho}\}$, we obtain:

\bigskip
\noindent
{\bf 4.7 PROPOSITION.} {\em Let $(X_j,\rho_j)$ be ordinary compact metric
spaces for $j = 1,2$, with associated quantum compact metric spaces
$(A_j,L_j)$.  Then}
\[
\dist_q(A_1,A_2) \le \dist_{GH}(X_1,X_2).
\]

\bigskip
But the above inequality can fail
to be an equality. A specific example, found by Hanfeng Li, is
presented in Appendix 1.

Let us explore a bit more what is happening here.  View everything inside
$S(A_1 \oplus A_2)$.  Let $L$ be any Lip-norm on $A_1 \oplus A_2$ which
induces $L_1$ and $L_2$ on $A_1$ and $A_2$, respectively; and let $\rho_L$
be the corresponding metric on $S(A_1 \oplus A_2)$.  Then $\rho_L$ certainly
restricts to a metric on $X_1 \dotcup X_2$ which agrees with $\rho_j$ on
$X_j$.  Thus, $\dist_H^{\rho_L}(X_1,X_2)$ makes sense, and the infimum of
these Hausdorff distances
over all $L$ will be $\dist_{GH}(X_1,X_2)$.  But here we are asking
how close each extreme point of $S(A_1)$ is to an extreme point of $S(A_2)$,
and conversely. However, as the example in Appendix 1 shows, it can happen
that each extreme point of $X_1$ is close to some point of $S(A_2)$,
but that some extreme points of $X_1$ are relatively far from all 
of the extreme
points of $S(A_2)$. And similarly for $X_2$. The consequence will
be that $\dist_q(A_1,A_2) < \dist_{GH}(X_1,X_2)$.

Let us make clear that this can only happen because we admit Lip-norms
which need not come from metrics on $X_1 \dotcup X_2$. To be specific:

\bigskip
\noindent
{\bf 4.8 PROPOSITION.} {\em Let $(X_j,\rho_j)$ for $j = 1,2$ be ordinary
compact metric spaces, and let $(A_j,L_j)$ be their associated quantum
compact metric spaces. Let $L$ be a Lip-norm on $A_1 \oplus A_2$
which comes from a metric, $\rho$, on $X_1 \dotcup X_2$ whose 
restriction to $X_j$ is $\rho_j$ for $j=1, 2$. Then
\[
\dist_{GH}(X_1,X_2) \le \dist_H^{\rho_L}(S(A_1), S(A_2)).
\]}

\bigskip
\noindent
{\it\bfseries Proof.} Let $d = \dist_{GH}(X_1,X_2)$. Then there must
be at least one point, $p$, in one of $X_1$ or $X_2$, say $X_1$, 
such that $\rho( p, y) \ge d$ for all $y \in X_2$. Define $h$ on
$X_1 \dotcup X_2$ by $h(z) = \rho(p, z)$. Then $L(h) = 1$. (This is
the step which fails if $L$ does not come from a metric, as can
be seen by examining example 7.1 of \cite{R5}.) Then for
any $\mu \in S(A_2)$ we have
\[
\rho_L(\d_p, \mu) \ge |(\d_p - \mu)(h)| = \mu(h) 
= \int \rho(p, y)d\mu(y) \ge d.
\]
Thus $\dist_H^{\rho_L}(S(A_1), S(A_2)) \ge d$. \qed

\bigskip
Let us put this another way.  For any order-unit space $A$, let $S^e(A)$
denote the set of extreme points of $S(A)$.  This can be a quite strange
set.  In particular, often it is not closed.  Nevertheless, given quantum
compact metric spaces $(A_j,L_j)$ for $j = 1,2$, and a Lip-norm $L$ on $A_1
\oplus A_2$ inducing $L_1$ and $L_2$, it makes sense to consider
$\dist_H^{\rho_L}(S^e(A_1),S^e(A_2))$.  We will get the same result if we
use the closures, $(S^e(A_j))^-$.  Because $L_{\rho}$ is a convex metric, it
follows just as in the proof of Proposition $4.7$ that:

\bigskip
\noindent
{\bf 4.9 PROPOSITION.} {\em With notation as above,}
\[
\dist_H^{\rho_L}(S(A),S(B)) \le \dist_H^{\rho_L}(S^e(A),S^e(B)).
\]

\bigskip
Thus for arbitrary compact quantum metric spaces $(A_j,L_j)$ we could
define $\dist_q^e$ by
\[
\dist_q^e(A_1,A_2) = \inf_L(\dist_H^{\rho_L}(S^e(A_1),S^e(A_2)).
\]
We would then have a definition which clearly agrees with the classical
definition for ordinary compact metric spaces.  But for many $C^*$-algebras
the set $S^e(A)$ is quite elusive.  
This is certainly true for the quantum
tori (which constitute our main example), when there is 
some irrationality in the
structure constants.  (Though Ed Effros has reminded me
that a theorem of Glimm (see lemma 11.2.1 of \cite{D}) says that for
any unital simple infinite-dimensional $C^*$-algebra, such as the
irrational rotation algebras, $S^e(A)$ is dense in $S(A)$.) 
Thus, it is not clear to me whether $\dist_q^e$ can be
useful.

Only experience with more examples will reveal whether the fact that
$\dist_q$ does not always agree with $\dist_{GH}$ for
ordinary compact metric spaces is ``a feature
or a bug''.

\section{Bridges}

Before continuing with the general theory it seems appropriate to give a
few simple examples.  When dealing with specific examples, the challenge, of
course, is to construct Lip-norms on $A \oplus B$ which induce $L_A$ and
$L_B$ and for which $\dist_H^{\rho_L}(S(A),S(B))$ is
appropriately small. In this section we will formulate a somewhat general
approach to doing this.  We then use it to discuss some simple examples.
But we will find this approach very useful later in connection with our main
examples.

If $L \in {\mathcal M}(L_A,L_B)$, then for $(a,b) \in A \oplus B$ we must
have $L(a,b) \ge L_A(a) \vee L_B(b)$.  Thus, we can look for $L$ of the form
\[
L(a,b) = L_A(a) \vee L_B(b) \vee N(a,b)
\]
for some seminorm $N$ on $A \oplus B$.  Every $L$ is trivially of this form
by setting $N = L$.  But if we think of $N$ as only seeing distances between
$A$ and $B$, not within $A$ or within $B$, then $N$ can have a technically
attractive feature.  Intuitively, because we are viewing $A$ and $B$ as
``disjoint'', for any given $L$ there should be a gap between $S(A)$ and
$S(B)$, that is, a strictly positive lower bound on distances between points
of $S(A)$ and points of $S(B)$.  This means that $N$ should be bounded with
respect to the norm on $A \oplus B$, unlike $L_A$ and $L_B$.  Note that we
must have $N(e_A,e_B) = 0$ if $L$ is to have the same (required) property.
But note also that we cannot have $N(e_A,0) = 0$, for otherwise we would
have $L(e_A,0) = 0$, which is not permitted as it gives infinite distances.
(Equivalently, we can require $N(0,e_B) \ne 0$.)  Of course, the most
important requirement --- the one which is most difficult to arrange for
specific examples --- is that $N$ be such that $L$ induces $L_A$ and $L_B$.
These observations are summarized in:

\bigskip
\noindent
{\bf 5.1 DEFINITION.} Let $(A,L_A)$ and $(B,L_B)$ be compact quantum metric
spaces.  By a {\em bridge} between $(A,L_A)$ and $(B,L_B)$ we mean a
seminorm, $N$, on $A \oplus B$ such that:

\begin{itemize}
\item[1)] $N$ is continuous for the norm on $A \oplus B$.
\item[2)] $N(e_A,e_B) = 0$ but $N(e_A,0) \ne 0$.
\item[3)] For any $a \in A$ and $\d > 0$ there is a $b \in B$ such that
\[
L_B(b) \vee N(a,b) \le L_A(a) + \d,
\]
and similarly for $A$ and $B$ interchanged.
\end{itemize}

\bigskip
\noindent
{\bf 5.2 THEOREM.} {\em Let $N$ be a bridge between compact quantum metric
spaces $(A,L_A)$ and $(B,L_B)$.  Define $L$ on $A \oplus B$ by}
\[
L(a,b) = L_A(a) \vee L_B(b) \vee N(a,b).
\]
{\em Then $L$ is a Lip-norm which induces $L_A$ and $L_B$.  That is, $L \in
{\mathcal M}(L_A,L_B)$.  If $L_A$ and $L_B$ are lower semi-continuous,
then so is $L$.}

\bigskip
\noindent
{\it\bfseries Proof.} From Condition~2 of Definition $5.1$ it is clear that
$L(e_A,e_B) = 0$ and that if $L(a,b) = 0$ then $(a,b) \in {\mathbb
R}(e_A,e_B)$.  From Condition~1 of Definition $5.1$ it is easily seen that
$L$ is lower semi-continuous if $L_A$ and $L_B$ are.  From Condition~3 of
Definition $5.1$ it follows immediately that $L$ induces $L_A$ and $L_B$.
Thus, the one property which is not quickly evident is that $\rho_L$ gives
the $w^*$-topology on $S(A \oplus B)$.  We now verify this property.  For
this purpose we use the criterion given in Theorem $4.5$.  In contrast to
what usually happens, it is the finite radius condition which is slightly
tricky to verify here, while the total boundedness is easy.  We deal first
with the finite radius condition.

Let $\g = 1/N(e_A,0)$.  We call $\g$ the {\em gap} of $N$ because we have:

\bigskip
\noindent
{\bf 5.3 PROPOSITION.} {\em For all $\mu \in S(A)$ and $\nu \in S(B)$ we
have $\rho_L(\mu,\nu) \ge \g$.}

\bigskip
\noindent
{\it\bfseries Proof.} Since $L((\g e_A,0)) = 1$, we have
\[
\rho_L(\mu,\nu) \ge |\mu(\g e_A,0) - \nu(\g e_A,0)| = \g\mu(e_A) = \g.
\]
\qed

\bigskip
Let $(a,b) \in A \oplus B$.  As usual (see section~1 of \cite{R5}) we let
$\max(a) = \inf\{r: a \le re_A\}$, and similarly for $\min(a)$.  Set $m_a =
(\max(a) + \min(a))/2$, viewed as either a scalar or as $m_ae_A \in A$.
Then
\[
\|a\|^{\sim} = \|a - m_a\| = (\max(a) - \min(a))/2.
\]
Define $m_b$ similarly.  For any $t \in {\mathbb R}$ we have
\begin{eqnarray*}
\|(a,b) - t(e_A,e_B)\| &= &\|(a-m_a,b-m_b) + (m_a-te_A,m_b-te_B)\| \\
&\le
&\|a-m_a\| \vee \|b-m_b\| + |m_a - t| \vee |m_b - t|.
\end{eqnarray*}
We minimize the right-most summand by setting $t = (m_a + m_b)/2$.  Then the
above expression is
\[
\le \|a\|^{\sim} \vee \|b\|^{\sim} + |(m_a-m_b)/2|.
\]
Consequently,
\[
\|(a,b)\|^{\sim} \le \|a\|^{\sim} \vee \|b\|^{\sim} + |m_a-m_b|/2.
\]
But
\begin{eqnarray*}
|m_a-m_b| &= &\g N((m_a-m_b)e_A,0) = \g N(m_a,m_b) \\
&\le &\g (N(a,b) +
N(a-m_a,b-m_b)) \\
&\le &\g (L(a,b) + \|N\| \|(a-m_a,b-m_b)\|) \\
&= &\g (L(a,b) + \|N\|(\|a\|^{\sim} \vee \|b\|^{\sim})).
\end{eqnarray*}
Putting this together with the previous inequality, we obtain:
\[
\|(a,b)\|^{\sim} \le (1+\g\|N\|/2)(\|a\|^{\sim} \vee \|b\|^{\sim}) +
(\g/2)L(a,b).
\]
Let $r_A$ and $r_B$ denote the radii of $A$ and $B$, so that $\|a\|^{\sim}
\le r_AL_A(a) \le r_AL(a,b)$, and similarly for $b$.  From this and the
above inequality we obtain
\[
\|(a,b)\|^{\sim} \le ((1+\g\|N\|/2)r_A \vee r_B + \g/2)L(a,b).
\]
This says exactly that
\[
\mbox{radius}(A \oplus B,L) \le (1+\g\|N\|/2)r_A \vee r_B + \g/2,
\]
and so $(A \oplus B,L)$ has finite radius.

We now verify that Condition ii) of Theorem $4.5$ holds in the present
situation.  We denote by ${\mathcal B}_1^A$ the ${\mathcal B}_1$ of Theorem
$4.5$ for our $A$, and similarly for ${\mathcal B}_1^B$.  Let $(a,b) \in
{\mathcal B}_1$ (for $(A \oplus B,L)$).  We know that $L_A(a) \le L(a,b) =
1$ and $\|a\| \le \|(a,b)\| \le 1$, so that $a \in {\mathcal B}_1^A$.  In
the same way $b \in {\mathcal B}_1^B$.  That is, ${\mathcal B}_1 \subseteq
{\mathcal B}_1^A \x {\mathcal B}_1^B$.  But ${\mathcal B}_1^A$ and
${\mathcal B}_1^B$ are totally bounded since $L_A$ and $L_B$ are Lip-norms.
It follows easily that ${\mathcal B}_1$ is totally bounded.\qed

\bigskip
We now give four simple applications of bridges.  The method for
constructing bridges which we use in the first application will play an
important role later.  (See Proposition $11.1$.)  This first application
shows that the distance between any two compact quantum metric spaces is
always finite.

\bigskip
\noindent
{\bf 5.4 PROPOSITION.} {\em Let $(A,L_A)$ and $(B,L_B)$ be compact quantum
metric spaces.  Then}
\[
\dist_q(A,B) \le \mbox{diameter}(A) + \mbox{diameter}(B).
\]

\bigskip
\noindent
{\it\bfseries Proof.} We can think of a proof of the corresponding fact for
ordinary metric spaces as follows.  The two spaces are islands, with
distance being (dry) minimal walking distance.  We can choose a point on
each of these islands and build a bridge between these two points, of
arbitrarily small strictly positive length (the ``gap'').  This will not
change the minimal walking distance on each island individually.

Accordingly, given compact quantum metric spaces $(A,L_A)$ and 
\linebreak
$(B,L_B)$,
we choose arbitrarily $\mu_0 \in S(A)$, $\nu_0 \in S(B)$, and $\g > 0$.  As
bridge we set
\[
N(a,b) = \g^{-1}|\mu_0(a) - \nu_0(b)|.
\]
It is evident that $N$ satisfies the first two conditions of Definition
$5.1$.  To check the third condition, given $a \in A$ it suffices to chose
$b$ to be $\mu_0(a)e_B$; and similarly if we are given $b \in B$.

We now find a bound for the corresponding distance between $S(A)$ and
$S(B)$.  Let $\mu \in S(A)$ and $\nu \in S(B)$.  For any $(a,b) \in A
\oplus B$ with $L(a,b) \le 1$ we have $|\mu_0(a) - \nu_0(b)| \le \g$, so
that
\begin{eqnarray*}
|\mu(a,b) - \nu(a,b)| &= &|\mu(a) - \nu(b)| \\
&\le &|\mu(a) - \mu_0(a)| +
|\mu_0(a) - \nu_0(b)| + |\nu_0(b) -\nu(b)| \\
&\le &\rho_{L_A}(\mu,\mu_0) + \g
+
\rho_{L_B}(\nu_0,\nu).
\end{eqnarray*}
Since $\g$ is arbitrarily small, it follows that
\[
\dist_q^L(A,B) \le \mbox{diameter}(A) + \mbox{diameter}(B),
\]
where $\mbox{diameter}(A)$ is just the usual diameter of the
ordinary metric space
$(S(A),\rho_{L_A})$, or
equivalently, twice the radius defined earlier.\qed

\bigskip
We consider next distances from the one-point space.  That is, we take $B =
{\mathbb R}$ and $L_B \equiv 0$.  Just as in the classical case we have:

\bigskip
\noindent
{\bf 5.5 PROPOSITION.}  {\em Let $B$ be the one-point order-unit space
${\mathbb R}$.  Then for any compact quantum metric space $(A,L_A)$ we have}
\[
\dist_q(A,{\mathbb R}) = \mbox{radius}(A).
\]

\bigskip
\noindent
{\it\bfseries Proof.} Of course, $S(B)$ consists of just one point.  Let $r
= \mbox{radius}(A)$.  Since we saw that $r$ is the ordinary radius of
$S(A)$ for $\rho_{L_A}$, and since $\dist_q$ is defined in terms of
ordinary Hausdorff distance between $S(A)$ and $S(B)$, it is easily seen
that $\dist_q(A,{\mathbb R})$ cannot be strictly smaller than $r$.  But
define a bridge by
\[
N(a,b) = r^{-1}\sup\{|\mu(a) - b|: \mu \in S(A)\}.
\]
It is easily seen that $N$ satisfies the first two conditions of Definition
$5.1$.  We check the third.  Given $b \in B = {\mathbb R}$, it suffices to
choose $a = be_A$.  Suppose instead that we are given $a \in A$.  
Define $m_a$ as done shortly after the proof of Proposition 5.3,
so that, as seen there, $\|a\|^\sim = \|a - m_a\|$. Set $b = m_a$. For
any $\mu \in S(A)$ we have
$$
|\mu(a) - b| = |\mu(a - m_a)| \le \|a\|^\sim   .
$$
Thus $N(a,b) \le r^{-1}\|a\|~ \le L_A(a)$,
as needed.

For the corresponding $L$ we estimate the Hausdorff distance.  Let $\eta$
denote the unique element of $S(B)$.  Then for $(a,b) \in A \oplus B$ with
$L(a,b) \le 1$ we have, for any $\mu \in S(A)$,
\[
|\mu(a,b) - \eta(a,b)| = |\mu(a) - b| \le r.
\]
\qed

\bigskip
\noindent
{\bf 5.6 EXAMPLE.} Given a compact quantum metric space $(A,L)$, it is at
times convenient to form a space consisting of two copies of it at Hausdorff
distance some given $\e > 0$ from each other.  (A typical application will
be given in the next proposition.)  To do this we construct a bridge, $N$, on
$A \oplus A$ by
\[
N(a,b) = \e^{-1}\|a - b\|.
\]
It is easily seen that this is indeed a bridge.  To see that it gives the
desired distance, let $\mu \in S(A)$ viewed as in the first copy, then
choose $\nu = \mu$ viewed as in the second copy.  If $L$ is the Lip-norm on
$A \oplus A$ using $N$, then when $L(a,b) \le 1$ we have $\|a-b\| \le \e$,
and so
\[
|\mu(a,b) - \nu(a,b)| = |\mu(a) - \nu(b)| = |\mu(a-b)| \le \|a-b\| \le \e.
\]
The roles of $\mu$ and $\nu$ can be reversed.  Thus, indeed the two copies
are at distance $\le \e$ from each other.  But the gap $\g$ is clearly $\e$,
and so by Proposition $5.3$ the two copies are at distance exactly $\e$.

\bigskip
We will later have use for the following fact about quotients, as defined in
Section~3, especially as viewed in Proposition $3.5$ and the discussion
which followed.

\bigskip
\noindent
{\bf 5.7 PROPOSITION.} {\em Let $(A,L_A)$ be a compact quantum metric space,
and let $K_1$ and $K_2$ be compact convex subsets of $S(A)$.  Let
$(B_j,L_j)$ for $j = 1,2$ be the corresponding quotients.  Then}
\[
\dist_q(B_1,B_2) \le \dist_H^{\rho_{L_A}}(K_1,K_2).
\]

\bigskip
\noindent
{\it\bfseries Proof.} We remark first that the inequality can be strict, for
example for one-point subsets.  We also remark that $K_1$ and $K_2$ need not
be disjoint here.

We need to construct suitable Lip-norms on $B_1 \oplus B_2$.  Actually, we
will construct a Lip-norm for each $\e > 0$, using the construction of
Example $5.6$.  Let $\e > 0$ be given, and define $N$ on $A \oplus A$ as in
that example.  View $K_j$ as $K_j \x \{0\}$ in $S(A \oplus A)$, and set
$K'_j = \{0\} \x K_j$ in $S(A \oplus A)$, for $j = 1,2$.  
Let $K = \mbox{co}(K_1
\cup K'_2)$, and let $\pi$ be the restriction map from $A \oplus A$ into
$Af(K)$ as discussed before Proposition $3.5$.  Let $C$ be the image of $A
\oplus A$ under $\pi$, so that it is the corresponding quotient.  By
Proposition $3.5$ we can canonically identify $S(C)$ with $K$.  Now $K_1$
and $K'_2$ are disjoint subsets of $K$ (in fact, ``split faces'' as
discussed in \cite{A}), and it is evident that the kernel of $\pi$ is just
the direct sum of the kernels of the restriction maps of $A$ onto $K_1$ and
$K'_2$, respectively.  In other words, there is a canonical identification
$C = B_1 \oplus B_2$.

Let $M$ denote the Lip-norm on $A \oplus A$ for $N$, and let $L$ denote its
quotient on $C$.  According to Proposition $3.1$ the metric $\rho_L$ on
$S(C) = K$ coincides with the restriction to $K$ of $\rho_M$.

Of course, $S(B_1)$ and $S(B_2)$ are naturally included in $S(C)$, and
coincide with the inclusions of $K_1$ and $K'_2$ into $K$.  Thus,
$\dist_q(B_1,B_2) \le \dist_H^{\rho_L}(K_1,K'_2)$.  Now from the
justification given in Example $5.6$ it is clear that
$\dist_H^{\rho_M}(K_2,K'_2) \le \e$.  But $\rho_M$ coincides on $S(A)$ with
$\rho_{L_A}$, 
so $\dist_H^{\rho_M}(K_1,K_2) = \dist_H^{\rho_{L_A}}(K_1,K_2)$.  By
the ordinary triangle inequality
\[
\dist_H^{\rho_M}(K_1,K'_2) \le \dist_H^{\rho_{L_A}}(K_1,K_2) + \e.
\]
But $\dist_H^{\rho_M}(K_1,K'_2) = \dist_H^{\rho_L}(K_1,K'_2)$.  Thus,
\[
\dist_H^{\rho_L}(S(B_1),S(B_2)) \le \dist_H^{\rho_{L_A}}(K_1,K_2) + \e.
\]
Since $\e$ is arbitrary, we obtain the desired result.\qed

\section{Isometries}

One of the basic facts about the classical Gromov--Hausdorff distance is
that if the distance between two compact metric spaces is zero, then the
spaces are isometric.  In the next section we will prove the corresponding
fact for compact quantum metric spaces.  But in preparation for this we need
to discuss what we mean by isometries in the quantum case.

Let $(A,L_A)$ and $(B,L_B)$ be compact quantum metric spaces.  Let $\var$ be
an order-unit isomorphism from $A$ onto $B$ such that $L_B(\var(a)) =
L_A(a)$ for all $a \in A$.  Let $\var'$ denote the corresponding
$w^*$-continuous affine bijection of $S(B)$ onto $S(A)$.  Then it is easily
seen that $\var'$ is an ordinary isometry from $(S(B),\rho_{L_B})$ onto
$(S(A),\rho_{L_A})$.  Thus, we certainly want to consider $\var$ to be
an isometry from $(A,L_A)$ to $(B,L_B)$.

But this is not quite sufficient for our purposes.  For example, let $T$ be
the unit circle, and let $A$ consist of the smooth functions on $T$.  View
$A$ as consisting of functions on ${\mathbb R}$ periodic of period~1.
Define $L$ on $A$ just as in example $3.5$ of \cite{R5}, that is, by
\[
L(f) = \|f'\|_{\infty} + |f'(0)|.
\]
It is not difficult to show that $L$ is a Lip-norm for our present
definition, and that $\rho_L$ is the usual metric on $T$ coming from that on
${\mathbb R}$.  Thus, we would expect that every rotation of $T$ gives an
isometry of $(A,L)$.  But it is clear that the only rotation leaving $L$
invariant is the identity rotation.  Now as seen in example $3.5$ of
\cite{R5}, $L$ is not lower semi-continuous.  It is not difficult to show
that the largest lower semi-continuous seminorm, $L^s$, smaller than $L$,
whose existence is assured by theorem $4.2$ of \cite{R5}, is given by
$$
L^s(f) = \|f'\|_{\infty} .   \leqno*)   
$$
In accordance with theorem $4.2$ of \cite{R5}, $L^s$ gives the same metric as
$L$, and we see that $L^s$ is invariant under all rotations of $T$.  This
suggests that we need to use lower semi-continuous Lip-norms.

But this is not enough.  Choose some fixed $t_0 \in T$, and define $A$ to
be, instead, the algebra of all smooth functions on $T$ which are constant
on some neighborhood of $t_0$, where the neighborhood depends on the
function.  Define $L$ on this $A$ by $*)$ above.  It is easy to see that $L$
is still lower semi-continuous, and that $\rho_L$ still gives the usual
metric on $T$.  But it is clear that the only rotation of $T$ which carries
$A$ into itself is the identity rotation.

It appears that the nicest solution to these difficulties is to require that
$(A,L)$ be closed, as defined in definition $4.5$ of \cite{R5}, and recalled
above in connection with Proposition $3.3$.  By proposition $4.4$ of
\cite{R5} every lower semi-continuous Lip-norm extends uniquely to a closed
Lip-norm.  We have the following attractive description of the closure, which
was not given in \cite{R5}:

\bigskip
\noindent
{\bf 6.1 PROPOSITION.} {\em Let $(A,L)$ be a compact quantum metric space.
If $L$ is closed, then $A$ corresponds to the subspace of $Af(S(A))$
consisting of all affine functions on $S(A)$ which are ordinary Lipschitz
functions for $\rho_L$.  Conversely, if $L$ is lower semi-continuous and if
every affine function on $S(A)$ which is Lipschitz for $\rho_L$ corresponds
to an element of $A$, then $L$ is closed.}

\bigskip
\noindent
{\it\bfseries Proof.} The proof of the first assertion is very similar to
the proof of theorem $4.2$ of \cite{R5}.  We include the argument here for
the reader's convenience.  Let ${\mathcal L}_1$ denote the ``unit-ball'' for
$L$, that is, $\{a \in A: L(a) \le 1\}$.  We view its elements as functions
contained in $Af(S(A))$, and examine its bipolar \cite{Cw} there, where
$Af(S(A))$ is equipped with its order-unit norm.  As in \cite{R5}, we denote
the norm-dual of $Af(S(A))$, which coincides with that of $A$, by $A'$.
Since ${\mathbb R}e \subseteq {\mathcal L}_1$, the polar, $({\mathcal
L}_1)^O$, of ${\mathcal L}_1$ will be contained in $A'{}^{\circ}$, the set of
$\l \in A'$ for which $\l(e) = 0$.  In fact, essentially by definition,
\[
({\mathcal L}_1)^O = \{\l \in A'{}^{\circ}: |\l(a)| \le 1 \mbox{ for all $a$
with } L(a) \le 1\}.
\]
But this just says that $({\mathcal L}_1)^O$ is the unit-ball for the ``dual
seminorm'', $L'$, to $L$ on $A'{}^{\circ}$.  Thus the bipolar,
${}^O(({\mathcal L}_1)^O)$, of ${\mathcal L}_1$ is given by
\[
{}^O(({\mathcal L}_1)^O) = \{b \in Af(S(A)): |\l(b)| \le 1 \mbox{ for all
$\l$ with } L'(\l) \le 1\}.
\]
Equivalently, $b \in {}^O(({\mathcal L}_1)^O)$ exactly if $|\l(b)| \le
L'(\l)$ for all $\l \in A'{}^{\circ}$.  As in \cite{R5}, let $D_2 = \{\l \in
A'{}^{\circ}: \|\l\| \le 2\}$.  Then, clearly, $b \in {}^O(({\mathcal
L}_1)^O)$ exactly if $|\l(b)| \le L'(\l)$ for all $\l \in D_2$.  But by
lemma $2.1$ of \cite{R5} every $\l \in D_2$ is (not uniquely) of the form $\l
= \mu - \nu$ for $\mu,\nu \in S(A)$, and conversely.  Furthermore, by lemma
$4.3$ of \cite{R5} we then have $L'(\l) = \rho_L(\mu,\nu)$.  Thus, we see
that $b \in {}^O(({\mathcal L}_1)^O)$ exactly if $|\mu(b) - \nu(b)| \le
\rho_L(\mu,\nu)$ for all $\mu,\nu \in S(A)$.  But this says exactly that
$L_{\rho_L}(b) \le 1$.  Thus, the bipolar of ${\mathcal L}_1$ in $Af(S(A))$
coincides with the elements of $\rho_L$-Lipschitz norm $\le 1$.  Now,
${\mathcal L}_1$ is already convex and balanced, so the bipolar theorem
\cite{Cw} says that the bipolar of ${\mathcal L}_1$ is the norm-closure of
${\mathcal L}_1$.  But we are assuming that ${\mathcal L}_1$ is
norm-closed.  Thus, ${\mathcal L}_1$ coincides with the set of elements of
$Af(S(A))$ of $\rho_L$-Lipschitz norm $\le 1$.  Consequently, $A$ coincides
with the set of elements of $Af(S(A))$ which are Lipschitz for $\rho_L$, as
desired.

Conversely, suppose $A$ corresponds to the set of all affine functions on
$S(A)$ which are Lipschitz for $\rho_L$ (and so necessarily continuous).
It is easily seen that the set of $a \in Af(S(A))$ with $L_{\rho_L}(a) \le
1$ is norm-closed.  Since we assume now that $L$ is lower semi-continuous,
we have $L = L_{\rho_L}$ by theorem $4.2$ of \cite{R5}.  Thus ${\mathcal
L}_1$ for $L$ is closed in ${\bar A}$, and this means that $L$ is closed.\qed

\bigskip
\noindent
{\bf 6.2 THEOREM.} {\em Let $(A,L_A)$ and $(B,L_B)$ be compact quantum
metric spaces which are closed.  For every affine map $\a$ from $S(B)$ onto
$S(A)$ which is isometric for $\rho_{L_B}$ and $\rho_{L_A}$ there is a
(unique) order isomorphism, $\var$, of $A$ onto $B$ such that $\a = \var'$
(that is, $\a(\nu)(a) = \nu(\var(a))$ for all $\nu \in S(B)$ and $a \in A$)
and such that $L_A = L_B \circ \var$.  The converse is also true.}

\bigskip
\noindent
{\it\bfseries Proof.} We already saw at the beginning of this section that
the converse holds.  Now let $\a$ be as above.  Since $\a$ is isometric, it
is continuous.  Thus, composing with $\a$ gives an order isomorphism,
$\var$, of $Af(S(A))$ onto $Af(S(B))$, and so from $A$ onto $B$ since we are
assuming that $L_A$ and $L_B$ are closed (and we use Proposition $6.1$).
Since $\a$ is isometric, we will have $L_{\rho_{L_B}}(\var(a)) =
L_{\rho_{L_A}}(a)$ for all $a \in A$.  But because $L_A$ and $L_B$ are
closed, they are lower semi-continuous, so that $L_{\rho_{L_A}} = L_A$ on
$A$, and similarly for $B$.  Thus, $L_A = L_B \circ \var$ as desired.\qed

\bigskip
For any compact quantum metric space $(A,L_A)$ let $L_A^s$ denote, as
before, the greatest lower semi-continuous Lip-norm smaller than $L_A$, and
let $(A^c,L_A^c)$ denote the closure of $(A,L_A^s)$.

\bigskip
\noindent
{\bf 6.3 DEFINITION.} Let $(A,L_A)$ and $(B,L_B)$ be compact quantum metric
spaces.  By an {\em isometry} from $(A,L_A)$ to $(B,L_B)$ we mean an order
isomorphism $\var$ from $A^c$ onto $B^c$ such that $L_A^c = L_B^c \circ
\var$.

\bigskip
  From Theorem $6.2$ we immediately obtain:

\bigskip
\noindent
{\bf 6.4 COROLLARY.} {\em The isometries from $(A,L_A)$ to $(B,L_B)$ are in
natural bijective correspondence with the affine isometries from
\linebreak
$(S(B),\rho_{L_B})$ onto $(S(A),\rho_{L_A})$.}

\bigskip
Isometries can be composed in the evident way.  There is, of course, a good
argument for saying that $\var$ as in Definition $6.3$ is ``from $(B,L_B)$
onto $(A,L_A)$'' rather than the other way around.

By the comments at the beginning of this section, if we happen to be given an
order isomorphism $\var$ from $A$ to $B$ such that $L_A = L_B \circ \var$,
then it will lead to an isometry as defined in Definition $6.3$.  Of course,
the reason we do not always work with Lip-norms which are closed is that
the Lip-norms which arise in examples are often not closed, e.g., the case
considered earlier in which $A$ consists of the smooth functions on $T$.  In
fact, we see that to a large extent our theory is about compact convex sets
(e.g., $S(A)$) and metrics on them which are suitably convex in the sense
discussed in section~9 of \cite{R5}.  But again, this is not the form in
which the examples arise.

We remark that our definition of isometries is more general than the
definition given in \cite{P1}, \cite{P2}, since that definition assumes that
$L$ comes from a specific spectral triple, and that the isometry is unitarily
implemented on the Hilbert space of that spectral triple. A definition of
isometries in the context of spectral triples which 
is closer to our definition
makes a brief appearance in the third paragraph after corollary III.1.5
of \cite{Lo}.

We now give a class of examples of isometries which is closely related to
the considerations which we will see in Section~8.

\bigskip
\noindent
{\bf 6.5 EXAMPLE.} As in \cite{R4}, let $G$ be a compact group and let
$\ell$ be a continuous length-function on $G$.  Let $\a$ be an ergodic
action of $G$ on a unital $C^*$-algebra ${\bar A}$.  For $a \in {\bar A}$ let
\[
L(a) = \sup\{\|\a_x(a) - a\|/\ell(x): x \ne e_G\},
\]
and let $A = \{a \in {\bar A}: L(A) < \infty\}$ (or, more precisely,
the real subspace of self-adjoint elements).  Then $(A,L)$ is a compact
quantum metric space for which $L$ is closed (Proposition 8.1).  
As remarked in section~2 of
\cite{R4}, it is not so clear whether $A$ is carried into itself by $\a$.
However, if $\ell$ satisfies the additional hypothesis that $\ell(zxz^{-1})
= \ell(x)$ for all $x,z \in G$, then a little calculation given in section~2
of \cite{R4} shows that $\a$ does carry $A$ into itself; and furthermore,
that $\a$ leaves $L$ invariant.  In other words, $G$ acts as a group of
isometries of $(A,L)$ when $\ell$ satisfies this extra condition.

\bigskip
Let $(A,L)$ be a compact quantum metric space, and let $G = \Iso(A,L)$ be
its group of isometries.  It is natural to define a length function,
$\ell$, on $G$ by
\[
\ell(\var) = \sup\{\|\var(a) - a\|: a \in A^c \mbox{ and } L^c(a) \le 1\}.
\]
A quick calculation then shows that $\ell$ has the special property used in
Example $6.5$, namely, $\ell(\psi\var\psi^{-1}) = \ell(\var)$ for all
$\var,\psi \in G$.  Thus, this special condition is quite natural in our
context. The following proposition extends a well-known fact about the
isometry groups of compact manifolds, and the comments surrounding
equation 1.30 of \cite{Cnr}.

\bigskip
\noindent
{\bf 6.6 PROPOSITION.} {\em Let $(A,L)$ be a compact quantum metric space 
for which $L$ is complete. Let $G = Iso(A,L)$ equipped with the length
function $\ell$ defined above, and the corresponding left-invariant
metric, $\rho$. Then $G$ is a compact group for the topology defined
by $\rho$.} 

\bigskip
\noindent
{\it\bfseries Proof.} Let $r$ denote the radius of $(A,L)$ as defined
at the end of Section 2, and let
$$
{\mathcal B}_r = \{a \in A:L(a) \le 1 \mbox { and } \|a\| \le r\}.
$$
Because $\|a\|^{\sim} \le rL(a)$ for all $a \in A$, and
$\var(e) - \psi(e) = 0$ for all $\var, \psi \in G$, an equivalent
definition of $\rho$ is
$$
\rho(\var, \psi) = \sup\{\|\var(a) - \psi(a)\|: a \in {\mathcal B}_r\}.
$$
Note that ${\mathcal B}_r$ is carried into itself by any element of $G$.
Because $L$ is closed, ${\mathcal B}_r$ is closed in $\bar A$, and so
is compact for the norm topology. But each element of $G$ is an
isometry for the norm, and so $G$ is a bounded equicontinuous collection
of maps from ${\mathcal B}_r$ into itself. Thus $G$ is totally bounded
by the Arzela-Ascoli theorem for the point-norm topology, which is
exactly the topology determined by $\rho$. But it is easily verified 
that $G$ is closed in the set of all continuous maps from ${\mathcal B}_r$
into itself. Thus $G$ is compact. Since $\rho$ comes from a length function
on $G$, its topology is compatible with the group structure of $G$. 
(A slightly less elementary proof can be given by using Corollary 6.4.)
\qed
\bigskip

It would be interesting to develop and study the notion of ``quantum 
isometry group" for quantum metric spaces as quantum subgroups of the quantum 
symmetry groups studied by Wang \cite{Wng}.

\section{Distance zero}

The aim of this section is to show that if $\dist_q(A,B) = 0$ then there is
an isometry from $A$ to $B$, in the sense defined in the previous section.
We note, conversely, that if there is an isometry, $\psi$, from $A$ to $B$
(so defined on $A^c$), then indeed $\dist_q(A,B) = 0$.  To see this, for any
$\g > 0$ define a bridge from $A^c$ to $B^c$ by
\[
N(a,b) = \g^{-1}\|\psi(a)-b\|.
\]
It is easily seen that for the corresponding Lip-norm $L$ on $A^c \oplus
B^c$ we have $\dist_{GH}^{\rho_L}(S(A^c), S(B^c)) \le \g$.  
We also need to clarify the
situation with:

\bigskip
\noindent
{\bf 7.1 PROPOSITION.} {\em Let $(A,L)$ be a compact quantum metric space,
let $L^s$ denote the largest lower semi-continuous seminorm smaller than $L$
and let $(A^c,L^c)$ denote the closure of $(A,L)$.  Then}
\[
\dist_q((A,L),(A,L^s)) = 0 = \dist_q((A,L),(A^c,L^c)).
\]

\bigskip
\noindent
{\it\bfseries Proof.} We consider first the case of $L^s$.  Let $\g > 0$ be
given, and define $N$ on $A \oplus A$ by
\[
N(a,b) = \g^{-1}\|a-b\|.
\]
We must show that $N$ is a bridge for $L$ and $L^s$.  
The first two conditions of Definition
$5.1$ are clearly satisfied.  We now check Condition~3.  Let $A_1$ and $A_2$
denote $A$ as first or second copy in $A \oplus A$.  For $a \in A_1$ set $b
= a$ in $A_2$.  This works since $L^s(a) \le L(a)$.  Suppose, however, that
we are given $b \in A_2$.  We recall from theorem $4.2$ of \cite{R5} that
the ``unit $L^s$-ball'' is the norm closure in $A$ (not ${\bar A}$) of the
``unit $L$-ball''.  Thus we can find a sequence $\{a_n\}$ of elements in
$A$ such that $L(a_n) \le L^s(b)$ and $\{a_n\}$ converges to $b$ in norm.
Consequently, we can find an $a_n$, viewed as in $A_1$, such that
$\g^{-1}\|a_n-b\| \le L^s(b)$.  This shows that $N$ is a bridge.  For the
corresponding Lip-norm on $A \oplus A$ for $N$, say $M$ (briefly, $M = (L
\vee L^s) \vee N)$, it is now routine to verify that
\[
\dist_H^{\rho_M}(S(A_1),S(A_2)) \le \g.
\]
Since $\g$ is arbitrary, we see that $\dist_q(A_1,A_2) = 0$.

The proof for the case of $(A^c,L^c)$ is very similar, but now we must
recall from proposition $4.4$ of \cite{R5} that the ``unit $L^c$-ball'' is
the closure in ${\bar A}$ of the ``unit $L$-ball'', and notice that the norm
in the definition above of $N$ must be that for $A^c$.\qed

\bigskip
In developing our proof that if $\dist_q(A,B) = 0$ then there is an isometry
between $A$ and $B$, we find it useful to first give a proof of the
corresponding theorem for ordinary compact metric spaces, but along lines
which then can be used for the proof of the quantum case.  A significant
obstacle to the proof in either case is that there is, in general, no
uniqueness to the isometry --- it can be composed with isometries of either
of the individual spaces to get other isometries.  This means that some
choice principle must be used, usually in the form of compactness.  The
compactness which we find convenient to use here is that involved in the
Arzela--Ascoli theorem \cite{Cw}.  In my browsing through the extensive
literature concerning Gromov--Hausdorff distance I have not noticed the use
of the Arzela--Ascoli theorem for this specific purpose.  But I would not be
surprised if this use does appear somewhere, since I have seen the
Arzela--Ascoli theorem used for closely related purposes in the literature.

Let $(X,\rho_X)$ and $(Y,\rho_Y)$ be ordinary compact metric spaces.  We
wish to work on one fixed space containing them, namely the disjoint union
$X \dotcup Y$.  To compensate for this we must admit semi-metrics, that is,
continuous functions which satisfy all the properties of a metric except
that they may take value $0$ on some pairs of distinct points.  We denote by
$\Sigma(\rho_X,\rho_Y)$ the family of all semi-metrics on $X \dotcup Y$
whose restrictions to $X$ and $Y$ are $\rho_X$ and $\rho_Y$, respectively.

\bigskip
\noindent
{\bf 7.2 LEMMA.} {\em The family $\Sigma(\rho_X,\rho_Y)$ of functions on $(X
\dotcup Y)^2$ is equicontinuous.}

\bigskip
\noindent
{\it\bfseries Proof.} Both $X \x X$ and $Y \x Y$ are open subsets of $(X
\dotcup Y)^2$, and all the functions in $\Sigma(\rho_X,\rho_Y)$ agree on
these two subsets, so equicontinuity is obvious there.  Suppose 
we are given $x_0 \in X$ and $y_0 \in Y$, as well as $\e > 0$.  Let
\[
{\mathcal N}_{x_0} = \{x \in X: \rho_X(x,x_0) < \e/2\},
\]
and similarly for ${\mathcal N}_{y_0}$, so that ${\mathcal N}_{x_0} \x
{\mathcal N}_{y_0}$ is a neighborhood of $(x_0,y_0)$.  If $(x,y) \in
{\mathcal N}_{x_0} \x {\mathcal N}_{y_0}$, then for any $\s \in
\Sigma(\rho_X,\rho_Y)$ we have
\begin{eqnarray*}
|\s(x,y) - \s(x_0,y_0)| &\le &|\s(x,y) - \s(x,y_0)| + |\s(x,y_0) -
\s(x_0,y_0)| \\
&\le &\rho_Y(y,y_0) + \rho_X(x,x_0) < \e.
\end{eqnarray*}
\qed

\bigskip
\noindent
{\bf 7.3 LEMMA.} {\em The family $\Sigma(\rho_X,\rho_Y)$ is uniformly closed
in $C((X \dotcup Y)^2)$.}

\bigskip
\noindent
{\it\bfseries Proof.} All of the conditions in the definition of a metric
are closed conditions except the condition that the distance between
distinct points must be non-zero.  It is for this reason that we drop this
condition and allow semi-metrics.\qed

\bigskip
\noindent
{\bf 7.4 LEMMA.} {\em Let $\s \in \Sigma(\rho_X,\rho_Y)$.  For each $x \in
X$ there is at most one $y \in Y$ such that $\s(x,y) = 0$; and similarly for
each $y\in Y$.}

\bigskip
\noindent
{\it\bfseries Proof.} If $\s(x,y) = 0 = \s(x,y')$, then $\rho_Y(y,y') \le
\s(y,x) + \s(x,y') = 0$, so that $y' = y$.\qed

\bigskip
The condition that two points of $X \dotcup Y$ have $\s$-distance zero is an
equivalence relation, $\sim$, on $X \dotcup Y$.  By Lemma $7.4$ each
equivalence class contains either one point, or a pair $(x,y) \in X \x Y$.
We set $X \cup_{\s} Y = (X \dotcup Y)/\sim$.  Then $\s$ drops to a genuine
metric, $\rho_{\s}$, on $X \cup_{\s} Y$.  Furthermore, the map $i_X: X
\rightarrow X \cup_{\s} Y$ defined by $i_X(x) = {\tilde x}$ is clearly an
isometry from $X$ into $X \cup_{\s} Y$, and similarly for $i_Y$.  We clearly
have $X \cup_{\s} Y = i_X(X) \cup i_Y(Y)$, not necessarily disjointly, and
$X \cup_{\s} Y$ is clearly compact.

Suppose, conversely, that we have a compact metric space $(Z,\rho)$ together
with isometries $j_X$ and $j_Y$ of $(X,\rho_X)$ and $(Y,\rho_Y)$ into $Z$
such that $Z = j_X(X) \cup j_Y(Y)$.  It is evident how to define from this a
semi-metric, $\s$, on $X \dotcup Y$ such that there is an isometry (unique)
from $(X \cup_{\s} Y,\rho_{\s})$ onto $(Z,\rho)$ which carries $(i_X,i_Y)$
to $(j_X,j_Y)$.  From these considerations we obtain:

\bigskip
\noindent
{\bf 7.5 PROPOSITION.} {\em Let $(X,\rho_X)$ and $(Y,\rho_Y)$ be compact
metric spaces.  There is a natural bijection between $\Sigma(\rho_X,\rho_Y)$
and the set of equivalence classes of objects $(Z,\rho,j_X,j_Y)$ where
$(Z,\rho)$ is a compact metric space and $j_X$ and $j_Y$ are isometries of
$X$ and $Y$ into $Z$ such that $Z = j_X(X) \cup j_Y(Y)$, under the
equivalence relation of surjective isometries which preserve $(j_X,j_Y)$.
Under this bijection $\s \in \Sigma(\rho_X,\rho_Y)$ is sent to the
equivalence class of $(X \cup_{\s} Y,\rho_{\s})$.}

\bigskip
\noindent
{\bf 7.6 THEOREM.} {\em Let $(X,\rho_X)$ and $(Y,\rho_Y)$ be compact metric
spaces, and let $d = \dist_{GH}(X,Y)$.  Then there is a compact metric space
$(Z,\rho)$, and isometric injections $i_X$ and $i_Y$ of $X$ and $Y$ into
$Z$, such that}
\[
\dist_H^{\rho}(i_X(X),i_Y(Y)) = d,
\]
{\em and $Z = i_X(X) \cup i_Y(Y)$.}

\bigskip
\noindent
{\it\bfseries Proof.} From the definition of $d$, there is a sequence,
$\{\rho_n\}$, of metrics in $\Sigma(\rho_X,\rho_Y)$ such that within $X
\dotcup Y$
\[
\dist_H^{\rho_n}(X,Y) \le d + 1/n.
\]
This sequence is easily seen to be uniformly bounded by
\[
(\mbox{diameter}(X) \vee \mbox{diameter}(Y)) + d + 1.
\]
Since the sequence is also equicontinuous by Lemma $7.2$, we can apply the
Arzela--Ascoli theorem to conclude that there will be a uniformly convergent
subsequence.  For simplicity we still denote this subsequence by
$\{\rho_n\}$.  Let $\s$ be its limit.  Since $\Sigma(\rho_X,\rho_Y)$ is
closed, $\s$ must be a semi-metric.  Set $Z = X \cup_{\s} Y$ and $\rho =
\rho_{\s}$ as constructed above, with $i_X$ and $i_Y$ the corresponding
isometric inclusions.

Let us determine $\dist_H^{\rho}(i_X(X),i_Y(Y))$.  It cannot be smaller than
$d$ by the definition of $d$.  But let $x \in X$.  Then for each $n$ there
is a $y_n \in Y$ such that $\rho_n(x,y_n) \le d +1/n$.  The sequence
$\{y_n\}$ has a convergent subsequence, which for simplicity, we denote
again by $\{y_n\}$.  Relabel the $\rho_n$'s accordingly.  Let $y_0$ be the
limit of $\{y_n\}$.  Given $\e > 0$, choose $N$ such that if $n > N$ then
$\|\rho_n - \s\|_{\infty} < \e/3$, and $1/n < \e/3$, and $\rho_Y(y_n,y_0) <
\e/3$.  Then for $n > N$ we have
\begin{eqnarray*}
\s(x,y_0) &\le &\s(x,y_n) + \s(y_n,y_0) \\
&\le &\rho_n(x,y_n) + \e/3 + \rho_Y(y_n,y_0) \\
&\le &d + \e/3 + 2\e/3 = d + \e.
\end{eqnarray*}
Since $\e$ is arbitrary, $\s(x,y_0) \le d$.  In the same way, for each $y
\in Y$ there is an $x_0 \in X$ such that $\s(y,x_0) \le d$.  These
inequalities pass to $\rho$ on $Z$, and show that
$\dist_H^{\rho}(i_X(X),i_Y(Y)) = d$, as desired.\qed

\bigskip
If $d = 0$ then
$\dist_H^{\rho}(i_X(X),i_Y(Y)) = 0$ under the circumstances of 
the above theorem, so that $i_X(X) = i_Y(Y)$
since they are both closed.  Thus, $i_X$ and $i_Y$ are surjective, and
$i_Y^{-1}
\circ i_X$ is an isometry from $X$ onto $Y$.  In this way we obtain:

\bigskip
\noindent
{\bf 7.7 THEOREM.} {\em Let $(X,\rho_X)$ and $(Y,\rho_Y)$ be compact metric
spaces.  If}
\[
\dist_{GH}(X,Y) = 0,
\]
{\em then there is an isometry from $X$ onto $Y$.}

\bigskip
We now turn to the quantum case.  For the  next theorem it is essential that
by ``isometry'' we mean in the sense discussed in Section~6, involving
closures.

\bigskip
\noindent
{\bf 7.8 THEOREM.} {\em Let $(A,L_A)$ and $(B,L_B)$ be compact quantum
metric spaces.  If}
\[
\dist_q(A,B) = 0,
\]
{\em then there is an isometry between $(A,L_A)$ and $(B,L_B)$.}

\bigskip
\noindent
{\it\bfseries Proof.} If $\dist_q(A,B) = 0$, then there is a sequence,
$\{L_n\}$, of Lip-norms on $A \oplus B$ inducing $L_A$ and $L_B$ such that
\[
\dist_H^{\rho_{L_n}}(S(A),S(B)) < 1/n.
\]
For notational simplicity set $\rho_n = \rho_{L_n}$, and $\rho_A =
\rho_{L_A}$, $\rho_B = \rho_{L_B}$.  As seen in Proposition $3.1$, the
restrictions of $\rho_n$ to $S(A)$ and $S(B)$ are $\rho_A$ and $\rho_B$,
respectively.  Of course, $S(A)$ and $S(B)$ are disjoint subsets of $S(A
\oplus B)$.  Thus, if we view $\rho_n$ as restricted to their disjoint
union, this says that $\rho_n \in \Sigma(\rho_A,\rho_B)$ in the notation
introduced early in this section.  Exactly as in the proof of Theorem $7.6$,
there is a subsequence, which we still denote by $\{\rho_n\}$, which
converges uniformly on $S(A) \dotcup S(B)$ to a semi-metric, $\s$.  Exactly
as in the proof of Theorem $7.7$ we see that $\s$ determines an isometry,
$\a$, from $S(A)$ onto $S(B)$, by the condition that $\s(\mu,\a(\mu)) = 0$.
According to Theorem $6.2$, in order to show that $\a$ gives an isometry
from $B^c$ onto $A^c$ it suffices to show that $\a$ is affine.  Thus, let
$\mu_1,\mu_2 \in S(A)$ and let $t \in [0,1]$.  Let $\e > 0$ be given, and
find $N$ such that if $n \ge N$ then $\|\s - \rho_n\|_{\infty} < \e/2$.
Then for any $n \ge N$ we have
\begin{eqnarray*}
&&\s(t\mu_1 + (1-t)\mu_2,t\a(\mu_1) + (1-t)\a(\mu_2)) \\
&&\le \rho_n(t\mu_1 +
(1-t)\mu_2,t\a(\mu_1) + (1-t)\a(\mu_2)) + \e/2 \\
&&= L'_n(t(\mu_1-\a(\mu_1)) + (1-t)(\mu_2-\a(\mu_2))) + \e/2 \\
&&\le tL'_n(\mu_1 - \a(\mu_1)) + (1-t)L'_n(\mu_2 - \a(\mu_2)) + \e/2 \\
&&= t\rho_n(\mu_1,\a(\mu_1)) + (1-t)\rho_n(\mu_2,\a(\mu_2)) + \e/2 \\
&&\le t(\s(\mu_1,\a(\mu_1)) + \e/2) + (1-t)(\s(\mu_2,\a(\mu_2)) + \e/2) +
\e/2 \\
&&\le \e.
\end{eqnarray*}
Since $\e$ is arbitrary, we conclude that
\[
\a(t\mu_1 + (1-t)\mu_2) = t\a(\mu_1) + (1-t)\a(\mu_2),
\]
as needed.\qed

\bigskip
We remark that a $C^*$-algebra and its opposite algebra need not be
isomorphic; but their corresponding order-unit spaces of self-adjoint
elements will, nevertheless, be isomorphic as order-unit spaces.  Thus, in
the above theorem, if $A$ and $B$ happen to be the order-unit spaces of
self-adjoint elements of two $C^*$-algebras, it does not follow that those
$C^*$-algebras must be isomorphic. Interesting examples of $C^*$-algebras
which are not isomorphic to their opposite algebra can be found in \cite{Ph}.
Well after this paper was submitted for publication, two approaches to
dealing with this unsatisfactory feature have been given. David Kerr
\cite{Ker} has extended much of the present paper to operator systems 
with their matricial structure, where he can take advantage of the fact 
that if two $C^*$-algebras 
are {\em completely} order-isomorphic, in terms of their matricial norms,
then they are isomorphic as $C^*$-algebras. (See corollary 5.2.3 
of \cite{ER}.) My former doctoral student Hanfeng Li has shown in
his doctoral thesis \cite{Li} how to define a Gromov-Hausdorff
distance between $C^*$-algebras with Lip-norms which explicitly uses
the product in the $C^*$-algebras, and has the desired property.

For simplicity of notation in the next corollaries, we will not 
distinguish between $C(X)$ and its dense subalgebras of Lipschitz
functions, but all this should be interpreted as in the discussion
following the proof of Lemma 4.6.

\bigskip
\noindent
{\bf 7.9 Corollary.} {\em Let $X$ be a compact space, and let $L$ be
a Lip-norm on $C(X)$, not necessarily coming from a metric on $X$.
Let $(B,L_B)$ be a compact quantum metric space. If}
\[
\dist_q((C(X), L), (B, L_B)) = 0,
\]
{\em then the completion of $B$ is order-isomorphic to $C(X)$.}

\bigskip
\noindent
{\it\bfseries Proof.} By Theorem 7.8 there will be an isometry between the
two spaces. Now by Definition 6.3 an isometry is, in particular, an
order-isomorphism between the closures. But any order isomorphism is 
norm-continuous, and so extends to an order-isomorphism of the norm 
completions. \qed

\bigskip
\noindent
{\bf 7.10 Corollary.} {\em Let $(X, \rho_X)$ and $(Y, \rho_Y)$ be 
compact metric spaces. Let $L_X$ and $L_Y$ be the (closed) Lip-norms 
determined by $\rho_X$ and $\rho_Y$, so that $(C(X), L_X)$ and $(C(Y), L_Y)$
can be viewed as compact quantum metric spaces. If}
\[
\dist_q((C(X), L_X), (C(Y), L_Y)) = 0,
\]
{\em then}
\[
\dist_{GH}((X, \rho_X), (Y, \rho_Y)) = 0.
\]

\bigskip
\noindent
{\it\bfseries Proof.} From Theorem 7.7 there is an isometry from
$(C(X), L_X)$ to $(C(Y), L_Y)$. By the definition of an isometry, it must,
in particular, extend as above to 
an order-isomorphism from $C(X)$ to $C(Y)$. 
But by corollary 3.4.8 of \cite{KdR}, 
such an order-isomorphism is an algebra isomorphism, 
which then corresponds to a 
homeomorphism between $X$ and $Y$. It is easily seen that this 
homeomorphism is an isometry, so that $\dist_{GH}(X, Y) = 0$ as desired. \qed

For further results along these lines see Theorem 13.16.

\section{Actions of compact groups}

Our aim now is to apply to quantum tori the theory which we have developed.
But one of our main steps in doing this can be treated in a somewhat more
general framework, useful in other situations.  We will discuss this step in
this section.

Much as in \cite{R4} and in Example $6.5$, we consider a compact group $G$,
equipped with a length function $\ell$ which satisfies the additional
condition that $\ell(xyx^{-1}) = \ell(y)$ (which we saw in Example $6.5$
insures that $G$ will act by isometries).  We let $\a$ be an ergodic
action of $G$ on a unital $C^*$-algebra ${\bar A}$ (where ``ergodic'' means
simply that the fixed-point subalgebra is ${\mathbb C}1_A$).  We define $L$
on ${\bar A}$ by
\[
L(a) = \sup_G\{\|\a_x(a) - a\|/\ell(x): x \ne e_G\},
\]
and we set $A = \{a \in {\bar A}: L(a) < \infty\}$.  As shown in \cite{R4},
$A$ is a dense $*$-subalgebra of ${\bar A}$, and $L$ is a Lip-norm on $A$.
To remain strictly within Definition $2.1$ we should restrict further to the
order-unit space of all self-adjoint elements of $A$.  Since the above $L$
clearly satisfies $L(a^*) = L(a)$, this is equivalent, by the comments made
just before Definition $2.1$.  Consequently, we will be a bit careless here
about this distinction.

Notice further that the above $L$ is defined as a supremum of functions on
${\bar A}$ which are continuous, and so $L$ is lower semi-continuous, though
it may take value $+\infty$ on ${\bar A}$.  From this we easily see:

\bigskip
\noindent
{\bf 8.1 PROPOSITION.} {\em Let $A$ and $L$ be defined as above in terms of
an ergodic action $\a$ of the compact group $G$ on the unital $C^*$-algebra
${\bar A}$ and a length function $\ell$ on $G$.  Then $(A,L)$ is closed.}

\bigskip
We will now use $\a$ to single out certain order-unit subspaces of $A$.  For
this we must assume that $G$ has a faithful finite-dimensional unitary
representation, say $\pi_0$.  Thus $G$ is a Lie group, possibly
disconnected (or even a finite group).  Adjust $\pi_0$ so that it contains
the trivial representation, and set $\pi = \pi_0 \otimes {\bar \pi}_0$,
where ${\bar \pi}_0$ denotes the contragradient representation.  Thus $\pi$
is a faithful representation containing the trivial representation, and the
character, $\chi$, of $\pi$ is a non-negative real-valued function on $G$.
Let ${\hat G}$ denote the dual of $G$, that is, the set of equivalence
classes of irreducible unitary representations of $G$.  Let ${\hat G}_n$
denote the finite subset of ${\hat G}$ consisting of those irreducible
representations which occur in $\pi^{\otimes n}$.  Notice that ${\hat G}_n$
is closed under taking contragradient representations, and contains the
trivial representation.

For each $n$ let $B_n$ denote the direct sum of the (finite number of)
$\a$-isotypic subspaces of ${\bar A}$ corresponding to all the elements of
${\hat G}_n$.  Because ${\hat G}_n$ is closed under taking contragradients,
$B_n$ is closed under taking adjoints.  Because ${\hat G}_n$ contains the
trivial representation, $B_n$ contains the identity element of ${\bar A}$.
Thus, $B_n$ (or more precisely the real part of $B_n$) is an order-unit
subspace of (the real part of) ${\bar A}$.

We now need the quite unobvious fact \cite{HLS}, \cite{Wa} that because $\a$
is ergodic, each isotypic component of ${\bar A}$ is finite-dimensional.
Thus each $B_n$ is finite-dimensional.  Now $\a$ carries each $B_n$ into
itself, and the Lipschitz elements for this action (for $\ell$) will form a
dense subspace, and so be all of $B_n$.  
(The proof is much like that of proposition
$2.2$ of \cite{R4}.)  Thus, each $B_n$ is in $A$, not just in ${\bar A}$.
This says that our $L$ above is finite on $B_n$.  Thus $(B_n,L)$ is a
compact (in fact ``finite'') quantum metric space.  The main theorem of this
section is:

\bigskip
\noindent
{\bf 8.2 THEOREM.} {\em Let notation be as above.  Then there is a sequence,
$\{\d_n\}$, of non-negative numbers, converging to $0$, depending only on
$\pi$ and $\ell$, but not on $A$ and $\a$, such that}
\[
\dist_q(A,B_n) \le \d_n
\]
{\em for each $n$.}

\bigskip
\noindent
{\it\bfseries Proof.} We begin by specifying the $\d_n$'s.  With $\chi$ the
character of $\pi$ as before, $\chi^n$ (pointwise product) is the character
of the inner tensor-power $\pi^{\otimes n}$.  Fix a Haar measure on $G$, and
set
\[
\var_n = \chi^n/\|\chi^n\|_1,
\]
where the norm on the right is the norm of $L^1(G)$.  Set
\[
\d_n = \int_G \var_n(x)\ell(x)dx.
\]
We show that this choice of $\d_n$'s has the properties given in the
statement of the theorem.  The sequence $\{\d_n\}$ clearly depends only on
$\pi$ and $\ell$.

Since $\pi$ was chosen so that $\chi$ is non-negative, clearly $\var_n \ge
0$ and $\int_G \var_n(x)dx = 1$ for each $n$.  In particular, $\d_n \ge 0$.
Because $\pi$ is a faithful representation of $G$, it assigns the identity
operator only to $e_G$, the identity element of $G$.  For any other $x \in
G$ the unitary operator $\pi(x)$ will have some eigenvalues different from
$1$.  It follows that $\chi(e_G)  > \chi(x)$ for $x \ne e_G$.  Set 
$d = \chi(e_G)$. Note that $d \ge 2$ as long as $G$ has more than one
element. From this it
follows that the ``mass'' of $\var_n$ is increasingly concentrated
near $e_G$ as $n$ increases.  For completeness and for future use, we sketch
the simple proof of this, following the ideas in the proof of 
lemma 4.1 of \cite{Hnd}. It does not use the fact that $\chi$ is a 
character, but only the observations immediately above. Let $\epsilon > 0$ 
be given, and let $W_\epsilon = \{x:\chi(x) < (d-2\epsilon)\}$ and $N_\epsilon
= \{x:\chi(x) > (d-\epsilon)\}$. Let $|W_\epsilon|$ and $|N_\epsilon|$
denote the measures of $W_\epsilon$ and $N_\epsilon$. Then
$$
\|\chi^n\|_1 \geq \int_{N_\epsilon} \chi^n \geq |N_\epsilon|(d-\epsilon)^n,
$$
while
$$\int_{W_\epsilon} \chi^n \leq |W_\epsilon|(d-2\epsilon)^n  .
$$
Thus
$$
\int_{W_\epsilon}(\chi^n/\|\chi^n\|_1) \leq 
(|W_\epsilon|/|N_\epsilon|)((d-2\epsilon)/(d-\epsilon))^n,
$$
which goes to $0$ as $n$ grows. But the complements of
the $W_\epsilon$'s form a neighborhood base for $e_G$ as $\epsilon$
goes to 0. 

Thus $\{\var_n\}$ is an approximate identity
(of norm $1$) for the convolution algebra $L^1(G)$.  It is an analogue of
the classical Fejer kernel of harmonic analysis.  Note that each $\var_n$ is
central in $L^1(G)$ since characters are always central.  (The construction
of $\{\var_n\}$ basically appears in \cite{Hnd}, and probably in 
other places, 
but I have not seen it used in conjuction with Lipschitz seminorms
as we do below.)

Because $\ell(e_G) = 0$ and $\ell$ is continuous, it follows from the above
properties of $\{\var_n\}$ that $\{\d_n\}$ converges to $0$.

Fix $n$.  We define a bridge, $N$, between $(A,L)$ and $(B_n,L)$ by
\[
N(a,b) = \d_n^{-1}\|a-b\|,
\]
where the norm is that of $A$.  We check that $N$ satisfies the conditions
of Definition $5.1$.  It is clearly norm-continuous, and satisfies the
required properties with respect to $e_A = e_B$.  Let $b \in B_n$ be given.
Then we can choose $a = b$ to show that this part of Condition~3 is
satisfied.

So the challenge is to see, given $a \in A$, how to choose a corresponding
$b \in B$.  For any $f \in L^1(G)$ we define the operator $\a_f$ on ${\bar
A}$ as usual by $\a_f(a) = \int_G f(x)\a_x(a)dx$.  Set $P_n = \a_{\var_n}$.
Now
$\var_n$ is a (finite) linear combination of the characters of the
irreducible representations  in ${\hat G}_n$.  Consequently, the range of
$P_n$ is contained in $B_n$, by well-known properties \cite{D}, \cite{FD} of
characters.  (In fact, $P_n$ is a unit-preserving completely positive map of
${\bar A}$ onto
$B_n$.)  We now view $P_n$ as a map from $A$ into $B_n$.  We choose $b$ for
$a$ by $b =P_n(a)$.  We now start to show that this works.

\bigskip
\noindent
{\bf 8.3 LEMMA.} {\em For any $a \in A$ we 
have $\|a - P_n(a)\| \leq \d_n L(a)$.}

\bigskip
\noindent
{\it\bfseries Proof.} For $a \in A$ we have
\begin{eqnarray*}
\|a-P_n(a)\| &= 
&\left\| a \int \var_n(x)dx -\int \var_n(x)\a_x(a)dx\right\| \\
&\le &\int \var_n(x)\|a-\a_x(a)\|dx 
\le \int \var_n(x)\ell(x)L(a)dx \\
&= &\d_nL(a).
\end{eqnarray*} \qed

\bigskip
Consequently, for our choice of $b = P_n(a)$ we have 
$\d_n^{-1}\|a-b\| \le L(a)$ as needed.

\bigskip
\noindent
{\bf 8.4 LEMMA.} {\em The map $P_n$ is $\a$-equivariant, and $\|P_n\| \le
1$.  Consequently, $L(P_n(a)) \le L(a)$ for all $a \in A$.}

\bigskip
\noindent
{\it\bfseries Proof.} Because $\|\var_n\|_1 \le 1$ we have $\|P_n\| \le 1$.
A simple computation shows that because $\var_n$ is central $P_n$ is
equivariant.  Thus,
\[
\|\a_x(P_n(a)) - P_n(a)\| = \|P_n(\a_x(a) - a)\| \le \|\a_x(a) - a\|
\]
for any $x \in G$.  From this we get $L(P_n(a)) \le L(a)$.\qed

\bigskip
Consequently, for our choice of $b$ we have $L(b) \le L(a)$, which is 
the other required conditions.  Thus $N$ is a
bridge.

We can carry out the rest of the proof of Theorem $8.2$ within the more
general framework of order-unit spaces as follows:

\bigskip
\noindent
{\bf 8.5 PROPOSITION.} {\em Let $(A,L_A)$ be a compact quantum metric space,
and let $B$ be a subspace of $A$ which contains $e_A$.  Let $L_B$ denote the
restriction of $L_A$ to $B$, so that $(B,L_B)$ is a compact quantum metric
space.  Let $P$ be a function (not necessarily even linear or continuous)
from $A$ into $B$ for which there is a $\d > 0$ such that:}

\begin{itemize}
\item[1)] {\em $L_B(P(a)) \le L_A(a)$ for all $a \in A$.}
\item[2)] {\em $\|a-P(a)\| \le \d L_A(a)$ for all $a \in A$.}
\end{itemize}
{\em Then $\dist_q(A,B) \le \d$.}

\bigskip
\noindent
{\it\bfseries Proof.} It is easily verified that $L_B$ is a Lip-norm on
$B$.  Exactly as above, we define a bridge, $N$, between $A$ and $B$ by
\[
N(a,b) = \d^{-1}\|a-b\|.
\]
The conditions above insure that $N$ is indeed a bridge, where we use $P$
just as we used $P_n$ above.  We let $L$ denote the corresponding Lip-norm
on $A\oplus B$.

Suppose we are given $\mu \in S(A)$.  Let $\nu \in S(B)$ be the restriction
of $\mu$ to $B$.  Suppose that $(a,b) \in A \oplus B$ is such that $L(a,b)
\le 1$.  Then $\|a-b\| \le \d$, and so
\[
|\mu(a,b) - \nu(a,b)| = |\mu(a)-\nu(b)| = |\mu(a-b)| \le \|a-b\| \le \d.
\]
Thus, we see that $S(A)$ is in the $\d$-neighborhood of $S(B)$.

Suppose instead that we are given $\nu \in S(B)$.  Then we can extend it by
the Hahn--Banach theorem to $\mu \in S(A)$.  Then the above argument works
again.\qed

\bigskip

It is shown in section 3 of \cite{R4} that if $G$ happens to be a
connected Lie group, then we can use norms on the Lie algebra,
${\mathfrak g}$, of $G$ to construct Lip-norms on $A$. We now consider
how this construction fits into the present context.

Let $\| \cdot\|_{\mathfrak g}$ be 
some norm on ${\mathfrak g}$.  As in \cite{R4} we
now take $A$ to be the dense $*$-subalgebra of ${\bar A}$
consisting of the elements which are once-differentiable for the action
$\a$.  For $X \in {\mathfrak g}$ we let $\a_X$ denote the derivation of
$A$ corresponding to $X$.   For any $a \in A$ we
denote by $da$ the linear map $X \mapsto \a_X(a)$ from ${\mathfrak g}$ to
${\bar A}$.  Since both ${\mathfrak g}$ and ${\bar A}$ have a
norm, the operator norm $\|da\|$ is defined.  Then for $a \in A$ we
set $L(a) = \|da\|$.  It is shown in
\cite{R4} that $L$ is a Lip-norm.  (But $L$ is now not closed,
though it is lower semi-continuous.)  

By using $\| \cdot\|_{\mathfrak g}$ in the usual way to define
the lengths of smooth paths in $G$, we obtain a continuous length function,
$\ell$, on $G$, to which we can apply the discussion of the present
section. But we need to require that $\| \cdot\|_{\mathfrak g}$ is
Ad-invariant in order that $\ell$ satisfy the extra condition
$\ell(yxy^{-1}) = \ell(x)$ needed above. In \cite{R4} the argument is recalled 
for the fact that $\|\a_x(a) -a\| \leq \|da\|\ell(x)$,
so that $L_0(a) \leq L(a)$, where $L_0$ is the Lip-norm from $\ell$.
But for any $X \in {\mathfrak g}$ we 
have $\ell(\exp(tX)) = |t|\|X\|_{\mathfrak g}$ for small $t$, and so if
$\|X\|_{\mathfrak g} = 1$ we have
$$
     \|(\a_{\exp(tX)}(a) -a)/t\| =  \|(\a_{\exp(tX)}(a) -a)\|/\ell(\exp(tX)).
$$
From this we see that $\|d_X a\| \leq L_0(a)$
for all $X$ with $\|X\|_{\mathfrak g} = 1$. Thus, for later use, we
can record:

\bigskip
\noindent
{\bf 8.6 PROPOSITION.} {\em For $L$ defined as above by a norm
on ${\mathfrak g}$, and for $L_0$ defined in terms of the corresponding
length function on $G$, we have  $L_0 = L$.}

\bigskip
Consequently Theorem 8.2 applies to $L$ defined by a norm on ${\mathfrak g}$.

In section 4 of \cite{R4} it is shown how to construct Lip-norms
on $A$ in terms of an inner product on the dual, ${\mathfrak g}'$, of
${\mathfrak g}$ and of a corresponding Dirac operator. Presumably
the results of the present section can be applied to that situation too,
as long as the inner product is invariant under the coadjoint 
representation, so that $G$ again acts by isometries for the Lip-norm.
But I have not verified this.

\section{Quantum tori}

We can now treat the main example of this paper.  We first recall the basic
facts about the quantum tori which we need, and establish our notation.  The
positive integer $d$ will denote the dimension of the tori we work with.  We
let $\Theta$ denote the vector space of all real skew-symmetric $d \x d$
matrices.  For any $\th \in \Theta$ we let ${\bar A}_{\th}$ denote the
corresponding quantum torus \cite{R1}, \cite{R3}.  It is defined as
follows.  Let $\s_{\th}$ denote the skew bicharacter on ${\mathbb Z}^d$
defined by
\[
\s_{\th}(p,q) = \exp(i\pi p \cdot \th q).
\]
Equip $C_c({\mathbb Z}^d)$, the space of complex-valued functions on
${\mathbb Z}^d$ of finite support, with the product consisting of
convolution twisted by $\s_{\th}$.  That is, for $f,g \in C_c({\mathbb
Z}^d)$ we have
$$
(f*g)(p) = \Sigma f(q)g(p-q)\s_{\th}(q,p). \leqno(9.1)
$$
We also equip
$C_c({\mathbb Z}^d)$ with the involution $f^*(p) = {\bar f}(-p)$, and the
norm of $\ell^1({\mathbb Z}^d)$, so that $C_c({\mathbb Z}^d)$ is a
$*$-normed algebra.  It is important for us that the involution and norm do
not depend on $\th$.  We let $\pi_{\th}$ denote the $*$-representation of
$C_c({\mathbb Z}^d)$ on the Hilbert space $\ell^2({\mathbb Z}^d)$ given by
formula $(9.1)$ except with $g$ replaced by $\xi \in \ell^2({\mathbb
Z}^d)$.  We let $\|\cdot\|_{\th}$ be the $C^*$-norm on $C_c({\mathbb Z}^d)$
defined by $\|f\|_{\th} = \|\pi_{\th}(f)\|$.  Then ${\bar A}_{\th}$ is
defined to be the completion of $C_c({\mathbb Z}^d)$ for this norm.

We let $G = {\mathbb T}^d$, where ${\mathbb T}$ is the circle group.  Thus
$G$ is the dual group of ${\mathbb Z}^d$.  We denote the duality by
$\<p,x\>$.  For us, a crucial fact is that $G$ has a natural action, $\a$,
on ${\bar A}_{\th}$, defined on $C_c({\mathbb Z}^d)$ by
\[
(\a_x(f))(p) = \<p,x\>f(p).
\]
Furthermore, this action is ergodic, so that we are in the setting of the
previous section.  It is important for us that the definition of the action
$\a$ on functions is independent of $\th$.

In \cite{R4} we gave three different methods for defining Lip-norms
on quantum tori, namely, by length functions on the dual group,
by norms on the Lie algebra of the dual group, and by a Dirac operator
construction. The simplest technically is in terms of length functions.
We discuss this case first.

We fix a continuous length function, $\ell$, on $G$.  (Because $G$ is
Abelian, the condition $\ell(xyx^{-1}) = \ell(y)$ is automatic.)  For each
$\th \in \Theta$ we define the Lip-norm $L_{\th}$ on ${\bar A}_{\th}$ in
terms of $\a$ and $\ell$ exactly as we did early in Section~8.  We let
$A_{\th}$ denote the dense $*$-subalgebra where $L_{\th}$ is finite.  Thus
each $(A_{\th},L_{\th})$ is a compact quantum metric space (at least if we
take its self-adjoint part), and the $L_{\th}$'s are defined consistently
in terms of one fixed length function $\ell$ as $\th$ varies.

We will actually obtain uniform estimates of quantum Gromov--Hausdorff
distance between the $A_{\th}$'s in terms of the parameter $\th$.  For this
purpose we fix an arbitrary norm on the vector space $\Theta$.  We can now
state our main theorem.

\bigskip
\noindent
{\bf 9.2 THEOREM.} {\em Let notation be as above.  For any $\e > 0$ there is
a $\d > 0$ such that if $\|\th - \psi\| < \d$ for $\th,\psi \in \Theta$,
then}
\[
\dist_q(A_{\th},A_{\psi}) < \e.
\]

\bigskip
\noindent
{\it\bfseries Proof.} Choose a finite-dimensional faithful unitary
representation, $\pi$, of $G$ as done in the previous section, so that, in
particular, its character, $\chi$, is non-negative.  Notice that the dual
group, ${\hat G}$, of $G$ is just ${\mathbb Z}^d$.  Thus, each ${\hat G}_n$,
as defined in the previous section, is just a finite subset of ${\mathbb
Z}^d$, containing $0$, closed under $p\rightarrow -p$, and generating
${\mathbb Z}^d$.  For each $p \in {\hat G}$ the $\a$-isotypic component of
$A_{\th}$ for $p$ clearly consists of the $1$-dimensional subspace of
functions on
${\mathbb Z}^d$ supported at $p$.  Thus the spaces $B_n$ of the previous
section all consist exactly of the functions supported on ${\hat G}_n
\subset {\mathbb Z}^d$, independently of $\th$.  What depends on $\th$ is
how these functions act as operators on $\ell^2({\mathbb Z}^d)$, and so
their operator norms and whether they are positive (but not whether they are
self-adjoint).  We denote these spaces with that structure by $B_{\th}^n$.

For each $n$ define $\var_n \in L^1(G)$ in terms of $\pi$ as done in the
previous section, and then define $\d_n$ in terms of $\var_n$ and $\ell$ as
done there.  Choose $n$ such that $\d_n < \e/3$.  We hold $n$ fixed for the
rest of the proof.  Then for simplicity of notation we write $B$ for $B_n =
L^1({\hat G}_n)$, and we write $B_{\th}$ instead of $B_{\th}^n$.  Each
$B_{\th}$ is equipped with the restriction to it of $\|\ \|_{\th}$ and
$L_{\th}$.  We know from Theorem $8.2$ that for each $\th$
\[
\dist_q(A_{\th},B_{\th}) < \e/3.
\]
  From the triangle inequality of Theorem $4.3$ it follows that, in order to
complete the proof, it suffices to find a $\d > 0$ such that if $\|\th -
\psi\| < \d$ then
\[
\dist_q(B_{\th},B_{\psi}) \le \e/3.
\]
We are in a favorable situation for doing this because each $B_{\th}$ is the
same finite-dimensional space $B$ of functions, but with different operator
norms and Lip-norms (but same self-adjoint part and order-unit).

We need to construct a bridge, $N$, between $B_{\th}$ and $B_{\psi}$.  For
this purpose we let $\pi_{\th}$ denote the representation of $A_{\th}$ on
$\ell^2({\mathbb Z}^d)$ defined earlier, restricted to $B_{\th}$; and
similarly for $\pi_{\psi}$.  We now collect the facts which we need.  What
we will be dealing with is the subject of continuous fields of compact
quantum metric spaces.  In the interest of brevity we will not develop its
general theory here.  In particular, we will not strive for minimal
hypotheses.  But we remark that much of the theory of continuous fields of
$C^*$-algebras, as presented, for example, in \cite{Bl}, \cite{D}, extends
to a theory of continuous fields of order-unit spaces.  One must then add
to that the topic of continuous fields of Lip-norms.  Many of our maneuvers
below can be placed in this general framework.

The first fact which we need is that the $C^*$-algebras ${\bar A}_{\th}$
form a continuous field of $C^*$-algebras over $\Theta$.
Specifically, we need the fact that for any $f \in \ell^1({\mathbb Z}^d)$
the function $\th
\rightarrow \|f\|_{\th}$ is continuous.  See corollary $2.8$ of \cite{R2};
but note from the proof of theorem $2.5$ of \cite{R2} that while $\th
\rightarrow \pi_{\th}(f)$ is strongly continuous, it cannot be expected to
be norm continuous.  (Also, the hypotheses of theorem $2.5$ should have
included the requirement that $\Omega$ be first countable.)  We denote by
$\|\cdot\|_*$ (rather than $\|\cdot\|_1$) the usual norm on $\ell^1({\mathbb
Z}^d)$.  It is a standard fact \cite{FD} that $\|f\|_{\th} \le \|f\|_*$ for
$f \in \ell^1({\mathbb Z}^d)$.

Let $S$ denote some finite subset of ${\mathbb Z}^d$, such as our earlier
${\hat G}_n$'s.  Note that $\a$ carries the finite-dimensional vector space
$\ell^1(S)$ into itself.  By finite-dimensionality we have $L_{\th}(f) <
\infty$ for all $f \in \ell^1(S)$.  (See the proof of proposition $2.2$ of
\cite{R4}.)  Of course, each $f \in \ell^1(S)$ is ``Lipschitz'' for $\|\
\|_*$ also.  We denote the corresponding Lip-norm by $L_*$.

\bigskip
\noindent
{\bf 9.3 LEMMA.} {\em Let $f \in \ell^1({\mathbb Z}^d)$.  If $f$ has finite
support, then the function $\th \rightarrow L_{\th}(f)$ is continuous.
Furthermore, $L_{\th}(f) \le L_*(f)$ for each $\th \in \Theta$.}

\bigskip
\noindent
{\it\bfseries Proof.} Let $f$ be supported on the finite set $S$.  Let
\[
D_f = \{(\a_x(f) - f)/\ell(x): x \ne e_G\}.
\]
Note that $D_f \subseteq \ell^1(S)$, and that $D_f$ is bounded for $\|\
\|_*$ since $L_*(f) < \infty$.  Now $L_{\th}(f)$ is the supremum of
$\|g\|_{\th}$ for $g \in D_f$, and similarly for $L_*(f)$.  Thus
$L_{\th}(f) \le L_*(f)$, since $\|\cdot \|_{\th} \le \|\cdot \|_*$.

For any $g \in \ell^1(S)$ let $F_g$ be the 
continuous function on $\Theta$ defined by
$F_g(\th) = \|g\|_{\th}$.  Since $L_{\th}(f) = \sup\{\|g\|_{\th}: g \in
D_f\}$, it follows that $\th \rightarrow L_{\th}(f)$ is the supremum of the
functions $F_g$ for $g \in D_f$.  As a supremum of continuous functions,
$\th \mapsto L_{\th}(f)$ must be lower semi-continuous.

We now show that $\th \mapsto L_{\th}(f)$ is
actually continuous.  For any $h,k \in\ell^1(S)$ we have
\[
|F_h(\th) - F_k(\th)| = |\|h\|_{\th} - \|k\|_{\th}| \le \|h-k\|_{\th} \le
\|h-k\|_*.
\]
Let $U$ be any compact subset of $\Theta$, and view $F_h$ and
$F_k$ as restricted to $U$, and so as elements of $C(U)$ with the usual
norm $\|\cdot\|_{\infty}$.  The above inequality says that $\|F_h -
F_k\|_{\infty} \le \|h-k\|_*$.  Thus $F$, as a mapping from $\ell^1(S)$
into $C(U)$, is Lipschitz.  Now $D_f$, as a bounded subset of a
finite-dimensional normed vector space, is totally bounded.  Consequently,
$F(D_f)$ is totally bounded in $C(U)$.  But by part of the Arzela--Ascoli
theorem \cite{Cw} it follows that $F(D_f)$ is equicontinuous.  In other
words, $\{F_g: g \in D_f\}$ is an equicontinuous family of functions on
$U$.  But the supremum of an equicontinuous family is continuous.  Thus
$\th \mapsto L_{\th}(f)$ is continuous on $U$.\qed

\bigskip
For each $\th$ we let $r_{\th}$ denote the radius of $(A_{\th},L_{\th})$.
As an immediate consequence of lemma $2.4$ of \cite{R4} we have:

\bigskip
\noindent
{\bf 9.4 LEMMA.} {\em There is a constant, $R$, which depends only on the
length function $\ell$ and the choice of Haar measure on $G$, such that
$r_{\th} \le {\mathbb R}$ for all $\th \in \Theta$.}

\bigskip
We now need an ample supply of continuous fields of states.  We let
${\mathcal T}_1^+$ denote the space of positive trace-class operators of
trace $1$ on $\ell^2({\mathbb Z}^d)$, that is, the ``density matrices''.  We
denote the trace by tr.  For each $T \in {\mathcal T}_1^+$ and 
each $\th \in
\Theta$ we let $\omega_{\th}^T$ denote the state on ${\bar A}_{\th}$ defined
by
\[
\omega_{\th}^T(a) = \mbox{tr}(\pi_{\th}(a)T).
\]
Because $\pi_{\th}$ is a faithful representation of ${\bar A}_{\th}$, it
follows from proposition VII$.5.4$ of \cite{FD} that this set of states is
$w^*$-dense in $S(A_{\th})$.

\bigskip
\noindent
{\bf 9.5 LEMMA.} {\em For each $f \in \ell^1({\mathbb Z}^d)$ and $T \in
{\mathcal T}_1^+$ the function $\th \mapsto \omega_{\th}^T(f)$ is
continuous.}

\bigskip
\noindent
{\it\bfseries Proof.} As mentioned above (with caveat), the function $\th
\mapsto \pi_{\th}(f)$ is strong-operator continuous by the proof of theorem
$2.5$ of \cite{R2}.  It follows easily that $\th \mapsto \omega_{\th}^T(f)$
is continuous when $T$ has rank~$1$; and hence when $T$ is of finite
rank.  But the finite rank operators in ${\mathcal T}_1^+$ are dense for the
trace-norm.  A uniform convergence argument completes the proof.\qed

\bigskip
We have now collected the facts which we need.  Before continuing with the
proof, we turn to considering briefly the 
method discussed at the end of Section 8 
for defining Lip-norms in terms of norms on the Lie algebra.
Accordingly, let ${\mathfrak g}$ denote the Lie algebra of ${\mathbb T}^d$,
so that ${\mathfrak g}$ can be identified with ${\mathbb R}^d$.  
Let $\|\cdot\|_{\mathfrak g}$ 
denote some norm on ${\mathfrak g}$.  As done
near the end of Section 8, we
now take $A_{\th}$ to be the dense $*$-subalgebra of ${\bar A}_{\th}$
consisting of the elements which are once-differentiable for the action
$\a$.  For $X \in {\mathfrak g}$ we let $\a_X$ denote the derivation of
$A_{\th}$ corresponding to $X$.  Thus for $f \in C_c({\mathbb Z}^d)$ and $p
\in {\mathbb Z}^d$ we have
\[
(\a_X(f))(p) = 2\pi i(p \cdot X)f(p),
\]
where $p \cdot X$ denotes the standard inner product on ${\mathbb R}^d$ and
we use ${\mathbb Z}^d \subseteq {\mathbb R}^d$. We define $da$ as 
near the end of Section 8, and
set $L_{\th}(a) = \|da\|$.  In the same way we define $L_*$.   
The proof of Lemma $9.3$ for this setting is actually easier than the proof 
given above:

\bigskip
\noindent
{\bf 9.6 LEMMA.} {\em Let $A_{\th}$ and $L_{\th}$ be defined as just above
in terms of a norm on ${\mathfrak g}$.  For each $f \in C_c({\mathbb Z}^d)$
the function $\th \mapsto L_{\th}(f)$ is continuous.  Furthermore,
$L_{\th}(f) \le L_*(f)$ for each $\th \in \Theta$.} (This is also
true if $f$ is a Schwartz function, or if $p \mapsto \|p\|f(p)$
is in $\ell^1({\mathbb Z}^d)$.)

\bigskip
\noindent
{\it\bfseries Proof.} Changing the notation given in the proof of Lemma
$9.3$, we now set
\[
D_f =\{\a_X(f): \|X\|_{\mathfrak g} \le 1\},
\]
so that $L_{\th}(f) = \sup\{\|g\|_{\th}: g \in D_a\}$.  From this we
immediately have $L_{\th}(f) \le L_*(f)$.  Define $F_g$ exactly as in the
proof of Lemma $9.3$, so that again $\th \mapsto L_{\th}(f)$ is the supremum
of the functions $F_g$ for $g \in D_f$.  Since $\{X: \|X\|_{\mathfrak g} \le
1\}$ is a compact subset of a finite-dimensional vector space, and $X
\mapsto \a_X(f)$ is linear, it is clear that $D_f$ is a compact subset of
$\ell^1({\mathbb Z}^d)$.  Then the argument in the proof of Lemma $9.3$
shows that $F(D_f)$ is an equicontinuous family of functions on compact
subsets of $\Theta$, so that $\th \mapsto L_{\th}(f)$ is continuous.\qed

\bigskip

From the discussion at the end of Section 8 it is clear 
that the analogue of Lemma $9.4$ holds.

\bigskip
In order to try to make the rest of 
the proof of Theorem 9.2 as clear as possible, we 
find it useful to treat the situation developed in this section
in a somewhat axiomatic framework. We carry this out in
the next two sections.  But we remark first that
because the $(A_{\th},L_{\th})$'s are unchanged when integers are added to
the entries of $\th$, it is sufficient to prove that our main theorem holds
for any compact subset, $\Theta_0$, of $\Theta$.  For simplicity of
notation, in the next sections we will let $\Theta$ denote a compact metric
space.

\section{Continuous fields of order-unit spaces}

As said earlier, we will not strive for a full theory of continuous fields
here.  Rather we will take a fairly direct route to what we need.
Throughout this section $V$ will be a finite-dimensional real vector space
equipped with a distinguished element $e$ (for example, the self-adjoint
part of $\ell^1({\hat G}_n)$ with $e = 1 = \d_0$).  We let $\Theta$ be a
compact set with metric $d$.  We assume that we are given, for each $\th \in
\Theta$, a norm, $\|\cdot\|_{\th}$, on $V$ such that for each $v \in V$ the
function $\th \mapsto \|v\|_{\th}$ is continuous.  We will call such a
family of norms a {\em continuous field of norms}.  We will shortly make
further hypotheses on the $\|\cdot\|_{\th}$'s.

If we pick a basis $\{v_j\}$ for $V$ and use compactness and continuity to
find a constant $K$ such that $\|v_j\|_{\th} \le K$ for all $\th \in
\Theta$ and all $j$, we see quickly that there is a norm, $\|\cdot\|_*$,
such that $\|\cdot\|_{\th} \le \|\cdot\|_*$ for all $\th$.  It will be
convenient for us to fix such a norm.

We let $V'$ denote the vector-space dual of $V$, and we let
$\|\cdot\|'_{\th}$ and $\|\cdot\|'_*$ denote the corresponding dual
norms.  We let $\Sigma_{\th}$ denote the unit $\|\cdot\|_{\th}$-sphere
(not ball), and similarly for $\Sigma_*$, etc.

\bigskip
\noindent
{\bf 10.1 LEMMA.} {\em The function $(\th,v) \mapsto \|v\|_{\th}$ is
jointly continuous on $\Theta \x V$.  The family $\{\|\cdot\|'_{\th}\}$ is
a continuous field of norms on $V'$.  There is a strictly positive
constant, $k$, such that for all $\th \in \Theta$}
\[
k\|\cdot\|_* \le \|\cdot\|_{\th} \le \|\cdot\|_*
\]
{\em and}
\[
k^{-1}\|\cdot\|'_* \ge \|\cdot\|'_{\th} \ge \|\cdot\|'_*.
\]

\bigskip
\noindent
{\it\bfseries Proof.} The joint continuity at $(w,\psi)$ follows from the
inequalities
\begin{eqnarray*}
|\|v\|_{\th} &- &\|w\|_{\psi}| \le |\|v\|_{\th} - \|w\|_{\th}| +
|\|w\|_{\th} -  \|w\|_{\psi}| \\
&\le &\|v-w\|_{\th} + |\|w\|_{\th} - \|w\|_{\psi}| \le \|v-w\|_* +
|\|w\|_{\th} - \|w\|_{\psi}|.
\end{eqnarray*}
Let $\l \in V'$, and for each $v \in \Sigma_*$ define a function,
$F_{v}$, on
$\Theta$ by
\[
F_{v}(\th) = |\<v,\l\>|/\|v\|_{\th}.
\]
  From the above joint continuity we see that the function $(\th,v) \mapsto
F_{v}(\th)$ is jointly continuous on $\Theta \x \Sigma_*$.  It follows
from the compactness of $\Theta \x \Sigma_*$ that this function is
uniformly continuous.  But this implies that the family $\{F_{v}\}$ is
equicontinuous.  Thus, $\th \mapsto \|\l\|'_{\th}$, being the supremum of
this equicontinuous family, is continuous.

The function $(\th,v) \mapsto \|v\|_{\th}$ is continuous on $\Theta \x
\Sigma_*$ and never takes value $0$ there.  Thus, by compactness there is a
strictly positive constant, $k$, such that $k \le \|v\|_{\th}$ for all
$\th$ and all $v$ with $\|v\|_* = 1$, as needed.  The inequalities for
the dual norms follow.\qed

\bigskip
The metric on $V$ from the norm $\|\cdot\|_*$ give a corresponding Hausdorff
metric on the compact subsets of $V$, which we denote by $\dist_H^*$.  We
will not explicitly need the following lemma and its corollary later.  But
we will need the calculation which constitutes most of the proof of the
lemma; and the corollary is a good preview of what comes later.

\bigskip
\noindent
{\bf 10.2 LEMMA.} {\em Let $\|\cdot\|_1$ and $\|\cdot\|_2$ be norms on $V$,
and let $k$ be a constant such that $k\|v\|_* \le \|v\|_j$ for $j = 1,2$
and all $v \in V$.  Let $\e > 0$ be given, and suppose that $|\|v\|_1 -
\|v\|_2| < \e k^2\|v\|_*$ for all $v$.  Then
$\dist_H^*(\Sigma_1,\Sigma_2) < \e$.}

\bigskip
\noindent
{\it\bfseries Proof.} Let $u \in \Sigma_1$ be given.  Set $v = u/\|u\|_2$,
so that $v \in \Sigma_2$. Then
\begin{eqnarray*}
\|u-v\|_* &= &\|(1-\|u\|_2^{-1})u\|_* \\
&= &|\|u\|_2 - 1|\|u\|_2^{-1}\|u\|_* \le |\|u\|_2 - \|u\|_1|k^{-1} \\
&\le &\e k\|u\|_* \le \e\|u\|_1 = \e.
\end{eqnarray*}
Thus, $\Sigma_1$ is in the $\e$-neighborhood of $\Sigma_2$ for
$\|\cdot\|_*$.  But we can reverse the roles of $\Sigma_1$ and
$\Sigma_2$.\qed

\bigskip
\noindent
{\bf 10.3 COROLLARY.} {\em For a continuous field $\{\|\cdot\|_{\th}\}$ of
norm as above, the function $\th \mapsto \Sigma_{\th}$ is (uniformly)
continuous for $\dist_H^*$.}

\bigskip
\noindent
{\it\bfseries Proof.} Let $\e > 0$ be given.  Let $k$ be as in Lemma $10.1$,
and let $E_*$ denote the unit $\|\cdot\|_*$-ball of $V$.  According to Lemma
$10.1$ the function $(v,\th) \mapsto \|v\|_{\th}$ restricted to the compact
set $\Theta \x E_*$ is jointly continuous, and so uniformly continuous.
Consequently, the family of functions $\th \mapsto \|v\|_{\th}$ for $v \in
E_*$ is uniformly equicontinuous, that is, we can find $\d > 0$ such that if
$d(\th,\psi) < \d$, then for each $v \in E_*$ we have $|\|v\|_{\th} -
\|v\|_{\psi}| \le \e k^2$.  We are then exactly in position to use Lemma
$10.2$ to conclude that $\dist_H^*(\Sigma_{\th},\Sigma_{\psi}) < \e$.\qed

\bigskip
In accordance with the situation for $B$ and the norms $\|\cdot\|_{\th}$ of
the previous section, we will now assume further that each $\|\cdot\|_{\th}$
is an order-unit norm on $(V,e)$, where:

\bigskip
\noindent
{\bf 10.4 DEFINITION.} Let $V$ be a vector space and let $e$ be a
distinguished element of $V$.  Let $\|\cdot\|$ be a norm on $V$, and set
\[
V^+ = \{v \in V: \|(\|v\|e - v)\| \le \|v\|\}.
\]
We say that $\|\cdot\|$ is an {\em order-unit norm} on $(V,e)$ if $\|e\| =
1$ and $V^+$ is the positive cone for an ordering on $V$ for which $(V,e)$
is an order-unit space with order-unit norm equal to $\|\cdot\|$.  When
$\|\cdot\|$ is an order-unit norm, we will feel free to refer to
$(V,e,\|\cdot\|)$ as an order-unit space, with the ordering from $V^+$
understood.

\bigskip
We remark that the above definition of $V^+$ is motivated, for example, by
corollary VI$.7.7$ of \cite{FD}, or by part~2 of the proof of proposition
II$.1.3$ of \cite{A}.

Since the norms $\|\cdot\|_{\th}$ are now order-unit norms, it is reasonable
to consider the corresponding state spaces.  They will all consist of linear
functionals on our fixed vector space $V$.  We denote the state space for
$\|\cdot\|_{\th}$ by $S_{\th}$.  Our eventual aim is, roughly speaking, to
show that the $\e$-density of a family of continuous fields of states
propagates.  This will permit us to define suitable bridges $N$.  For this
purpose we need to examine the continuous-field structure of the duals of
the order-unit spaces $(V,e,\|\cdot\|_{\th})$.  By a theorem of Ellis
(theorem II$.1.15$ of \cite{A}) the dual of an order-unit space is a
base-norm space.  This means the following.  Let $(V,e,\|\cdot\|)$ be an
order-unit space.  Let $V'$ be its Banach-space dual.  Let $\eta$ denote the
linear functional which $e$ defines on $V'$, and let $\|\cdot\|'$ denote the
dual norm.  The {\em base}, $S$, is defined by
\[
S = \{\l \in V': \eta(\l) = 1 = \|\l\|'\}.
\]
We recognize this as just the state space.  Then the unit $\|\cdot\|'$-ball
of
$V'$ is equal to
$\mbox{co}(S
\cup -S)$.  (See proposition II$.1.7$ of \cite{A}.)  We emphasize that
because
$S$ is
$w^*$-compact, ``co'' here means ``convex hull'', not ``closed convex
hull''.  It is easy to see then that if $\l \in V'$ with $\|\l\|' \le 1$,
then there are $\mu,\nu \in S$ (not unique) and $t \in [0,1]$ such that $\l
= t\mu - (1-t)\nu$.  This is closely related to lemma $2.1$ of \cite{R5}
mentioned in the last paragraph of Section~1.  Let $C$ be the cone in $V'$
generated by
$S$ (so that
$S$ is a ``base'' for
$C$ because $\eta(S) = 1$).  Then $C$ is the positive cone for an order on
$V'$ which satisfies the requirements for $(V',C,S)$ to be a base-norm
space, as defined just after proposition II$.1.12$ of \cite{A}.

All this suggests the following:

\bigskip
\noindent
{\bf 10.5 DEFINITION.} Let $X$ be a vector space, and let $\eta$ be a
distinguished linear functional on $X$.  By a {\em base-norm} on $(X,\eta)$
we mean a norm, $\|\cdot\|$, on $X$ such that when we set
\[
S = \{\l \in X: \eta(\l) = 1 = \|\l\|\},
\]
the unit $\|\cdot\|$-ball coincides with $\mbox{co}(S \cup -S)$.

\bigskip
Much as above, it is easily seen that if $\|\cdot\|$ is a base-norm on
$(X,\eta)$, then $(X,S)$ is a base-norm space \cite{A} when $X$ is ordered
by the cone generated by $S$.

Returning to the situation where we have a continuous field of order-unit
norms on a finite-dimensional vector space $V$, we see that
$\{\|\cdot\|'_{\th}\}$ is a continuous field of base-norm on $(V',\eta)$.
We can now momentarily forget where $V'$ and the $\|\cdot\|'_{\th}$ come
from, and just consider a finite-dimensional vector space, $X$, with linear
functional $\eta$, and a continuous family $\{\|\cdot\|_{\th}\}$ of
base-norms.  In view of Lemma $10.1$, we now assume that we have a norm
$\|\cdot\|_*$ and a constant $k$ such that
\[
k^{-1}\|x\|_* \ge \|x\|_{\th} \ge \|x\|_*
\]
for $x \in X$.  The metric from $\|\cdot\|_*$ defines a Hausdorff metric on
the compact subsets of $X$, which we again denote by $\dist_H^*$.

\bigskip
\noindent
{\bf 10.6 THEOREM.} {\em Let $X$ be a finite-dimensional vector space with
distinguished linear functional $\eta$.  Let $\{\|\cdot\|_{\th}\}$ be a
continuous field of base-norms for $(X,\eta)$, and let $S_{\th}$ denote the
base for $(V,\eta,\|\cdot\|_{\th})$.  Let $\|\cdot\|_*$ be a norm on $X$
such that $\|\cdot\|_{\th} \ge \|\ \|_*$ for all $\th$.  Then the function
$\th \mapsto S_{\th}$ on $\Theta$ is (uniformly) continuous for $\dist_H^*$.}

\bigskip
\noindent
{\it\bfseries Proof.} It suffices to show that for every $\e > 0$ such that
$\e < 2$ there is a $\d > 0$ such that if $d(\th,\psi) < \d$ then for every
$\mu \in S_{\psi}$ there is a $\nu \in S_{\th}$ with $\|\mu-\nu\|_*  < \e$,
since then we can reverse the roles of $\th$ and $\psi$.  For each $\th$ let
$\Sigma_{\th}$ denote the unit $\|\cdot\|_{\th}$-sphere of $X$, much as we
did earlier in the dual situation.  Much as in the proof of Corollary
$10.3$, we let $E_*$ denote the unit $\|\cdot\|_*$-ball in $X$, so that
$\Sigma_{\th} \subset E_*$ for all $\th$.  Note that $S_{\th} \subset
\Sigma_{\th}$.

Given $\e > 0$, we choose $\d$ as follows.  The function $(\th,x) \mapsto
\|x\|_{\th}$ restricted to the compact set $\Theta \x E_*$ is jointly
continuous by Lemma $10.1$, and so is uniformly continuous.  Thus, the
family of functions $\th \mapsto \|x\|_{\th}$ for $x \in E_*$ is uniformly
equicontinuous, and so we can find $\d > 0$ such that if $d(\th,\psi) < \d$
then for each $x \in E_*$ we have
\[
|\|x\|_{\th} - \|x\|_{\psi}| < \e/4.
\]
This choice of $\d$ works, as we see from the following lemma, which we will
also need later.

\bigskip
\noindent
{\bf 10.7 KEY LEMMA.} {\em Let $\|\cdot\|_1$ and $\|\cdot\|_2$ be base-norms
on $(X,\eta)$, and let $\|\cdot\|_*$ be a norm on $X$ such that $\|x\|_* \le
\|x\|_j$ for $j= 1,2$ and all $x \in X$.  Let $\e > 0$ be given such that
$\e < 2$, and suppose that
\[
|\|x\|_1 - \|x\|_2| < (\e/4)\|x\|_*
\]
for all $x \in X$.  Then $\dist_H^*(S_1,S_2) < \e$.}

\bigskip
\noindent
{\it\bfseries Proof.} Let $\mu \in S_1$.  Let $x = \mu/\|\mu\|_2$, so that
$x \in \Sigma_2$.  By the calculation in the proof of Lemma $10.2$ we see
that $\|\mu-x\|_* < \e/4$.  Since $\|\cdot\|_2$ is a base-norm and $\|x\|_2
= 1$, we can find $\nu,\nu_0 \in S_2$ such that $x = t\nu - (1-t)\nu_0$ for
some $t \in [0,1]$.  When we apply $\eta$ to this equation we obtain
$\|\mu\|_2^{-1} = 2t - 1$.  Thus,
\[
\|\mu\|_2^{-1} - 1 = 2(t-1).
\]
But $\e < 2$ and $\|\mu\|_* \le \|\mu\|_1 = 1$, and so
\[
|1 - \|\mu\|_2| = |\|\mu\|_1 - \|\mu\|_2| < (\e/4)\|\mu\|_* \le 1/2.
\]
Thus $\|\mu\|_2^{-1} < 2$, so that, because $1 = \|\mu\|_1$,
\[
|\|\mu\|_2^{-1} - 1| = |\|\mu\|_1 - \|\mu\|_2|\|\mu\|_2^{-1} <
(\e/2)\|\mu\|_* \le \e/2.
\]
Consequently, $2|t-1| < \e/2$, and so
\[
\|x-\nu\|_* = \|(t-1)(\nu+\nu_0)\|_* \le |t-1|\|\nu+\nu_0\|_2 \le 2|t-1| <
\e/2.
\]
It follows that $\|\mu-\nu\|_* < \e$.  Thus $S_1$ is in the
$\e$-neighborhood of $S_2$ for $\dist_H^*$.  But we can reverse the roles of
$S_1$ and $S_2$.\qed

\bigskip
If the requirement that the $\|\cdot\|_j$'s be base-norms is dropped, the
conclusions of the above lemma and theorem can easily fail.  An example
could involve strictly convex norms converging to a base-norm.

In our earlier concrete situation of $\ell^1({\mathbb Z}^d)$ we saw that we
had a large supply of continuous fields of states, namely the
$\omega_{\th}^T$'s.  For each $\th$ they were $w^*$-dense in $S(A_{\th})$.
But when restricted to $B = \ell^1({\hat G}_n)$ they will then fill out all
of $S(B_{\th})$.  Thus, in our present abstract setting of a
finite-dimensional vector space $V$ with distinguished element $e$ and
order-unit norms $\|\cdot\|_{\th}$ we can assume that we have a large
family, ${\mathcal S}$, of continuous fields of states, where:

\bigskip
\noindent
{\bf 10.8 DEFINITION.} With notation as above, by a {\em continuous field of
states} we will mean a function, $\Omega$, from $\Theta$ to $V'$, such that
$\Omega_{\th} \in S_{\th}$ for each $\th \in \Theta$, and such that $\th
\mapsto \<v,\Omega_{\th}\>$ is continuous for each $v \in V$.  When $V$ is
finite-dimensional, the latter condition is equivalent to $\Omega$ being
continuous for $\|\cdot\|'_*$ on $V'$.

\bigskip
We remark that related definitions of continuous fields of states for
continuous fields of $C^*$-algebras are given in \cite{Bl}, \cite{L},
\cite{N}, \cite{N2}.  If one does not already know that one has enough continuous
fields of states, one can, in the infinite-dimensional case, try to prove
their existence by means of selection theorems, as done for theorem $3.3$ of
\cite{Bl}.  Although we do not need it for our application to quantum tori,
it seems appropriate to point out that in our finite-dimensional situation
there are always plenty of continuous fields of states.  We can again work
directly with fields of base-norms.  In that setting, by a ``continuous
field of states'' we will mean simply a continuous function $\Omega$ from
$\Theta$ to $X$ such that $\Omega_{\th} \in S_{\th}$ for each $\th \in
\Theta$.

\bigskip
\noindent
{\bf 10.9 PROPOSITION.} {\em Let $X$ be finite-dimensional, and let
$\{\|\cdot\|_{\th}\}$ be a continuous field of base-norms on $(X,\eta)$.
Then for each $\psi \in \Theta$ and each $\mu \in S_{\psi}$ there is a
continuous field of states, $\Omega$, such that $\Omega_{\psi} = \mu$.}

\bigskip
\noindent
{\it\bfseries Proof.} Because $X$ is finite-dimensional, we can choose a
norm, $\|\cdot\|$, on $X$ which comes from an inner-product, and so
satisfies the parallelogram law.  We denote the Hausdorff metric which it
defines on compact subsets of $X$ by $\dist_H^P$.  Again, by
finite-dimensionality, the norm $\|\cdot\|$ will be equivalent to
$\|\cdot\|_*$, and from this it is easily seen that the metric $\dist_H^P$
is equivalent to $\dist_H^*$.  From Theorem $10.6$ we conclude that the
function $\th \mapsto S_{\th}$ is continuous on $\Theta$ for $\dist_H^P$.

Let $\psi \in \Theta$ and $\mu \in S_{\psi}$ be given.  Each $S_{\th}$ is a
compact subset of $X$, and so we can find a point of $S_{\th}$ which is
closest to $\mu$ for $\|\cdot\|$.  We denote this point by $\Omega_{\th}$.
Clearly $\Omega_{\psi} = \mu$.  We wish to show that the function $\Omega$
is continuous.  But this is an immediate consequence of the following
proposition, which may well appear in the literature somewhere.

\bigskip
\noindent
{\bf 10.10 PROPOSITION.} {\em Let $X$ be a Hilbert space with norm
$\|\cdot\|$.  Let ${\mathcal C}$ denote the collection of compact convex
subsets of $X$, with corresponding Hausdorff metric $\dist_H^P$.  Fix $z \in
X$, and for each $C \in {\mathcal C}$ let $\Omega_C$ denote an element of
$C$ closest to $z$.  (This point is unique, but we don't explicitly need
that fact.)  Then the function $\Omega$ is continuous from ${\mathcal C}$ to
$X$ for $\dist_H^P$ and $\|\cdot\|$.}

\bigskip
\noindent
{\it\bfseries Proof.} By translation we can, and do, assume that $z = 0$.
Let $D \in {\mathcal C}$.  We show continuity of
$\Omega$ at $D$.  Let $m$ be a constant such that $\|x\| \le m$ for all $x
\in D$.  Let $\d > 0$ be given, and suppose that $C \in {\mathcal C}$ with
$\dist_H^P(C,D) < \d$.  Then there is a $c \in C$ with $\|c - \Omega_D\| <
\d$, and a $d \in D$ with $\|d - \Omega_C\| < \d$.  By the parallelogram law
\begin{eqnarray*}
2\|c\|^2 + 2\|\Omega_C\|^2 &= &\|c + \Omega_C \|^2 + \|c
-\Omega_C\|^2 \\
&= &4\|(c+\Omega_C)/2 \|^2 + \|c-\Omega_C\|^2.
\end{eqnarray*}
Now $(c+\Omega_C)/2 \in C$, and so $\|(c+\Omega_C)/2 \| \ge \|\Omega_C\|$,
while
\[
\|c\| \le \|\Omega_D\| + \d \le \|d\| + \d \le \|\Omega_C\| + 2\d.
\]
When we combine these two inequalities with the above application of the
parallelogram law, we obtain
\[
2(\|\Omega_C\| + 2\d)^2 + 2\|\Omega_C\|^2 \ge 4\|\Omega_C\|^2 +
\|c-\Omega_C\|^2.
\]
Simplifying, we obtain $8\d(\|\Omega_C\| + \d) \ge \|c-\Omega_C\|^2$.
But $\|\Omega_C \| \le m + \d $, 
and so 
\[
\|c - \Omega_C\|^2 \le 8\d(m + 2\d).
\]
Since
\[
\|\Omega_C-\Omega_D\| \le \|\Omega_C-c\| + \|c-\Omega_D\| \le \|c-\Omega_C\|
+ \d,
\]
it is now clear that by making $\d$ sufficiently small we can arrange that
$\|\Omega_C-\Omega_D\|$ be as small as desired.\qed

\bigskip
\noindent
{\bf 10.11 LEMMA.} {\em With notation as above, let $\psi \in \Theta$ and
$\e > 0$ be given.  Suppose that ${\mathcal S}$ is a family of continuous
fields of states such that $\{\Omega_{\psi}: \Omega \in {\mathcal S}\}$ is
$\e$-dense in $S_{\psi}$ for $\|\ \|_*$.  Then there is a $\d > 0$ such that
if $d(\th,\psi) < \d$ then $\{\Omega_{\th}: \Omega \in {\mathcal S}\}$ is
$3\e$-dense in $S_{\th}$ for $\|\cdot\|_*$.}

\bigskip
\noindent
{\it\bfseries Proof.} If ${\mathcal S}$ is not already finite, we replace it
by a finite subset, since $S_{\psi}$ is compact.  Each $\Omega \in {\mathcal
S}$ is continuous for $\|\cdot\|_*$, and since there are only a finite
number of them, we can find $\d_1$ such that if $d(\th,\psi) < \d_1$ then
for each $\Omega \in {\mathcal S}$ we have $\|\Omega_{\th} -
\Omega_{\psi}\|_* < \e$.  By Theorem $10.6$ we can find $\d \le \d_1$ such
that for $d(\th,\psi) < \d$ we have $\dist_H^*(S_{\th},S_{\psi}) < \e$.  For
such a $\th$ let $\mu \in S_{\th}$ be given.  Then, by this Hausdorff
metric closeness, there is a $\nu \in S_{\psi}$ such that $\|\mu-\nu\|_* <
\e$.  By the condition on ${\mathcal S}$ there is an $\Omega \in {\mathcal
S}$ such that $\|\nu-\Omega_{\psi}\|_* < \e$.  In view of the choice of
$\d_1$ we have $\|\Omega_{\th} - \Omega_{\psi}\|_* < \e$.  Thus,
$\|\Omega_{\th} - \mu\|_* < 3\e$.\qed

\bigskip
\noindent
{\bf 10.12 PROPOSITION.} {\em Let $\{\|\cdot\|_\th\}$ be a continuous field of
base-norms on $X$.  Let ${\mathcal S}$ be a family of continuous fields of
states such that $\{\Omega_{\th}: \Omega \in {\mathcal S}\} = S_{\th}$ for
each $\th \in \Theta$.  Given any $\e > 0$ there is a finite subset,
${\mathcal S}_{\e}$, of ${\mathcal S}$ such that $\{\Omega_{\th}: \Omega
\in {\mathcal S}_{\e}\}$ is $\e$-dense for $\|\cdot\|_*$ in $S_{\th}$ for
every $\th \in \Theta$.}

\bigskip
\noindent
{\it\bfseries Proof.} Let $\e > 0$ be given.  For each $\psi \in \Theta$ we
can choose a finite subset, ${\mathcal S}_{\psi}$, of ${\mathcal S}$ such
that
$\{\Omega_{\psi}: \Omega \in {\mathcal S}_{\psi}\}$ is $\e/3$-dense in
$S_{\psi}$ for $\|\cdot\|_*$.  By Lemma $10.11$ there is a neighborhood,
$U_{\psi}$, of $\psi$ such that for $\th \in U_{\psi}$ we have
$\{\Omega_{\th}: \Omega \in {\mathcal S}_{\th}\}$ is $\e$-dense in
$S_{\th}$.  By the compactness of $\Theta$ a finite number of the
$U_{\psi}$'s cover
$\Theta$.  The union of the corresponding ${\mathcal S}_{\psi}$'s
provides the desired ${\mathcal S}_{\e}$.\qed

\bigskip
In view  of Proposition $10.7$ and the equivalence of norms on
finite-dimensional vector spaces, we then obtain:

\bigskip
\noindent
{\bf 10.13 THEOREM.} {\em Let $X$ be finite-dimensional, and let
$\{\|\cdot\|_{\th}\}$ be a continuous field of base-norms on $(X,\eta)$.
Let
$\|\cdot\|_*$ be a norm on $X$.  For each $\e > 0$ there is a finite set,
${\mathcal S}$, of continuous fields of states such that $\{\Omega_{\th}:
\Omega \in {\mathcal S}\}$ is $\e$-dense in $S_{\th}$ for $\|\cdot\|_*$ for
each
$\theta \in \Theta$.}

\section{Continuous fields of Lip-norms}

We must now bring Lip-norms into the picture.  Throughout this section we
again assume that $V$ is finite-dimensional with norm $\|\cdot\|_*$, and
that $\th \mapsto \|\cdot\|_{\th}$ is a continuous field of order-unit norms
for $(V,e)$ such that $\|\cdot\|_{\theta} \le \|\cdot\|_*$ for all $\theta
\in \Theta$.  We will denote the order-unit space
$(V,e,\|\cdot\|_{\th})$ by $V_{\th}$.  Motivated by the considerations of
Section~9, especially Lemma $9.3$, we now assume, in addition, that for each
$\th \in \Theta$ we are given a Lip-norm $L_{\th}$ on $(V,e)$ such that $\th
\mapsto L_{\th}(v)$ is continuous for each $v \in V$.  We will refer to such
a family
$\{L_{\th}\}$ as a {\em continuous field of Lip-norms}.  Since $V$ is
finite-dimensional, ``Lip-norm'' here means just that its null-space is
spanned by $e$.

\bigskip
\noindent
{\bf 11.1 PROPOSITION.} {\em Let ${\mathcal S}$ be a finite, non-empty, set
of continuous fields of states for the $V_{\th}$'s.  Let $\e > 0$ be given,
and for any $\th,\psi \in \Theta$ define $N_{\th\psi}$ on $V \oplus V$ by
\[
N_{\th\psi}(u,v) = \e^{-1}\max\{|\Omega_{\th}(u) - \Omega_{\psi}(v)|: \Omega
\in {\mathcal S}\}.
\]
Then there is a $\d > 0$ such that $N_{\th\psi}$ is a bridge between
$(V_{\th},L_{\th})$ and $(V_{\psi},L_{\psi})$ whenever $d(\th,\psi) < \d$.}

\bigskip
\noindent
{\it\bfseries Proof.} We remark that $N_{\th\psi}$ is constructed from a
finite number of spans each of which is a version of the bridge built in
Proposition $5.4$.  It is clear that $N$ is continuous for the norms.
Because $\Omega$ is not empty, Condition~2 in Definition $5.1$ for a bridge
is clearly satisfied.  We must show that we can choose $\d > 0$ such that
Condition~3 is satisfied.

Since the situation is symmetric in $\th$ and $\psi$, it suffices to show
that for any $u \in V_{\th}$ there is a $v \in V_{\psi}$ satisfying
Condition~3.  If $u \in {\mathbb R}e$ we can always take $v = u$.  So we now
show how to choose a $\d > 0$ which works for all $u$'s not in ${\mathbb
R}e$.  Let $W$ be a subspace of $V$ which is complimentary to ${\mathbb
R}e$, and let $\Sigma_W$ denote the unit $\|\cdot\|_*$-sphere of $W$.  Each
$L_{\th}$ restricts to a norm on $W$ satisfying the conditions of Lemma
$10.1$.  Thus, $(\th,w) \mapsto L_{\th}(w)$ is jointly continuous, and there
is a constant, $c$, such that $c \le L_{\th}(w)$ for all $\th \in \Theta$
and $w \in \Sigma_W$.  By joint continuity and compactness, we can choose
$\d_1$ such that if $d(\th,\psi) < \d_1$ then $|L_{\th}(w) - L_{\psi}(w)| <
\e c^2/2$ for all $w \in \Sigma_W$.  By uniform continuity and the
finiteness of ${\mathcal S}$ we can find $\d \le \d_1$ such that if
$d(\th,\psi) < \d$ then $|\Omega_{\th}(w) - \Omega_{\psi}(w)| < \e c/2$ for
all $\Omega \in {\mathcal S}$ and all $w \in \Sigma_W$.  We show that this
choice of $\d$ works.

To each $u \in V$, viewed as an element of $V_{\th}$, we must find a $v \in
V$, viewed as in $V_{\psi}$, satisfying the requirement of Condition~3 of
Definition $5.1$.  We assume that $u \notin {\mathbb R}e$.  We argue much as
in the proof of Lemma $10.2$.  We consider first the
special case in which $u = w \in \Sigma_W$.  Then clearly
$L_{\psi}(w)
\ne 0$, so that we can set $z = (L_{\th}(w)/L_{\psi}(w))w$.  Then
$L_{\psi}(z) = L_{\th}(w)$, so that we only need to check that
$N_{\th\psi}(w,z) \le L_{\th}(w)$.  We assume that $d(\th,\psi) < \d$.  Then
for each $\Omega \in {\mathcal S}$ we have
\begin{eqnarray*}
|\Omega_{\th}(w) - \Omega_{\psi}(z)| &\le &|\Omega_{\th}(w) -
\Omega_{\psi}(w)| + |\Omega_{\psi}(w) - \Omega_{\psi}(z)| \\
&\le &\e c/2 + |1 - (L_{\th}(w)/L_{\psi}(w))|
|\Omega_{\psi}(w)| \\
&= &\e c/2 + |L_{\psi}(w) - L_{\th}(w)|L_{\psi}(w)^{-1}|\Omega_{\psi}(w)|
\\ &\le &\e c/2 + (\e c^2/2)c^{-1}\|w\|_* = \e c.
\end{eqnarray*}
If we now take $u = te + w$ for some $t \in {\mathbb R}$, and set $v = te
+ z$, we again get $L_{\psi}(v) = L_{\th}(u)$ and
$|\Omega_{\th}(u) - \Omega_{\psi}(v)| \le \e c$. Thus
\[
N(u,v) \le  c \le L_{\th}(u).
\]
Now the two ends of this inequality are homogeneous in $(u,v)$ for positive
scalars, as is the equality $L_{\psi}(v) = L_{\th}(u)$.  But every
element of
$V$ not in ${\mathbb R}e$ is a positive scalar multiple of an element of our
special form
$te + w$.  Thus, for every $u \in V_{\th}$ we can find $v \in V_{\psi}$
satisfying these two relations.  Consequently, $N_{\th\psi}$ is a bridge.\qed

\bigskip
Our main theorem here, which applies immediately to the $B_{\th}$'s of
Section~9, is:

\bigskip
\noindent
{\bf 11.2 THEOREM.} {\em Let $V$ be a finite-dimensional vector space
with distinguished element $e$, and let $\Theta$ be a compact space with
metric $d$.  Let $\{\|\cdot\|_{\th}\}$ be a continuous field of order-unit
norms on $(V,e)$, and let $V_{\th}$ denote $V$ equipped with
$\|\cdot\|_{\th}$.  Let $\{L_{\th}\}$ be a continuous field of Lip-norms on
$(V,e)$.  Then for every $\e > 0$ there is a $\d > 0$ such that if
$d(\th,\psi) < \d$ then}
\[
\dist_q((V_{\th},L_{\th}),(V_{\psi},L_{\psi})) < \e.
\]

\bigskip
\noindent
{\it\bfseries Proof.} Let ${\tilde V} = (V/{\mathbb R}e)$.  Each $L_{\th}$
drops to a norm ${\tilde L}_{\th}$ on ${\tilde V}$.  It is simple to check
that ${\tilde L}_{\th}$ is a continuous field of norms.  According to Lemma
$10.1$ the field $\{{\tilde L}'_{\th}\}$ of dual norms is continuous.  But
the dual of ${\tilde V}$ is canonically identified with the subspace
$V'{}^{\circ}$ of $V'$ consisting of functionals annihilating ${\mathbb
R}e$.  Thus, each ${\tilde L}'_{\th}$ gives a norm $L'_{\th}$ on
$V'{}^{\circ}$, and the field $\{L'_{\th}\}$ is continuous.  Again,
according to Lemma $10.1$, there is a norm, $\|\cdot\|'_*$, on
$V'{}^{\circ}$ such that $L'_{\th} \le \|\cdot\|'_*$ for all $\th \in
\Theta$.  Thus, $\rho_{\th}(\mu,\nu) \le \|\mu-\nu\|'_*$ for all $\th \in
\Theta$ and $\mu,\nu \in S_{\th}$.

Let $\e > 0$ be given.  According to Theorem $10.11$ we can find a finite
family, ${\mathcal S}$, of continuous fields of states such that
$\{\Omega_{\th}: \Omega \in {\mathcal S}\}$ is $(\e/2)$-dense in $S_{\th}$
for $\|\ \|'_*$ for every $\th \in \Theta$.  Define $N_{\th\psi}$ as in the
statement of Proposition $11.1$ except using $\e/2$ instead of $\e$.  By
Proposition $11.1$ we can choose a $\d > 0$ such that $N_{\th\psi}$ is a
bridge for $(V_{\th},L_{\th})$ and $(V_{\psi},L_{\psi})$ whenever
$d(\th,\psi) < \d$.  We show that this $\d$ works.\

Let $\th,\psi \in \Theta$ with $d(\th,\psi) < \d$.  Let $L_{\th\psi}$ be
the Lip-norm on $V_{\th} \oplus V_{\psi}$ defined as at the beginning of
Section~5 using
$N_{\th\psi}$.  For simplicity write $L$ for $L_{\th\psi}$.  We show that
$\dist_H^{\rho_L}(S_{\th},S_{\psi}) < \d$.

Let $\mu \in S_{\th}$.  By the choice of ${\mathcal S}$ there is an
$\Omega^{\mu} \in {\mathcal S}$ such that
\[
\rho_{L_{\th}}(\mu,\Omega_{\th}^{\mu}) = L'_{\th}(\mu - \Omega_{\th}^{\mu})
\le \|\mu - \Omega_{\th}^{\mu}\|'_* < \e/2.\
\]
We now show that $\rho_L(\mu,\Omega^{\mu}_{\psi}) < \e$, thus showing that
$S_{\th}$ is in the $\e$-neighborhood of $S_{\psi}$ for $\rho_{L_{\th\psi}}$.

Let $u,v \in V$, with $(u,v)$ viewed as in $V_{\th} \oplus V_{\psi}$, and
suppose that $L(u,v) \le 1$.  Then $N_{\th\psi}(u,v) \le (\e/2)$, so that
\[
|\Omega_{\th}^{\mu}(u) - \Omega_{\psi}^{\mu}(v)| < \e/2.
\]
Since this is true for all $(u,v)$ with $L(u,v) \le 1$, this means that
\[
\rho_L(\Omega_{\th}^{\mu},\Omega_{\psi}^{\mu}) \le \e/2.
\]
Thus
\[
\rho_L(\mu,\Omega_{\psi}^{\mu}) \le \rho_L(\mu,\Omega_{\th}^{\mu}) +
\rho_L(\Omega_{\th}^{\mu},\Omega_{\psi}^{\mu}) < \e,
\]
as desired.

By reversing the roles of $\th$ and $\psi$ we see that $S_{\psi}$ is in the
$\e$-neighborhood of $S_{\th}$.\qed

\bigskip
When the above theorem is combined with the considerations of Section~9, we
see that we have completed the proof of Theorem $9.2$, our main example.

\section{Completeness}

One of the fundamental facts about the metric space of isometry classes of
ordinary compact metric spaces with the Gromov--Hausdorff metric is that it
is complete.  We now show that the quantum analogue of this fact is true.
(We will show in Theorem 13.15 that for the quantum case this metric
space is also separable.)
We parallel a common classical approach to this matter by showing that,
under suitable conditions, everything can be assembled inside one compact
quantum metric space where we can then use the fact that the set of closed
subsets of a compact metric space is complete for the Hausdorff metric.

We must consider Cauchy sequences of compact quantum metric spaces.  But in
preparation for this we first consider sequences of order-unit spaces.
(Much of what we do works for more general index sets than ${\mathbb N}$,
but we do not need this greater generality here.)  If $\{A_j\}$ is such a
sequence, we denote by $\Pi^b A_j$ the subspace of the full product which
consists of sequences $\{a_j\}$ for which the sequence $\{\|a_j\|\}$ is
bounded.  This certainly includes $e = \{e_j\}$.  We say that a sequence is
non-negative if each of its terms is non-negative.  Then $\Pi^bA_j$ is
clearly an order-unit space.  We let $\oplus A_j$ denote the subspace of
sequences which are $0$ except for a finite number of terms.  It is not an
order-unit space, and $e$ is not in it.

We consider ``compactifications'' of $\oplus A_j$, namely subspaces $C$ of
$\Pi^bA_j$ which contain both $\oplus A_j$ and $e$, and so are order-unit
spaces.  There is an evident projection of $\Pi A_j$, and so of $C$, onto
each $A_n$.  These projections give evident injections of each $S(A_j)$ into
$S(C)$.  Without further comment we will view each $S(A_j)$ as a subset of
$S(C)$ when convenient.

Let $Z_0 = \cup S(A_j)$ inside $S(C)$.  Notice that an element $c$ of $C$ is
non-negative exactly if $\mu(c) \ge 0$ for all $\mu \in Z_0$.  Let $Z =
\mbox{co}(Z_0)$ be the convex hull of $Z_0$ in $S(C)$.  Then we can apply
the bipolar theorem almost exactly as done in the proof of proposition
VII$.5.4$ of \cite{FD} to obtain:

\bigskip
\noindent
{\bf 12.1 LEMMA.} {\em With notation as above, $Z$ is $w^*$-dense in $S(C)$.}

\bigskip
For any positive $n$ let $B_n = \oplus_{j=1}^n A_j$.  In the evident way
each $B_n$ is an order-unit space.  There are evident projections of $C$
onto each $B_n$, and of $B_n$ onto each $A_j$ and $B_j$ for $j \le n$.
Thus, for the evident identifications, we have
\[
S(A_j) \subseteq S(B_j) \subseteq S(B_n) \subseteq S(C)
\]
for $j\le n$.  Furthermore, $S(B_n) = \mbox{co}(\cup^n S(A_j))$, and so $Z =
\cup^{\infty} S(B_n)$, an increasing union.

Suppose now that we have a sequence $\{(A_j,L_j)\}$ of compact quantum
metric spaces.  Suppose further that we have a family $M=\{M_j\}$
where each $M_j$ is a Lip-norm
on $A_j\oplus A_{j+1}$ which induces $L_j$ and $L_{j+1}$.  Define 
$L\equiv L_M$ on
$\Pi A_j$, the full product, by
\[
L(\{a_j\}) = \sup_j\{M_j(a_j,a_{j+1})\}.
\]
Of course $L(e) =0$. We set
\[
C_M = \{\{a_j\} \in \Pi^bA_j: L(\{a_j\}) < \infty\}.
\]
It is easy to check that $C_M$ is an order-unit subspace of
$\Pi^bA_j$ containing $\oplus^{\infty} A_j$, and that $L$ is a seminorm on
$C_M$ which takes value $0$ only on ${\mathbb R}e$.  But without further
hypotheses, $L$ need not induce the $L_j$'s, nor be a Lip-norm.

We remark that the above discussion, and much of what follows, generalizes
to the situation where one has a graph, with a compact quantum metric space
associated to each vertex, and a Lip-norm generalizing the above $M_j$'s
associated to each edge.  But we will not discuss this possibility
further here.

Much as for $L$, define a seminorm, $J_n$ on $B_n$ by
\[
J_n(\{a_j\}) = \sup\{M_j(a_j,a_{j+1}): 1 \le j\le n-1\}.
\]
By a simple adaptation of the proof of Lemma $4.4$ we obtain:

\bigskip
\noindent
{\bf 12.2 LEMMA.} {\em With notation as above, $J_n$ is a Lip-norm on
$B_n$.  For the evident projections of $B_n$ onto $A_m$ for $m \le n$, and
onto $A_m \oplus A_{m+1}$ and $B_m$ for $m\le n-1$, the Lip-norm $J_n$
induces $L_m$, $M_m$, and $J_m$.}

\bigskip
Let us denote the evident projection of $\Pi A_j$ onto $B_n$ by $\pi_n$.

\bigskip
\noindent
{\bf 12.3 LEMMA.} {\em Let $b \in B_n$, and let $\e > 0$ be given.  Then
there is a $c \in \Pi A_j$, the full product, such that $\pi_n(c)= b$ and
$L(c) \le J_n(b) + \e$.}

\bigskip
\noindent
{\it\bfseries Proof.} We choose $b_k \in B_k$ for $k \ge n$ by induction.
Set $b_n = b$.  Suppose that $b_k$ has been chosen.  Since $J_{k+1}$ induces
$J_k$ by Lemma $12.1$, we can find $b_{k+1} \in B_{k+1}$ such that
$\pi_k(b_{k+1}) = b_k$ and $J_{k+1}(b_{k+1}) < J_k(b_k) + \e/2^k$.
Consequently, for each $k > n$ we have $J_k(b_k) < J_n(b) + \e$.  We let $c$
be the unique element of $\Pi A_j$ such that $\pi_k(c) = b_k$ for each $k >
n$.  Then $L(c) \le J_n(b) +\e$ as desired.\qed

\bigskip
The above lemma goes in the direction of saying that $L$ induces $J_n$.
The difficulty is that in general the above $c$ will not be a bounded
sequence.  To ensure that it is bounded we need further hypotheses.

Because of the profusion of indices to follow, we will permit ourselves to
write $\rho(S,T)$ instead of $\dist_H^{\rho}(S,T)$ when $S$ and $T$ are
closed subsets of some metric space with metric $\rho$.

\bigskip
\noindent
{\bf 12.4 LEMMA.} {\em For each $m$ and $n$ with $m < n$ we have}
\[
\rho_{J_n}(S(B_m),S(B_n)) \le \sum_{j=m}^{n-1} \rho_{M_j}(S(A_j),S(A_{j+1})).
\]

\bigskip
\noindent
{\it\bfseries Proof.} This is an adaptation of a calculation done in the
proof of Lemma $4.4$, so we will be a bit succinct here.  Let $\mu_i \in
S(A_i)$ for some $i$ with $m < i \le n$.  Then we can inductively find
$\mu_{i-1},\dots,\mu_m$ with $\mu_j \in S(A_j)$ and
\[
\rho_{J_n}(\mu_j,\mu_{j+1}) \le \rho_{M_j}(S(A_j),S(A_{j+1}))
\]
for $m \le j \le i-1$.  Consequently,
\[
\rho_{J_n}(\mu_i,\mu_m) \le \sum_{j=m}^{i-1} \rho_{M_j}(S(A_j),S(A_{j+1})).
\]
Thus, $\rho_{J_n}(S(B_n),S(A_i))$ satisfies the same bound.  But $S(B_n)$ 
\linebreak 
$= \mbox{co}(\cup_{j=1}^{n-1} S(A_j))$, and $\rho_{J_n}$ is convex (definition
$9.1$ of \cite{R5}), and so $\rho_{J_n}(S(B_m),S(B_n))$ satisfies this bound
too for $i = n$.\qed

\bigskip
When we apply the above lemma twice with $m = 1$ so that $B_m = A_1$, we
obtain:

\bigskip
\noindent
{\bf 12.5 COROLLARY.} {\em For any $n$ we have}
\[
\diam(B_n,J_n) \le \diam(A_1,L_1) + 2 \sum_{j=1}^{n-1}
\rho_{M_j}(S(A_j),S(A_{j+1})).
\]

\bigskip
This lemma and corollary suggest correctly that the key condition for us (as
in the classical case) will be that $\sum_{j=1}^{\infty}
\rho_{M_j}(S(A_j),S(A_{j+1})) < \infty$.  To begin with, the corollary tells
us that in this case there is a number, say $\D$, such that $\diam(B_n,J_n)
\le \D$ for all $n$.  (There is no way that $(C_M,L)$ can have finite
diameter without such a bound.)  We now use this bound to control infinite
sequences.

\bigskip
\noindent
{\bf 12.6 LEMMA.} {\em Let $\{a_j\} \in \prod^{\infty} A_j$, the full
product, and suppose that $L(\{a_j\}) < \infty$.  If there is a constant,
$\D$, such that $\diam(B_n,J_n) \le \D$ for all $n$, then $\{a_j\} \in C_M$,
that is, $\sup_j\{\|a_j\|\} < \infty$.  Furthermore, $L$ on $C_M$ has finite
diameter, no greater than $\D$.  That is, for any $c \in C_M$ we have}
\[
\|c\|^{\sim} \le (\D/2)L(c).
\]

\bigskip
\noindent
{\it\bfseries Proof.} Set $h = L(\{a_j\})$, and let $b_n = (a_1,\dots,a_n)$
for each $n$.  Clearly, $J_n(b_n) \le h$.  Then by proposition $2.2$ of
\cite{R5}, recalled near the end of Section~2, $\|b_n\|^{\sim} \le h\D/2$.
Set $e_n = e_{B_n}$.  Then this means that there is a $t \in {\mathbb R}$
such that for each $n$
\[
\|b_n - te_n\| \le h\D/2.
\]
For each $n$ set $G_n= \{t: \|b_n - te_n\| \le h\D/2\}$.  Then $G_n$ is a
non-empty closed bounded subset of ${\mathbb R}$.  Furthermore, it is easily
seen that $G_n\supseteq G_{n+1}$.  By compactness there is a $t \in
\cap^{\infty} G_n$.  For this $t$ we have $\|a_j - te_j\| \le h\D/2$ for
all $j$, so that $\|a_j\| \le t + h\D/2$ for all $j$.  Thus, $\{a_j\} \in
C_M$, and $\|c\|^{\sim} \le (\D/2)L(c)$ where $c = \{a_j\}$.\qed

\bigskip
We can now apply this lemma to the context of Lemma $12.3$.  It assures us
that if $\diam(B_n,J_n) \le \D$ for all $n$, then the $c$ of Lemma $12.3$ is
bounded.  Of course, we always have $J_n(\pi_n(c)) \le L(c)$.  When all of
this is combined with Proposition $3.1$, we obtain:

\bigskip
\noindent
{\bf 12.7 PROPOSITION.} {\em Suppose that there is a constant, $\D$, such
that $\diam(B_j,J_j) \le \D$ for all $j$.  Then for any $b_n\in
B_n$ and
$\e >  0$ there is a $c \in C_M$ such that $\pi_n(c) = b_n$ and $L(c) \le
J_n(b_n) + \e$.  In other words, $L$ on $C_M$ induces $J_n$ via $\pi_n$.
Thus, the natural inclusion of $S(B_n)$ into $S(C_M)$ is isometric for
$\rho_{J_n}$ and $\rho_L$.}

\bigskip
Note that $L$ still need not be a Lip-norm.

\bigskip
\noindent
{\bf 12.8 LEMMA.} {\em Suppose now that $\sum^{\infty}
\rho_{M_j}(S(A_j),S(A_{j+1})) < \infty$.  Then the metric space $(Z,\rho_L)$
is totally bounded.}

\bigskip
\noindent
{\it\bfseries Proof.} Let $\e > 0$ be given.  Then there is an $m$ such that
$$
\sum_{j=m}^{\infty} \rho_{M_j}(S(A_j),S(A_{j+1})) < \e.
$$
It
follows from Lemma $12.4$ that for each $n \ge m$ we have
$$\rho_L(S(B_m),S(B_n)) \le \e.
$$
This says that $S(B_m)$ is $\e$-dense for $\rho_L$ in $Z$.  But
$S(B_m)$ is compact for the topology from $\rho_L$, since $\rho_{J_m}$ is a
Lip-norm and $\rho_L = \rho_{J_m}$ on $S(B_m)$.  Thus, $S(B_m)$ is totally
bounded for $\rho_L$, and a finite subset $\e$-dense in $S(B_m)$ will be
$2\e$-dense in $Z$.\qed

\bigskip
Let ${\hat Z}$ denote the abstract completion of $Z$ for $\rho_L$.  We let
$\rho_L$ denote also the extension of $\rho_L$ to ${\hat Z}$.  Under the
hypothesis of the above lemma, ${\hat Z}$ is compact.  Each $c \in C_M$ is
Lipschitz when viewed as a function on $Z$, and so extends to give a
Lipschitz function, ${\hat c}$, on ${\hat Z}$.  It is clear that if $c \ge
0$ in $C_M$ then ${\hat c} \ge 0$ as function, and that ${\hat e} \equiv 1$.
Thus, $C_M$ is realized isometrically isomorphically as an order-unit space of
continuous functions on ${\hat Z}$.  Each $x \in {\hat Z}$ determines an
element $\mu_x$ of $S(C_M)$ by evaluation.  By the definition of the
$w^*$-topology the map $x \mapsto \mu_x$ is continuous from the metric
topology on ${\hat Z}$.  But each element of $Z \subseteq {\hat Z}$ is
carried to itself in $S(C_M)$, and $Z$ is dense in $S(C_M)$.  Since ${\hat Z}$
is compact, its image is closed, and so is all of $S(C_M)$.

Let us see now that the map $x \mapsto \mu_x$ is injective.  Let $x,y \in
{\hat Z}$ with $x \ne y$, and set $\g = \rho_L(x,y)$.  We can find
$\mu,\nu\in Z$ such that $\rho_L(x,\mu) < \g/4$ and $\rho_L(y,\nu) < \g/4$.
Thus $\rho_L(\mu,\nu) > \g/2$, so that we can find $c \in C_M$ with $L(c) \le
1$ and $|\mu(c) - \nu(c)| > \g/2$.  But $L_{\rho_L}({\hat c}) \le 1$ so that
$|{\hat c}(x) - {\hat c}(\mu)| < \g/4$ and $|{\hat c}(y) - {\hat c}(\nu)| <
\g/4$.  Thus $|{\hat c}(x) - {\hat c}(y)| > 0$.

It follows that $x \mapsto \mu_x$ is a homeomorphism of ${\hat Z}$ onto
$S(C_M)$.  From this it is clear that the metric topology on $S(C_M)$ from
$\rho_L$ agrees with the $w^*$-topology.  We have thus obtained:

\bigskip
\noindent
{\bf 12.9 THEOREM.} {\em Let $\{(A_j,L_j)\}$ be a sequence of compact
quantum metric spaces, and for each $j$ let $M_j$ be a Lip-norm on $A_j
\oplus A_{j+1}$ which induces $L_j$ and $L_{j+1}$.  Let $L$ and $C_M$ be
defined as above.  If}
\[
\sum^{\infty} \rho_{M_j}(S(A_j),S(A_{j+1})) < \infty,
\]
{\em then $L$ is a Lip-norm on $C_M$.}

\bigskip
When we view the above situation within $S(C_M)$, we have 
$$
\sum_j^{\infty}
\rho_L(S(A_j),S(A_{j+1})) < \infty,
$$
which clearly implies that
$\{S(A_j)\}$ is a Cauchy sequence for $\rho_L$ ($= \dist_H^{\rho_L}$).
Since the space of compact convex subsets is complete for Hausdorff
distance, this sequence has a limit, say $K$, in $S(C_M)$.  We let
$\pi_{\infty}$ denote the process of restricting elements of $C_M$ to $K$, and
we let $D = \pi_{\infty}(C_M)$.  We let $L_D$ be the quotient of $L$ for
$\pi_{\infty}$.  Thus, $L_D$ is a Lip-norm by Proposition $3.1$.  Now,
according to Proposition $5.7$, we have
\[
\dist_q(A_j,D) \le \dist_H^{\rho_L}(S(A_j),K).
\]
  From this we obtain:

\bigskip
\noindent
{\bf 12.10 COROLLARY.} {\em Under the hypotheses of Theorem $12.9$ there is
a compact quantum metric space $(D,L_D)$ to which the sequence
$\{(A_j,L_j)\}$ converges for $\dist_q$.}

\bigskip
\noindent
{\bf 12.11 THEOREM.} {\em The metric space of isometry classes of compact
quantum metric spaces, with the metric $\dist_q$, is complete.}

\bigskip
\noindent
{\it\bfseries Proof.} Let $\{(A_j,L_j)\}$ be a Cauchy sequence.  It suffices
to show that a subsequence converges.  But we can choose a subsequence,
still denoted by $\{(A_j,L_j)\}$, such that $\sum^{\infty}_j
\dist_q(A_j,A_{j+1}) < \infty$.  Then for each $j$ we can find a Lip-norm
$M_j$ on $A_j \oplus A_{j+1}$ which induces $L_j$ and $L_{j+1}$, and is such
that
\[
\rho_{M_j}(S(A_j),S(A_{j+1})) \le \dist_q(A_j,A_{j+1}) + 1/2^j.
\]
It follows that $\sum \rho_{M_j}(S(A_j),S(A_{j+1})) < \infty$.  Then
according to Corollary $12.10$ there is a compact quantum metric space to
which $\{(A_j,L_j)\}$ converges for $\dist_q$.\qed

\section{Finite approximation and compactness}

An elementary property of ordinary compact metric spaces is that for every
$\e > 0$ there is a finite set which is $\e$-dense, and thus approximates 
within $\e$ for
Gromov--Hausdorff distance.  We begin this section by showing that the
quantum analogue of this is true.  (Of course, in Sections~8 to 11 we were
already using a special instance of this in which the ``quantum finite
sets'' were chosen in a particularly useful way.)  Afterward we use this
finite approximation and some of our earlier results to prove a quantum
analogue of a fundamental compactness criterion of Gromov for families of
ordinary compact metric spaces for Gromov--Hausdorff distance.  In
addition to its spectacular application in \cite{G1}, Gromov's compactness
theorem has been widely used in Riemannian geometry.  See, for example, the
exposition and references in \cite{G2}, \cite{Fu}, \cite{Sa}, \cite{Sh}.  It
remains to be seen what applications will be found for our quantum version.

Our finite approximation theorem is:

\bigskip
\noindent
{\bf 13.1 THEOREM.} {\em Let $(A,L)$ be a compact quantum metric space.  For
every $\e > 0$ there is a compact quantum metric space, $(B,L_B)$, such that
$B$ has finite dimension and $\dist_q(A,B) < \e$.  In fact, we can take
$(B,L_B)$ to be a quotient of $(A,L)$. Alternatively, or if $A$ is
infinite-dimensional, we can take $B$ to be order-isomorphic
to $C(X)$ for some finite set $X$.}

\bigskip
\noindent
{\it\bfseries Proof.} Since $S(A)$ with $\rho_L$ is a compact metric space,
we can find a finite subset, $F$, of $S(A)$ which is $\e$-dense in $S(A)$
for $\rho_L$.  Let $K = \mbox{co}(F)$, which is closed.  Let $B$ consist of
the restrictions to $K$ of elements of $A$, and let $L_B$ be the quotient
Lip-norm, as discussed in Section~3.  Note that $B$ has finite dimension
since its elements are determined by their values on $F$.  Also, $K$ is
naturally identified with $S(B)$ by Proposition $3.5$.  Then
\[
\dist_q(A,B) \le \dist_H^{\rho_L}(S(A),K) < \e
\]
by Proposition $5.7$. If $A$ is infinite-dimensional, then we can
perturb the elements of $F$ very slightly so that they are linearly
independent. Then $K$ will be a simplex, and $B$ will be 
order-isomorphic to $C(F)$. If $A$ is finite-dimensional, then we
must use a simple generalization of the doubling construction
of Example 5.6 to make sufficiently many very close copies of $A$
so that the elements of $F$ can be chosen to be linearly independent.
\qed

\bigskip

We emphasize that the Lip-norms on $C(X)$ will often not come from
metrics on $X$.
In contrast to this, see Theorem 13.16, where in effect the $L_B$'s are
required to come from ordinary metrics on the $X$'s.

A simple observation which seems important is that the topology
of a compact topological space can not be recovered from the data consisting
of an increasing dense union (with inclusion maps) of finite 
subsets with the relative topology (discrete), whereas the metric 
on a compact metric space (and so the space itself by completion)
can be recovered from an increasing dense
union of finite subsets with the relative metrics. It is easy to see
that this latter is still true for compact quantum metric
spaces if one uses finite-dimensional quotients.

We now recall the statement of the classical Gromov compactness theorem,
along the lines given in his original formulation of the theorem \cite{G1},
or as in theorem $6.3$ of the appendix of \cite{Sa}.  Let ${\mathcal M}$
denote the set of isometry equivalence classes of compact metric spaces,
equipped with the metric $\dist_{GH}$.  Our notation below will not
distinguish between equivalence classes and their representatives.

\bigskip
\noindent
{\bf 13.2 DEFINITION.} Let $(X,\rho)$ be a compact metric space.  We define
a function, $\Cov_{\rho}$, from the strictly positive real numbers,
${\mathbb R}^+$, to the positive integers, ${\mathbb N}^+$, by setting
$\Cov_{\rho}(\e)$ to be the smallest number of $\e$-balls for $\rho$
which are needed to cover $X$.  We call the function $\Cov_{\rho}$ 
the {\em covering growth} of $\rho$.

\bigskip
Of course, the most interesting aspect of $\Cov_{\rho}$ is how it
grows as $\e \rightarrow 0$.

\bigskip
\noindent
{\bf 13.3 GROMOV'S COMPACTNESS THEOREM.} {\em Let ${\mathcal S}$ be a subset
of ${\mathcal M}$.  Then ${\mathcal S}$ is totally bounded for $\dist_{GH}$
(hence has compact closure) if and only if}

\begin{itemize}
\item[1)] {\em there is a constant $D$ such that $\diam(X,\rho) \le D$ for
all $(X,\rho) \in {\mathcal S}$; and}

\item[2)] {\em there is a function, $G$, from ${\mathbb R}^+$ to ${\mathbb
N}^+$ such that $\Cov_{\rho}(\e) \le G(\e)$ for every $(X,\rho)$ in
${\mathcal S}$ and every $\e > 0$.}
\end{itemize}

\bigskip
We now formulate an analogue of this theorem for compact quantum metric
spaces.  For this purpose we need to mention that already in the literature
on ordinary metric spaces there are several alternative measures of growth
besides $\Cov_{\rho}$, which are nevertheless equivalent for the purposes of
the compactness theorem.  So it is to be expected that this happens also in
the quantum setting.  In fact, the proof of Theorem $13.1$ already suggests
three such alternatives, which we now formalize.

\bigskip
\noindent
{\bf 13.4 DEFINITION.} Let $(A,L)$ be a compact quantum metric space.  For
each $\e > 0$ we set:

\begin{itemize}
\item[1)] $\Fin_L(\e)$ is the smallest integer $n$ such that there is
a compact quantum metric space $(B,L_B)$ such that
$\dist_q(A,B) < \e$ and $\dim(B) \le n$ (where $\dim(B)$ is the vector-space
dimension of $B$).

\item[2)] $\Cov_L(\e)$ is the smallest integer $n$ such that there is a
compact quantum metric space $(B,L_B)$ and a surjection $\pi$ of $A$ onto
$B$ such that $L$ induces $L_B$, while $\dist_H^{\rho_L}(S(B),S(A)) < \e$ and
$\dim(B) \le n$.

\item[3)] $\Scv_L(\e)$ is the smallest integer $n$ such that there is
a subset of $S(A)$ which is $\e$-dense for $\rho_L$ and has only $n$ elements.
\end{itemize}
We will sometimes write $\Scv_A$ instead of $\Scv_L$, etc., when
this seems helpful and $L$ is understood.

\bigskip
Of course, $\Scv_L(\e)$ is always finite, and it is easily seen from
the proof of Theorem $13.1$ that
\[
\Scv_L(\e) \ge \Cov_L(\e) \ge \Fin_L(\e).
\]

All this suggests, of course, that it would be
interesting to study the Kolmogorov $\e$-entropy of $(S(A), \rho_L)$ for
specific examples, or ``quantum'' versions of it, say $\ln_2(Fin_L(\e))$.
See the discussion of fractal and Hausdorff-Besicovich dimension
in \cite{Brn}.

Let now ${\mathcal Q}$ denote the set of isometry equivalence classes of
compact quantum metric spaces, equipped with the metric $\dist_q$.

\bigskip
\noindent
{\bf 13.5 THE QUANTUM GROMOV COMPACTNESS THEOREM.} {\em Let ${\mathcal S}$
be a subset of ${\mathcal Q}$.  If ${\mathcal S}$ is totally bounded for
$\dist_q$ then}

\begin{itemize}
\item[1)] {\em there is a constant, $D$, such that $\diam(A,L) \le D$ for
all $(A,L) \in {\mathcal S}$; and}

\item[2)] {\em there is a function, $G$, from ${\mathbb R}^+$ to ${\mathbb
N}^+$ such that $ G(\e) \ge \Scv_L(\e) \  (\ge \Fin_L(\e))$ for 
every $(A,L) \in {\mathcal S}$
and every $\e > 0$.}
\end{itemize}

\noindent
{\em Conversely, if}

\begin{itemize}
\item[1)] {\em there is a constant, $D$, such that $\diam(A,L) \le D$ for
all $(A,L) \in {\mathcal S}$; and}

\item[2)] {\em there is a function, $G$, from ${\mathbb R}^+$ to ${\mathbb
N}^+$ such that $G(\e) \ge \Fin_L(\e)$ for every $(A,L) \in {\mathcal S}$
and every $\e > 0$,}
\end{itemize}
\noindent
{\em then $\mathcal S$ is totally bounded in $\mathcal Q$.}

\bigskip

Before launching into the proof of this theorem, we point out that
it is relevant to the situation studied in Sections 8--11. Fix a compact
Lie group $G$ (possibly not connected), and fix a length function,
$\ell$, on $G$. Let ${\mathcal E}(G, \ell)$ denote the subset of
$\mathcal Q$ consisting of the equivalence classes of pairs $(A, L_A)$
where $A$ is a unital $C^*$-algebra and $L_A$ is a Lip-norm 
defined, using $\ell$, by an ergodic action of $G$ on $A$, using the formula
at the beginning of Section 8. Then Theorem 13.5 will
tell us that ${\mathcal E}(G, \ell)$ is a
totally bounded subset of $\mathcal Q$. To see this, note that by lemma 2.4
of \cite{R4} every element of ${\mathcal E}(G, \ell)$ has radius no larger
than $\int_G \ell(x) dx$, where we use the Haar measure which gives
$G$ mass $1$. Thus condition 1 of Theorem 13.5 is satisfied. We now show 
how to define the function $G$ (not to be confused with the group $G$) 
for condition 2. Let $\e > 0$ be given. Choose a 
representation, $\pi$, of $G$ in the way discussed after Proposition 8.1,
and let $\{\delta_n\}$ be the sequence obtained in Theorem 8.2. Choose
$n$ such that $\delta_n < \e$. Let ${\hat G}_n$ be as defined after 
Proposition 8.1. We set
$$
G(\e) = \sum \{(\dim(\rho))^2 : \rho \in {\hat G}_n \}.
$$
This works for the following reasons. Let $\a$ be an ergodic action
of $G$ on a unital $C^*$-algebra $A$, and let $B_n$ be defined as 
above. Now a key assertion of the main theorem
of \cite{HLS} (see also \cite{Wa}) states that any irreducible 
representation of $G$ occurs in $A$ with multiplicity no greater
than its dimension. Thus $\dim(B_n) \leq G(\e)$, while
$\dist_q(A, B_n) \leq \e$ according to Theorem 8.2. Consequently 
condition 2 of Theorem 13.5 is satisfied. Theorem 9.2 is consistent 
with this observation, and the results we will give in \cite{R7}
will be also.

\bigskip
\noindent
{\it\bfseries Proof of Theorem 13.5.} Suppose first that $\mathcal S$
is totally bounded.
Although Condition~1 can be dealt with more directly, the following approach
seems interesting.

\bigskip
\noindent
{\bf 13.6 LEMMA.} {\em The function $(A,L) \mapsto \diam(A,L)$ from
${\mathcal Q}$ to ${\mathbb R}$ is Lipschitz.  In fact,}
\[
|\diam(A,L_A) - \diam(B,L_B)| \le 2 \dist_q(A,B)
\]
{\em for all $(A,L_A)$ and $(B,L_B)$ in ${\mathcal Q}$.}

\bigskip
\noindent
{\it\bfseries Proof.} Let $(A,L_A)$ and $(B,L_B) \in {\mathcal Q}$, and let
$d = \dist_q(A,B)$.  Given $\e > 0$ there is a Lip-norm $L$ on $A \oplus B$
inducing $L_A$ and $L_B$ and such that $\dist^{\rho_L}(S(A),S(B)) < d +
\e$.  If $\mu_1,\mu_2 \in S(A)$, then there are $\nu_1,\nu_2 \in S(B)$ with
$\rho_L(\mu_j,\nu_j) < d + \e$.  Thus
\[
\rho_{L_A}(\mu_1,\mu_2) \le \diam(B,L_B) + 2(d+\e).
\]
Because $\e$ is arbitrary, it follows that
\[
\diam(A,L_A) \le \diam(B,L_B) + 2d.
\]
But we can reverse the roles of $A$ and $B$ to obtain the desired
inequality.\qed

\bigskip
It follows that since ${\mathcal S}$ is a totally bounded subset 
of ${\mathcal Q}$, the set of diameters of elements of ${\mathcal S}$ 
will be a
totally bounded subset of ${\mathbb R}$, and so Condition~1 is satisfied.

We show now how to obtain a function $G$ for Condition~2 of the first part.  
Let $\e > 0$ be
given.  Since ${\mathcal S}$ is totally bounded, there is a finite subset,
$F$, of ${\mathcal S}$ which is $(\e/3)$-dense in ${\mathcal S}$ for
$\dist_q$.  We set
\[
G(\e) = \max\{\Scv_L(\e/3): (A,L) \in F\}.
\]
From the triangle inequality it follows that this choice of
$G(\e)$ works.

We now show that the conditions of the second part
of Theorem $13.5$ are sufficient.  Thus we
assume that $D$ and the function $G$ are given.  Let ${\mathcal S}_G^D$
consist of {\em all} (equivalence classes of) compact quantum metric spaces
$(A,L)$ of diameter $\le D$ for which $\Fin_L \le G$.  It suffices to show
that ${\mathcal S}_G^D$ is totally bounded.  Let $\e > 0$ be given.  Let
$(A,L) \in {\mathcal S}_G^D$.  Since $\Fin_L(\e) \le G(\e)$, there is a
compact quantum metric space $(B,L_B)$ with $\dim(B) \le G(\e)$ such that
$\dist_q(A,B) < \e$.  It follows from Lemma $13.6$ that
\[
\diam(B,L) \le D + 2\e.
\]
Let ${\mathcal S}_n(d)$ denote the subset of ${\mathcal Q}$ consisting of
elements of dimension $\le n$ and diameter $\le d$.  Thus we have just seen
that ${\mathcal S}_G^D$ is contained in the $\e$-neighborhood of ${\mathcal
S}_{G(\e)}(D + 2\e)$.  It thus suffices to show that each ${\mathcal
S}_n(d)$ is totally bounded.  By a simple scaling argument we can reduce to
the case of $d = 2$, so radius $\le 1$.  We set ${\mathcal S}_n = {\mathcal
S}_n(2)$.

Let ${\mathcal S}^j$ denote the subset of ${\mathcal Q}$ consisting of
elements of dimension exactly $j$ and radius $\le 1$.  Then ${\mathcal S}_n$
is the disjoint union of the ${\mathcal S}^j$ for $j \le n$.  Thus we see
that the crux of the matter is to show:

\bigskip
\noindent
{\bf 13.7 PROPOSITION.} {\em For every integer $n$ the set ${\mathcal S}^n$
is totally bounded for $\dist_q$.}

\bigskip
\noindent
{\it\bfseries Proof.} We fix $n$ for the rest of the discussion.  We proceed
by first finding  a kind of ``standard position'' for order-unit spaces of
dimension $n$.  Fix a real vector space $V$ of dimension $n$, an inner
product on it, and its corresponding Euclidean norm, $\|\cdot\|_E$.
According to a theorem of F.~John (proposition $9.12$ of \cite{TJ}) for any
other normed vector space, $(U,\|\cdot\|_U)$, of dimension $n$ there is a
linear operator, $T$, from $U$ onto $V$ such that $\|T\|\|T^{-1}\| \le
\sqrt{n}$.  (I am indebted to Ed~Effros for suggesting to me that this
theorem of F.~John might be useful for the present considerations.)  Fix
further a vector $e \in V$ such that $\|e\|_E = 1$.  Suppose now that
$(U,e_U,\|\cdot\|_U)$ is an order-unit space of dimension $n$.  Let $T: U
\rightarrow V$ be as above.  By multiplying $T$ by a constant we can arrange
that $\|Te_U\| = 1$.  By then composing $T$ with an orthogonal
transformation we can arrange that $Te_U = e$.  Since it follows that
$T^{-1}e = e_U$, we see that $\|T\| \ge 1$ and $\|T^{-1}\| \ge 1$, so that
$\|T\| \le \sqrt{n}$ and $\|T^{-1}\| \le \sqrt{n}$.  We can now use $T$ to
transfer to $V$ the norm of $U$.  Thus we see that every order-unit space of
dimension $n$ is (isometrically order) isomorphic to one coming from an
order-unit norm, $\|\cdot\|$, on $(V,e)$ such that
\[
(1/\sqrt{n})\|v\|_E \le \|v\| \le \sqrt{n}\|v\|_E.
\]
It is convenient to relate this to our earlier notation by defining a new
Euclidean norm $\|\cdot\|_*$ by $\|\cdot\|_* = \sqrt{n}\|\cdot\|_E$.  We
then summarize the above by:

\bigskip
\noindent
{\bf 13.8 LEMMA.} {\em Let $V$ be a vector space of dimension $n$, equipped
with a Euclidean norm $\|\cdot\|_*$ and a distinguished vector $e$ such that
$\|e\|_* = \sqrt{n}$.  Then every order-unit space of dimension $n$ is
isomorphic to one coming from an order-unit norm, $\|\cdot\|$, on $(V,e)$
such that for all $v \in V$}
\[
n^{-1}\|v\|_* \le \|v\| \le \|v\|_*.
\]

\bigskip
Let ${\mathcal O}_n$ denote the set of order-unit norms on $(V,e)$
satisfying the above inequalities with respect to $\|\cdot\|_*$.  Let $V'$
be the dual vector space to $V$ with dual norm $\|\cdot\|'_*$ and
distinguished functional $\eta$ from $e$.  Let ${\mathcal O}'_n$ denote the
set of base-norms on $(V',\eta)$ which are the duals of the order-unit norms
in ${\mathcal O}_n$.  For $\|\cdot\|' \in {\mathcal O}'_n$ we will have
\[
n\|x\|'_* \ge \|x\|' \ge \|x\|'_*
\]
for $x \in V'$.  Let $B'_*$ denote the unit ball in $V'$ for
$\|\cdot\|'_*$.  View the elements of ${\mathcal O}'_n$ as functions on
$B'_*$.  Because these functions are dominated by $n\|\cdot\|'_*$, it is
easily seen that ${\mathcal O}'_n$ is a bounded equicontinuous family of
functions on $B'_*$.  Thus, by the Arzela--Ascoli theorem, ${\mathcal O}'_n$
is totally bounded for the supremum norm.  From this we obtain:

\bigskip
\noindent
{\bf 13.9 LEMMA.} {\em For any $\d > 0$ we can find a finite subset,
$F_{\d}$, of ${\mathcal O}'_n$ such that if $\|\cdot\|_1 \in {\mathcal O}'_n$
then there is a $\|\cdot\|_2$ in $F_{\d}$ such that}
\[
|\|x\|'_1 - \|x\|'_2| \le \d\|x\|'_*
\]
{\em for all $x \in V'$.}

\bigskip
We are then in position to apply Lemma $10.7$, with $\d = \e/4$.  We obtain:

\bigskip
\noindent
{\bf 13.10 LEMMA.} {\em Let $\e > 0$ be given, with $\e < 2$.  Then for any
$\|\cdot\|_1 \in {\mathcal O}'_n$ there is a $\|\cdot\|_2 \in F_{\e/4}$
(chosen as in Lemma $13.9$) such that}
\[
\dist^*_H(S_1,S_2) < \e.
\]

\bigskip
We must now bring Lip-norms into the picture.  Let $L_0$ be a Lip-norm on
$(V,e)$.  Again, ``Lip-norm'' here simply means that the null-space of $L$
is spanned by $e$.  Let us consider two order-unit norms, $\|\cdot\|_1$ and
$\|\cdot\|_2$ in ${\mathcal O}_n$.  We will write $(V_1,L_1)$ when we think
of $V$ as equipped with $\|\cdot\|_1$ and with $L_0$ as Lip-norm, and
similarly for $(V_2,L_2)$.  The corresponding state spaces will be denoted
by $S_1$ and $S_2$, etc.

\bigskip
\noindent
{\bf 13.11 LEMMA.} {\em Let $\|\cdot\|_1$ and $\|\cdot\|_2 \in {\mathcal
O}_n$, and let $L_0$ be a Lip-norm on $(V,e)$ giving $V_1$ and $V_2$ radius
$\le 1$.  Suppose that}
\[
\dist_H^*(S_1,S_2) < \d.
\]
{\em Then $\dist_q(V_1,V_2) \le n\d$.}

\bigskip
\noindent
{\it\bfseries Proof.} Let $\g > 0$ be given.  By hypothesis we can find a
finite subset, $F$, of $S_1 \x S_2$ such that if $(\mu,\nu) \in F$ then
$\|\mu-\nu\|'_* < \d$, the first coordinates of elements of $F$ are
$\g$-dense in $S_1$ for $\rho_{L_1}$, and the second coordinates of elements
of $F$ are $\g$-dense in $S_2$ for $\rho_{L_2}$.  Define $N$ on $V \oplus V$
by
\[
N(u,v) = (n\d)^{-1}\max\{|\mu(u)-\nu(v)|: (\mu,\nu) \in F\}.
\]
We show that $N$ is a bridge between $(V_1,L_1)$ and $(V_2,L_2)$.  The first
two conditions of Definition $5.1$ clearly hold.  For the third condition,
let $u\in V_1$ be given.  Set $v = u$ viewed as element of $V_2$.  Clearly
$L_2(v) = L_1(u)$.  For $(\mu,\nu) \in F$ we have
\[
|\mu(u) - \nu(v)| = |(\mu-\nu)(u)| \le \|\mu-\nu\|'_* \|u\|_*^{\sim} \le
\d(n\|u\|_1^{\sim}) \le \d nL_1(u).
\]
Thus $N(u,v) \le L_1(u)$.  We can do the same calculation with the roles of
$V_1$ and $V_2$ reversed.  Thus $N$ is a bridge.

Define $L$ on $V_1 \oplus V_2$ as before in terms of $N$ (briefly, by $L =
(L_1 \vee L_2) \vee N$).  Let $(\mu,\nu) \in F$.  For every $(u,v) \in V_1
\oplus V_2$ such that $L(u,v) \le 1$ we have $|\mu(u) - \nu(v)| \le n\d$, so
that $\rho_L(\mu,\nu) \le n\d$.  For any $\mu_0 \in S_1$ we can find
$(\mu,\nu) \in F$ such that $\rho_L(\mu_0,\mu) = \rho_{L_1}(\mu_0,\mu) <
\g$.  Thus $\rho_L(\mu_0,\nu) < n\d + \g$.  In the same way we see that for
every $\nu_0 \in S_2$ there is a $\mu \in S_1$ such that $\rho_L(\mu,\nu_0)
< n\d + \g$.  Thus $\dist_H^{\rho_L}(S_1,S_2) < n\d + \g$.  Since $\g$ is
arbitrary, it follows that $\dist_q(V_1,V_2) \le n\d$.\qed

\bigskip
Let us now combine the above lemma with Lemma $13.10$.  Let $\e > 0$ be
given, with $\e < 2$, and find $F_{\e/4} \subset {\mathcal O}_n$ as in Lemma
$13.9$.  Let $\|\cdot\|_1 \in {\mathcal O}_n$, and let $L_1$ be a Lip-norm
on $V_1$ of radius $\le 1$.  By the definition of $F_{\e/4}$ we can find
$\|\cdot\|_2 \in F_{\e/4}$ such that both
\[
|\|x\|'_1 - \|x\|'_2| \le (\e/4)\|x\|'_*
\]
for all $x \in V'$, and $\dist_H^*(S_1,S_2) < \e$.  Now let $L_2$ be $L_1$,
but viewed as a Lip-norm on $V_2$.  The small difficulty with applying our
previous lemma is that $L_2$ may not have radius $\le 1$.  However, for any
$x \in V'$ we do have
\[
L'_2(x) = L'_1(x) \le \|x\|'_1 \le \|x\|'_2 + (\e/4)\|x\|'_* \le
(1+\e/4)\|x\|'_2.
\]
In other words, $L_2$ has radius $\le 1 + (\e/4)$.  Let $r = 1 + (\e/4)$,
and set $M = rL_1$.  Then $M$ is a Lip-norm, and when viewed on either
$V_1$ or $V_2$ it has radius $\le 1$.  Thus we can apply the previous lemma
to conclude that
\[
\dist_q((V_1,M_1),(V_2,M_2)) \le n\e.
\]
But changing from $M$ to $L$ clearly just multiplies distances by $r$.  We
thus obtain
\[
\dist_q((V_1,L_1),(V_2,L_2)) \le n\e(1 + \e/4).
\]
By changing the meaning of $\e$ and $F_{\e}$ accordingly, we see that we
obtain:

\bigskip
\noindent
{\bf 13.12 LEMMA.} {\em For any $\e > 0$ there is a finite subset, $F_{\e}$,
of ${\mathcal O}_n$ such that if $\|\cdot\|_1 \in {\mathcal O}_n$ and if
$L_1$ is a Lip-norm on $V_1$ of radius $\le 1$, then there is a $\|\cdot\|_2
\in F_{\e}$ and a Lip-norm $L_2$ on $V_2$ such that}
\[
\dist_q(V_1,V_2) < \e
\]
{\em (and $\mbox{\em radius}(L_2) \le 1 +\e$).}

\bigskip
If we combine this with the earlier discussion, we have arrived at the point
where we see that, given $G$ and $\e > 0$, there is a finite set $F$ of
order-unit spaces of dimension $\le G(\e)$ such that if $(A,L) \in {\mathcal
S}_G^1$ then there is a $B \in F$ and a Lip-norm $L_2$ on $B$ such that
$\dist_q(A,B) <\e$ and $\mbox{diameter}(B,L) \le 1 + \e$.  To conclude the
proof of Theorem $14.5$ it thus suffices to prove:

\bigskip
\noindent
{\bf 13.13 PROPOSITION.} {\em Let $(A,e,\|\cdot\|)$ be a finite-dimensional
order-unit space, and let $r \in {\mathbb R}^+$.  Let ${\mathcal Q}(A,r)$
denote the set of all elements of ${\mathcal Q}$ represented by $A$
equipped with a Lip-norm of $\mbox{\em radius} \le r$.  Then ${\mathcal
Q}(A,r)$ is totally bounded for $\dist_q$.}

\bigskip
\noindent
{\it\bfseries Proof.} By a simple scaling argument we can assume that $r =
1$.  The Lip-norms on $A$ of radius $\le 1$ correspond by duality to the
norms $L'$ on $A'{}^{\circ}$ such that $L' \le \|\cdot\|'$.  We will denote
this collection of norms by ${\mathcal N}_1(A)$.  For $L' \in {\mathcal
N}_1(A)$ and $\l_1,\l_2 \in A'{}^{\circ}$ we have
\[
|L'(\l_1) - L'(\l_2)| \le L'(\l_1-\l_2) \le \|\l_1-\l_2\|'.
\]
Thus, when the elements of ${\mathcal N}_1(A)$ are viewed as functions on
the unit  $\|\cdot\|'$-ball, ${\mathcal N}_1(A)$ is a bounded equicontinuous
family of functions.  Because the unit ball is compact since $A$ is of
finite dimension, we can apply the Arzela--Ascoli theorem to conclude that
${\mathcal N}_1(A)$ is totally bounded.  (We remark that usually ${\mathcal
N}_1(A)$is not closed, reflecting the fact that elements of ${\mathcal
Q}(A,1)$ can converge to quantum compact metric spaces of strictly lower
dimension.)  Thus, given $\e > 0$, we can find a finite subset, $F_{\e}$, of
${\mathcal N}_1(A)$ such that if $L' \in {\mathcal N}_1(A)$ then there is an
$L'_1 \in F_{\e}$ such that
\[
|L'(\l) - L'_1(\l)| \le \e\|\l\|'
\]
for all $\l \in A'{}^{\circ}$.  We now see that we need the following
quantum analogue of corollary $6.24$ of \cite{Be}.  For this analogue we do
not need our spaces to be finite-dimensional.

\bigskip
\noindent
{\bf 13.14 PROPOSITION.} {\em Let $(A,e,\|\cdot\|)$ be an order-unit space,
and let $L_1$ and $L_2$ be two Lip-norms on $A$.  If there is a $\d > 0$
such that}
\[
|L'_1(\l) - L'_2(\l)| \le \d\|\l\|'
\]
{\em for all $\l \in A'{}^{\circ}$, then}
\[
\dist_q((A,L_1),(A,L_2)) \le \d.
\]

\bigskip
\noindent
{\it\bfseries Proof.} We reduce first to the finite-dimensional situation.
Let $\e \ge 0$ be given. We argue along the lines of the proof
of Theorem 13.1. We can find a finite subset, $F$, of $S(A)$ which 
is $\e$-dense in $S(A)$ for both $\rho_{L_1}$ and $\rho_{L_2}$. We
let $K=\mbox{co}(F)$, and we let $B$ consist of the restrictions to $K$
of the elements of $A$. Thus $B$ is finite-dimensional. We
let $L_B^1$ and $L_B^2$ be the quotient Lip-norms on $B$ from $L_1$
and $L_2$. Then
\[
\dist_q((A,L_j),(B,L_B^j)) \le \e
\]
for $j = 1, 2$. Since $\e$ is arbitrary, it follows from the
triangle inequality that it suffices to show that
$\dist_q((B,L_B^1),(B,L_B^2)) \le \d$. Let $\pi$ be the projection 
of $A$ onto $B$. Then from Proposition 3.1 we know that $\pi'$ is 
an isometry for both $((L_B^1)', L_1')$ and $((L_B^2)', L_2')$. Thus
for any $\l \in B'{}^{\circ}$ we have
\[
|(L_B^1)'(\l) - (L_B^2)'(\l) | = |L_1'(\pi'(\l)) - L_2'(\pi'(\l))|
\le \d\|\pi'(\l)\|' \le \d\|\l\|'.
\]
Thus the hypothesis of our proposition is satisfied, and we see
that we have reduced matters to proving the proposition for $A$ 
finite-dimensional.

We now assume that $A$ is finite-dimensional, and
we attempt to define a bridge, $N$, on $A \oplus A$ by
\[
N(a,b) = \d^{-1}\|a-b\|.
\]
Accordingly,
let $L$ be defined as earlier by
\[
L(a,b) = L_1(a) \vee L_2(b) \vee N(a,b).
\]
Because our hypotheses are in terms of the dual seminorms $L'_j$, we must
examine $L'$.  Now $L'$ should be defined on $(A \oplus A)'{}^{\circ}$.
Notice that $A'{}^{\circ} \oplus A'{}^{\circ}$ is of codimension~$1$ in $(A
\oplus A)'{}^{\circ}$, and that for any $\xi \in A'$ we have $(\xi,-\xi) \in
(A \oplus A)'{}^{\circ}$.  Thus any element, $(\zeta_1,\zeta_2)$, of $(A
\oplus A)'{}^{\circ}$ can be expressed (in many ways) as $(\l_1,\l_2) +
(\xi,-\xi)$ for $\l_1,\l_2 \in A'{}^{\circ}$ and $\xi \in A'$.  We claim that
\[
L'(\zeta_1,\zeta_2) = \inf\{L'_1(\l_1) + L'_2(\l_2) + \d\|\xi\|'\},
\]
where the inf is taken over all such expressions of $(\zeta_1,\zeta_2)$.  To
see this, consider first the seminorm defined on $(A \oplus A) \oplus (A
\oplus A)$ by
\[
((L_1 \vee L_2) \vee N)(a,b,c,d) = L_1(a) \vee L_2(b) \vee N(c,d).
\]
The ``dual seminorm'' to $L_1 \vee L_2$ on $A \oplus A$ has value $+\infty$
off of the annihilator of the null-space of $L_1 \vee L_2$, and that
null-space is spanned by $(e,0)$ and $(0,e)$.  Thus $(L_1 \vee L_2)'$ is
finite exactly on $A'{}^{\circ} \oplus A'{}^{\circ}$, and there it is $L'_1
+ L'_2$.  The ``dual seminorm'' of $N$ has value $+\infty$ off of the
annihilator of the null-space of $N$, and that null-space is $\{(a,a): a \in
A\}$.  The annihilator is $\{(\xi,-\xi): \xi \in A'\}$, and on this
annihilator we have $N'(\xi,-\xi) = \d\|\xi\|'$.  Thus the dual of $(L_1
\vee L_2) \vee N$ is defined on $\{(\l_1,\l_2,\xi,-\xi)\}$, and
\[
((L_1 \vee L_1) \vee N)'(\l_1,\l_2,\xi,-\xi) = L'_1(\l_1) + L'_2(\l_2) +
\d\|\xi\|'.
\]
But $A \oplus A$ can be viewed as the subspace of $(a,b,a,b)$'s in $A^4$,
and our $L$ defined above is just the restriction of $(L_1 \vee L_1) \vee N$
to this subspace.  Because we are in the finite-dimensional
situation so that $L$ is continuous, it follows that $L'$ is the 
quotient of the above seminorm $((L_1\vee L_2) \vee N)'$.  (I thank 
Hanfeng Li for pointing out to me that this is difficult to
justify in the infinite-dimensional case.) 
But this gives exactly the formula for $L'$ given
above.

We must check that $L$ induces $L_1$ and $L_2$, and we wish to do this by
applying Corollary $3.10$. Because $A$ is finite-dimensional, $L$, $L_1$
and $L_2$ are all closed. Then
according to Corollary $3.10$ it suffices to show that the restriction
of $L'$ to $A'_1{}^{\circ} \subset (A_1 \oplus A_2)'{}^{\circ}$ coincides
with $L'_1$, and similarly for $L'_2$.  Let $\l \in A'{}^{\circ}$, and
suppose that we have an expression for $(\l,0)$ as
\[
(\l,0) = (\l_1,\l_2) + (\xi,-\xi)
\]
as above.  Then $\xi = \l_2$ and $\l = \l_1 + \l_2$, and so $\|\xi\|'$
becomes $\|\l_2\|'$.  But by hypothesis
\[
|L'_1(\l_2) - L'_2(\l_2)| \le \d\|\l_2\|'.
\]
Thus
\begin{eqnarray*}
L'_1(\lambda) &= &L'_1(\lambda_1 + \lambda_2) \le L'_1(\lambda_1) +
L'_1(\lambda_2) \\
&\le &L'_1(\lambda_1) + L'_2(\lambda_2) + \d\|\xi\|',
\end{eqnarray*}
and so $L'_1(\l) \le L'(\l,0)$.  But we can always take the decomposition
with $\l_1 = \l$ and $\l_2 = 0 = \xi$.  It follows that $L'_1(\l) =
L'(\l,0)$.  In the same way we see that $L'_2(\l) = L'(0,\l)$.  It now
follows that $N$ is a bridge and that $L$ is a
Lip-norm.

Suppose now that $\mu \in S(A)$.  Then
\[
\rho_L((\mu,0) - (0,\mu)) = L'(\mu,-\mu).
\]
In the formula above for $L'$ we can take $\l_1 = 0 = \l_2$ and $\xi =
\mu$.  We thus find that
\[
L'(\mu,-\mu) \le \d\|\mu\|' = \d.
\]
  From this it is clear that
\[
\dist_H^{\rho_L}(S_1,S_2) \le \d,
\]
so that $\dist_q((A,L_1),(A,L_2)) \le \d$.\qed

\bigskip
When we combine this with the earlier considerations, we see that we have
completed the proof of Theorem $13.5$.

We can use some of the facts accumulated above to prove:

\bigskip
\noindent
{\bf 13.15 THEOREM.} {\em The space ${\mathcal Q}$ of isometry equivalence
classes of compact quantum metric spaces, with the metric $\dist_q$, is
separable.}

\bigskip
\noindent
{\it\bfseries Proof.} From Theorem $13.1$ we see that the subset of
${\mathcal Q}$ consisting of finite dimensional spaces is dense.  It thus
suffices to show that for each integer $n$ the set of spaces of dimension
$n$ is separable.  But then it suffices to show that for each integer $D$
the set of spaces of dimension $n$ and diameter $\le D$ is separable.  But
this follows from Proposition $13.7$ by a scaling argument.\qed

\bigskip

As with Corollaries 7.9 and 7.10, our notation for the next corollary
will not destinguish between $C(X)$ and the domain of a Lip-norm
on it.

\bigskip
\noindent
{\bf 13.16 THEOREM.} {\em The set of (equivalence classes of) ordinary
compact metric spaces is a closed subset of $\mathcal Q$. That is, if
$\{(X_n, \rho_n)\}$ is a sequence of ordinary compact metric spaces,
with corresponding Lip-norms $L_{\rho_n}$, and if $(A, L)$ is 
a compact quantum metric space to which the sequence $(C(X_n), L_{\rho_n})$
converges for $\dist_q$, then there is an ordinary compact space $Y$ such
that the completion of $A$ is order-isomorphic to $C(Y)$}

\bigskip
\noindent
{\it\bfseries Proof.} From Lemma 13.6 we see that there is a constant
$D$ such that $\diam(C(X_n), L_{\rho_n}) \le D$ for all $n$. 
Since $X_n$ is identified with the extreme points of $S(C(X_n))$, it
follows that $\diam(X_n, \rho_n) \le D$ for all $n$. From Theorem 13.5
there is a function $G$ such that $\Scv_{C(X_n)}(\e) \le G(\e)$
for each $n$. A simple argument shows that when we view $X_n$ as 
the subset of extreme points of $S(C(X_n))$, 
we have $\Cov_{\rho_n}(\e) \le \Scv_{C(X_n)}(\e/2)$. From Gromov's 
compactness and completeness theorems it follows that a subsequence 
of the $C(X_n)$'s converges for $\dist_{GH}$ to some ordinary compact metric
space, say $(Y, \rho)$. From Proposition 4.7 it follows that as 
quantum metric spaces this subsequence converges to $(Y, \rho)$
also for $\dist_q$. Thus
\[
\dist_q((A, L_A), ((C(Y), L_\rho)) = 0 .
\]
Then the completion of $A$ is isomorphic to $C(Y)$ as order-unit spaces, 
and under this isomorphism we have $L_A = L_\rho$ by Theorem 7.7.   \qed

Finally, we give more quantitative relations between 
the three measures of growth defined in Definition $13.4$.  We begin with:

\bigskip
\noindent
{\bf 13.17 PROPOSITION.} {\em Let $(A,L)$ be a compact quantum metric space
of dimension $n$, and let $D = \diam(A,L)$.  Then for any $\e > 0$, there is
a subset of $S(A)$ which is $\e$-dense in $S(A)$ for $\rho_L$ and contains
no more than $((D/\e) + 1)^{n-1}$ points.}

\bigskip
\noindent
{\it\bfseries Proof.} Assume first that $D = 1$, so that $L' \le
(1/2)\|\cdot\|'$.  Fix $\mu_0 \in S(A)$ and set $K = S(A) - \mu_0$, a subset
of $A'{}^{\circ}$.  Let $B$ denote the open unit $L'$-ball about $0$ in
$A'{}^{\circ}$, and ${\bar B}$ its closure.  Thus $K \subseteq {\bar B}$
because $D = 1$.  We now use a standard argument.  (See the proof of lemma
$4.10$ of \cite{Pi}.  I am indebted to Bernd Sturmfels for steering me
toward this argument.)  Let $\{\l_1,\dots,\l_N\}$ be a subset of $K$ such
that $\|\l_j - \l_k\|' \ge 2\e$ for all $j\ne k$.  Then the open balls $\l_j
+ \e B$ are disjoint from each other, and are all contained in $(1 +
\e)B$.  Since $A'{}^{\circ}$ is finite-dimensional, there is a
translation-invariant volume on $A'{}^{\circ}$.  When we apply it to the
above situation, we see that
\[
N \mbox{ vol}(\e B) \le \mbox{vol}((1+\e)B).
\]
Since $A'{}^{\circ}$ has dimension $n-1$, it follows that
\[
N \le ((1/\e) + 1)^{n-1}.
\]
Now let $F$ be a maximal subset $K$ with respect to the above property that
$\|\l_j-\l_k\|' \le 2\e$.  Then $F$ is $2\e$-dense in $K$ for $\|\cdot\|'$,
and so
$\e$-dense in
$K$ for
$L'$.  It follows that $F + \mu_0$ is $\e$-dense in $S(A)$ for $\rho_L$.
Thus we obtain the desired conclusion when $D = 1$.

If $D \ne 1$, then $D^{-1}L'$ has diameter $1$, and so there is a subset $F$
of $S(A)$ which is $(\e/D)$-dense for $D^{-1}L'$, and has no more than
$((D/\e) - 1)^{n-1}$ elements.  Then $DF$ is $\e$-dense for $L'$.\qed

\bigskip
The important aspect for us of the bound given in Proposition $13.16$ is
that it depends only on $n$ and $D$, and not on other features of $A$ or $L$.

Suppose now that $(A,L_A)$ is a general compact quantum metric space, and
let $\e > 0$ be given.  Let $D = \diam(A,L_A)$, and let $n= \Fin_L(\e)$.
Thus we can find a compact quantum metric space $(B,L_B)$ such that
$\dist_q(A,B) < \e$ and $\dim(B)= n$.  Thus there is a Lip-norm $L$ on $A
\oplus B$ inducing $L_A$ and $L_B$ such that $\dist_H^{\rho_L}(S(A),S(B))
<\e$.  It is easily seen that $\diam(B,L_B) \le D + 2\e$.  Thus by
Proposition $13.16$ we can find a finite subset, $F$, of $S(B)$ which is
$\e$-dense in $S(B)$, and for which
\[
|F| \le (((D+2\e)/\e) + 1)^{n-1},
\]
where $|\cdot|$ denotes ``number of elements in''.  For each point in $F$
choose a point in $S(A)$ within distance $\e$ of it, and let $F_A$ denote
the set of these points.  It is easily seen that $F_A$ is $3\e$-dense in
$S(A)$.  Thus $\Scv_{L_A}(3\e) \le |F_A| \le |F|$.  Upon simplifying the
earlier bound for $|F|$, we obtain:

\bigskip
\noindent
{\bf 13.18 PROPOSITION.} {\em Let $(A,L)$ be a compact quantum metric space
of diameter $D$.  Then for every $\e > 0$ we have}
\[
\Scv_L(3\e) \le ((D/\e) + 5)^{\Fin_L(\e)-1}.
\]

\bigskip
  From this proposition we see that if ${\mathcal S}$ is a subset of
${\mathcal Q}$ for which there is a constant $D$ and function $G$ such that
every element of ${\mathcal S}$ is of diameter $\le D$, and that $\Fin_L(\e)
\le G(\e)$ for all $(A,L) \in {\mathcal S}$, and $\e > 0$, then there is a
function $H$ from ${\mathbb R}^+$ to ${\mathbb N}^+$ such that $\Scv_L(\e)
\le H(\e)$ for all $(A,L) \in {\mathcal S}$ and all $\e > 0$.  This is the
equivalence for the purposes of the compactness theorem which we had in mind.

\section*{Appendix 1. An example where $\dist_{GH}>\dist_q$ }
\noindent
by Hanfeng Li

\bigskip
Let $Y=\{y_1, y_2, y_3\}$, with metric $\rho_Y(y_1, y_2)=
1 = \rho_Y(y_2, y_3)$ and $ \rho_Y(y_1, y_3)=2$. 
Let $Z=\{z_1, z_2\}$ with metric $\rho_Z(z_1, z_2)=3$. 
Let $C(Y)$ and $C(Z)$ be the algebras of real-valued functions 
on $Y$ and $
Z$, with the Lip-norms $L_{\rho_Y}$ and $L_{\rho_Z}$. 
We will show that $\dist_{GH}(Y, Z)=1$ but that
$\dist_q((C(Y), L_{\rho_Y}), (C(Z), L_{\rho_Z}))= 1/2$. 
Hence this is an example for which $\dist_{GH}>\dist_q$. 

We first show that $\dist_{GH}(Y, Z)=1$. Let $\rho$ be a metric on the
union of $Y$ and $Z$ which restricts to the given metrics. Consider
the distance from $y_2$ to $Z$. Suppose
that $\rho(y_2, z_1) \leq 1$. Since for all $j$ we have
$$
3 = \rho(z_1, z_2) \leq \rho(z_1, y_2) + \rho(y_2, y_j) + \rho(y_j, z_2) ,
$$
we must then have $\rho(y_j, z_2) \geq 1$, that is, 
$\rho(z_2, Y) \geq 1$. In the same way, if 
$\rho(y_2, z_2) \leq 1$, then $\rho(z_1, Y) \geq 1$. Thus
$\dist_{GH}(Z, Y) \geq 1$. But we can isometrically embed $Y$ 
and $Z$ in $\mathbb R$ as $\{0, 1, 2\}$ and $\{0, 3\}$, from
which we see that $\dist_{GH}(Y, Z)=1$.

To show that $\dist_q((C(Y), L_{\rho_Y}), (C(Z), L_{\rho_Z}))= 1/2$, 
we need some preparation. 
Notice first that Proposition 5.7 says in effect that for compact 
quantum metric spaces $(B_i, L_i)$, if their state spaces 
are affinely isometrically embedded
into the state-space $S(A)$ of some other 
compact quantum metric space $(A, L)$, then
$$
 \dist_q(B_1, B_2)\le \dist^{\rho_L}_H(S(B_1), S(B_2)) .
$$
This provides a powerful way of getting upper bounds for quantum 
Gromov-Hausdorff distance. 

Here we will use the special case in which we are given $(A, L)$ and we just
enlarge the state space within $A'$. Let $H(A)$ be the
hyperplane $\{\eta \in A':\eta(e_A)=1\}$. 
Then $H(A)$ is a convex set containing $S(A)$.  
It is clear that $H(A)-H(A) \subseteq A'{}^{\circ}$. Thus the usual 
formula $\rho_L(\mu, \nu) = L'(\mu - \nu)$ actually 
defines $\rho_L$ as a metric on all of $H(A)$. It is easy to see 
that the proof of theorem 1.8 of \cite{R4} actually shows that 
the $\rho_L$-topology coincides with the $w^*$-topology on 
any $\|\cdot\|'$-bounded subset of $H(A)$.

Let $K$ be a closed convex  $\|\cdot\|'$-bounded (so $w^*$-compact)
subset of $H(A)$ which contains $S(A)$, 
and let $\pi$ denote the evident restriction map from $A$ into $Af(K)$, 
much as in Proposition 3.5. Let $B=\pi(A)$. 
Evidently $B$ is an order-unit space. From corollary I.1.5 
of \cite{A} it is easily seen that $B$ is dense in $Af(K)$ 
for the supremum norm. Note that $\pi$ is usually not a morphism, 
but $\pi$ is a bijection from $A$ onto $B$ because $K \supseteq S(A)$; 
and it is $\pi^{-1}$ which is a morphism. In particular, $B$ is just 
the vector space $A$ with the same distinguished element, 
but different order structure and norm. However, we will use 
the same $L$ on both $A$ and $B$.

Let $\|\cdot\|_B$ denote the supremum norm on $Af(K)$, and so on $B$. 
Then $\|\cdot\| \leq \|\cdot\|_B$ because $S(A) \subseteq K$. 
Since $K$ is  $\|\cdot\|'$-bounded, there is an $s \in \mathbb R$ 
such that $\|\cdot\|_B \leq s\|\cdot\|$. In particular, 
every state on $B$ will be continuous for $\|\cdot\|$, 
and so will be in $A'$. A simple argument quite similar to 
that in the proof of Proposition 3.5 now shows that $S(B)$ is 
naturally identified with  $K$. Since $K \supseteq S(A)$, 
we are then in the situation of Proposition 3.5 but with the 
roles of $A$ and $B$ reversed. For any $\mu, \nu \in K$ and 
any $a \in A$ we have
$$
|\mu(\pi(a))-\nu(\pi(a))| = |(\mu - \nu)(a)| \leq 
L'(\mu - \nu)L(a) = \rho_L(\mu,\nu)L(a).
$$
Thus each element of $B$ is Lipschitz for $\rho_L$, 
and $L_{\rho_L}(a) = L(a)$ as long as $L$ is lower semi-continuous.

Suppose now that $K_1$ is another compact convex subset of $K$, 
and let $(C, L_C)$ be the corresponding quotiont of $(B, L)$. 
Then from Proposition 5.7 we know that
$$
\dist_q(A, C) \leq \dist_H^{\rho_L}(S(A), K_1).
$$
We apply this fact to the example introduced at the beginning 
of this appendix.
Let $A = C(Y)$. Let $w_1$ and $w_2$ be the points of $H(A)$ defined by
$$
w_1 = y_1 + (y_2 - y_3)/2 , \quad
 \quad w_2 = y_3 +(y_2 - y_1)/2  ,
$$
for the evident meaning of the notation. Let $K \subset A'{}^{\circ}$
be defined by
$$
K = \mbox{co}\{y_1, y_2, y_3, w_1, w_2\} .
$$ 
Thus $K \supset S(A)$, and the above discussion applies. 
We let $B$ denote $A$ but with the order and norm coming 
from viewing its elements as affine functions on $K$. 
Thus for $f \in B=C(Y)$ we have $f \geq 0$ exactly 
if $f \geq 0$ on $Y$ and $f(w_j) \geq 0$ for $j = 1, 2$. 
We equip $B$ with the Lip-norm $L$. Thus $\rho_L$ is defined 
on $K$ by $L'$, and its restriction to $S(A)$ is the original 
metric coming from that on $Y$.

Now it is easy to see that for $f = \{f_1, f_2, f_3\} \in C(Y)$ we have
$$
L(f_1, f_2, f_3) = |f_1 - f_2| \vee |f_2 - f_3|.
$$
It follows that for $\lambda = 
\{\lambda_1, \lambda_2, \lambda_3\} \in  A'{}^{\circ}$ we have
$$
L'(\lambda_1, \lambda_2, \lambda_3) = |\lambda_1| + |\lambda_3| .
$$
Then we calculate that
$$
\rho_L(w_1, w_2) = L'((3/2)(y_1 - y_3)) = 3 .
$$
Consequently we can identify $Z$ isometrically with a subset of $K$ by 
sending $z_j$ to $w_j$ for $j = 1, 2$. We let $K_1$ be the closed 
line-segment in K joining $w_1$ and $w_2$. Then $C(Z)$ becomes 
identified with the space $C$ of affine functions on $K_1$, 
and $S(C) = K_1$. Because $C$ is 2-dimensional, there is only 
one possible Lip-norm on $C$ giving distance 3 between $w_1$ 
and $w_2$, and so it must be the quotient of $L$ on $C$. 
We are now in position to apply Proposition 5.7. From it we 
conclude that
$$
\dist_q(A, C) \leq \dist_H^{\rho_L}(S(A), K_1) .
$$
Now
$$
\rho_L(y_1, w_1) = L'((y_3 - y_2)/2) = 1/2 ,
$$
and similarly $\rho_L(y_3, w_2) = 1/2$. Furthermore, 
$(w_1 + w_2)/2 \in K_1$, and
$$
\rho_L(y_2, (w_1+w_2)/2) = 
L'((-y_1 + 2y_2 - y_3)/4) = 1/2 .
$$
Because $\rho_L$ is convex, it follows that 
$\dist_H^{\rho_L}(S(A), K_1) \leq 1/2$ .
Consequently $\dist_q(C(Y), C(Z)) \leq 1/2$. But from the fact 
that the diameter of $Y$ is 2 while the diameter of $Z$ is 3 
it follows easily that $\dist_q(C(Y), C(Z)) \geq 1/2$, so 
that  $\dist_q(C(Y), C(Z)) = 1/2$ as desired.

\section*{Appendix 2.  Dirac operators are universal}

In this brief section we answer question $11.1$ of \cite{R5}.  Namely, we
show that, in a suitable sense, every lower semi-continuous Lip-norm on an
order-unit space can be obtained from a ``Dirac'' operator.  This is seen as
follows.

Let $A$ be an order-unit space, and let $L$ be a lower semi-continuous
Lip-norm on $A$.  Let $\rho_L$ denote the corresponding metric on $S(A)$.
According to theorem $4.2$ of \cite{R5}, because $L$ is lower
semi-continuous, we can recover $L$ from $\rho_A$ by the formula
\[
L(a) = \sup\{|\mu(a) - \nu(a)|/\rho_L(\mu,\nu): \mu \ne \nu\}.
\]
(Notice that the right-hand side is always lower semi-continuous.)

Let $C(S(A))$ denote the algebra of complex-valued continuous functions on
the compact space $S(A)$, and let ${\mathcal L}$ denote the dense
$*$-subalgebra of $C(S(A))$ consisting of the Lipschitz functions for
$\rho_L$.  Let $L_0$ denote the ordinary Lipschitz seminorm on ${\mathcal
L}$ for $\rho_L$.  We view $A$ as a real subspace of $C(S(A))$ in the usual
way.  Then the formula above for $L$ shows that
$L$ is the restriction of $L_0$ to $A$.

Now $S(A)$, as a compact metric space, is separable, and so we can find (in
many ways) a positive finite Radon measure, $m$, on $S(A)$ whose support is
all of $S(A)$.  We now use the construction described late in section~11 of
\cite{R5}, or on page~274 of \cite{W2}.  We let $\D$ denote the diagonal of
$S(A) \x S(A)$ and set $Y = (S(A) \x S(A)) \smallsetminus \D$.
We restrict $m\x
m$ to $Y$, and let ${\mathcal H} = L^2(Y,m \x m)$.  We represent $C(S(A))$
on ${\mathcal H}$ by
\[
(f\xi)(\mu,\nu) = f(\mu)\xi(\mu,\nu)
\]
for $f \in C(S(A))$ and $\xi \in {\mathcal H}$.  We let $D$ denote the
(usually unbounded) operator on ${\mathcal H}$ defined by
\[
(D\xi)(\mu,\nu) = \xi(\nu,\mu)/\rho_L(\mu,\nu),
\]
with domain those $\xi$'s for which $D\xi \in {\mathcal H}$.  It is easily
seen that $D$ is self-adjoint.  Furthermore, as seen by simple calculations
(given in the references above), for any $f \in {\mathcal L}$, viewed as an
operator on ${\mathcal H}$, we find that $[D,f]$ is a bounded operator, and
that $L_0(f) = \|[D,f]\|$.

In particular, we see that $A$ is represented isomorphically as an
order-unit space of self-adjoint operators on ${\mathcal H}$.  Since $A
\subseteq {\mathcal L}$, each $[D,a]$ is a bounded operator on ${\mathcal
H}$, and $L(a) = L_0(a) = \|[D,a]\|$.  In this way $L$ is obtained from the
``Dirac'' operator $D$.

\end{document}